\documentclass[12pt, reqno, english]{amsart}
\usepackage{amsmath, amsthm, amssymb, color, xcolor}
\usepackage[colorlinks=true,citecolor=red,linkcolor=blue,urlcolor=blue]{hyperref}
\usepackage{graphicx}
\usepackage{comment}
\usepackage{caption}
\usepackage{bold-extra}
\usepackage{mathtools}
\usepackage{enumerate}
\usepackage{bm}
\usepackage{rotating}
\usepackage{mathrsfs}
\usepackage{verbatim}
\usepackage{tikz, tikz-cd, tikz-3dplot}
\usepackage{amssymb}
\usepackage{secdot}
\usepackage{algorithm}
\usepackage[noend]{algpseudocode}
\usepackage{caption}
\usepackage[normalem]{ulem}
\usepackage{subcaption}
\usepackage{multicol}
\usepackage{makecell}
\usepackage{array}
\usepackage{enumitem}
\usepackage{adjustbox}
\usepackage{blkarray}
\usepackage[top=25mm, bottom=25mm, left=25mm, right = 25mm]{geometry}
\usepackage{cleveref}
\usepackage{lineno}
\usepackage{enumitem}
\usepackage{titlesec}
\usetikzlibrary{matrix}
\usetikzlibrary{arrows}
\usetikzlibrary{decorations.pathmorphing}
\usetikzlibrary{patterns}

\titleformat{\section}
       {\centering \fontsize{12}{17} \large \bf \scshape }{\thesection}{0mm}{. \hspace{0.00mm}}

\titleformat{\subsection}[runin]
       {\fontsize{12}{17} \bf}{\thesubsection}{0mm}{. \hspace{0.00mm}}[.\\]

\newtheorem{theorem}{Theorem}[section]

\newtheorem{proposition}[theorem]{Proposition}
\newtheorem{lemma}[theorem]{Lemma}

\theoremstyle{definition}
\newtheorem{definition}[theorem]{Definition}

\newtheorem{remark}[theorem]{Remark}
\newtheorem{cor}[theorem]{Corollary}
\newtheorem{example}[theorem]{Example}

\newcommand{\Pf}{\mathrm{Pf}}

\newcommand{\PP}{\mathbb{P}}
\newcommand{\RR}{\mathbb{R}}

\newcommand{\CC}{\mathbb{C}}
\newcommand{\ZZ}{\mathbb{Z}}
\newcommand{\NN}{\mathbb{N}}
\renewcommand{\SS}{\mathbb{S}}

\newcommand{\Id}{\mathrm{Id}}

\newcommand{\Gr}{\mathrm{Gr}}
\newcommand{\OGr}{\mathrm{OGr}}

\newcommand{\Rcal}{\mathcal{R}}

\newcommand{\Span}{\mathrm{span}}

\newcommand{\SO}{\mathrm{SO}}
\newcommand{\Spin}{\mathrm{Spin}}
\newcommand{\SL}{\mathrm{SL}}

\newcommand{\LGP}{\mathrm{LGP}}

\newcommand{\rowspan}{\mathrm{rowspan}}

\renewcommand{\mod}{\mathrm{\ mod \ }}

\definecolor{dgreen}{HTML}{026a10}
\definecolor{dviolet}{HTML}{9109E3}
\definecolor{dorange}{HTML}{e55700}

\DeclareMathOperator{\sgn}{sgn}

\renewcommand{\tilde}{\widetilde}
\usepackage{nicematrix}

\title{\bf Totally positive skew-symmetric matrices}
\author[J. Boretsky]{Jonathan Boretsky}
\address{Jonathan Boretsky (MPI MiS)}
\email{jonathan.boretsky@mis.mpg.de}

\author[V. Calvo Cortes]{Veronica Calvo Cortes}
\address{Veronica Calvo Cortes (MPI MiS)}
\email{veronica.calvo@mis.mpg.de}

\author[Y. El Maazouz]{Yassine El Maazouz}
\address{Yassine El Maazouz (Caltech)}
\email{maazouz@caltech.edu}

\date{\today}


\begin{document}

\begin{abstract}
    A matrix is totally positive if all of its minors are positive. This notion of positivity coincides with the type A version of Lusztig's more general total positivity in reductive real-split algebraic groups. Since skew-symmetric matrices always have nonpositive entries, they are not totally positive in the classical sense.  The space of skew-symmetric matrices is an affine chart of the orthogonal Grassmannian \texorpdfstring{$\OGr(n,2n)$}{OGr(n,2n)}. Thus, we define a skew-symmetric matrix to be \emph{totally positive} if it lies in the \emph{totally positive orthogonal Grassmannian}. We provide a positivity criterion for these matrices in terms of a fixed collection of minors, and show that their Pfaffians have a remarkable sign pattern. The totally positive orthogonal Grassmannian is a CW cell complex and is subdivided into \emph{Richardson cells}. We introduce a method to determine which cell a given point belongs to in terms of its associated matroid.
\end{abstract}

\maketitle

\setcounter{tocdepth}{1}

\section{Introduction}

Let $n \geq 1$ be a positive integer and denote by $\SS_n = \SS_n(\RR)$ the $\binom{n}{2}$-dimensional real vector space of skew-symmetric $n \times n$ matrices  with real entries. This article studies the semi-algebraic set $\SS_n^{> 0}\subset \SS_n$ of \emph{totally positive} skew-symmetric matrices. The latter are defined using total positivity of partial flag varieties in the sense of Lusztig \cite{Lusztig1}, as follows.

\smallskip
Let $q$ be the non-degenerate symmetric bilinear form on $\RR^{2n}$ given by
\begin{equation}\label{eq:quadForm}
    q(x,y) = \sum_{i=1}^{n} x_{i} y_{n+i}  +  \sum_{i=1}^{n} y_{i} x_{n+i}, \quad \text{for } x,y \in \RR^{2n}.
\end{equation}

In the standard basis $(e_{1}, \dots, e_{n}, f_{1}, \dots, f_{n})$ of $\RR^{2n}$, this bilinear form is given by the matrix 
\[
    Q = \begin{bmatrix} 0 & \Id_n \\ \Id_n & 0 \end{bmatrix},
\]
where $\Id_n$ is the $n \times n$ identity matrix. The \emph{orthogonal Grassmannian} is the variety of $n$-dimensional vector subspaces $V$ of $\RR^{2n}$ that are \emph{$q$-isotropic}, meaning that $q(v,w) = 0$ for any $v,w \in V$. Two distinguished points in this variety are the vector spaces 
\begin{equation}\label{eq:EFspaces}
    E := \Span( e_1, \dots, e_n) \quad \text{and} \quad F:= \Span(f_{1}, \dots, f_{n}).
\end{equation}
The orthogonal Grassmannian is a smooth algebraic variety embedded in $\mathbb{RP}^{\binom{2n}{n}-1}$ by Pl\"ucker coordinates. It has two isomorphic irreducible connected components of dimension $\binom{n}{2}$:
\begin{align*}
    \OGr(n,2n) :=& \{ V \text{ $q$-isotropic} \colon \dim(V)=n,\; \dim(E \cap V) =  n \mod 2\},\\
    \OGr_{-}(n,2n) :=& \{ V \text{ $q$-isotropic} \colon \dim(V)=n,\; \dim(E \cap V) =  n+1 \mod 2 \}.
\end{align*}

The Zariski open set in $\OGr(n,2n)$ where the Pl\"ucker coordinate $\Delta^{1, \dots, n }$ does not vanish is isomorphic to the affine space $\SS_n$. This isomorphism identifies $A \in \SS_n$ with the rowspan of the $n\times 2n$ matrix $\begin{bmatrix} \Id_n | A \end{bmatrix}$.  

We may also view $\OGr(n,2n)$ as the connected component of the identity in a parabolic quotient of the real special orthogonal group ${\rm SO}(2n)$. This is a connected reductive $\RR$-split algebraic group and therefore admits a \emph{totally positive part} in the sense of Lusztig \cite{Lusztig1}.

\smallskip

A key example of Lusztig positivity is the case of ${\rm SL}(n)$. A parabolic quotient of ${\rm SL}(n)$ is a \textit{flag variety} whose points are flags of linear subspaces. Such flags can be represented as row spans of matrices in ${\rm SL}(n)$. Lusztig's total positivity then matches the classical notion of total positivity: a flag is totally positive (resp. nonnegative) if it can be represented by a totally positive (resp. nonnegative) matrix, that is, one whose minors are all positive (resp. nonnegative). In general, the totally nonnegative part of a flag variety admits a nice topological structure that interplays well with matroid theory \cite{PositiveGeometries,GKL_Ball22,postnikov06}. These notions have become increasingly important to understand for other real reductive groups as positivity and positroid combinatorics are gaining more relevance in the study of scattering amplitudes in quantum field theory \cite{ABCGPT, TheAmplituhedron,WilliamsICM}.

\begin{definition}
    A skew-symmetric matrix $A \in \SS_n$ is \textit{totally nonnegative} (resp. \textit{totally positive}) if the rowspan of $\begin{bmatrix} \Id_n | A \end{bmatrix}$ is a point in the totally nonnegative region $\OGr^{\geq 0}(n,2n)$ (resp. the totally positive region $\OGr^{> 0}(n,2n)$) of $\OGr(n,2n)$. See \Cref{def:LusztigPositive} for more details.
\end{definition}

Given a skew-symmetric matrix, or more generally, a point in any partial flag variety, it is difficult to determine directly from the definition whether it is totally positive. Accordingly, positivity tests for certain partial flag varieties have been developed, for example \cite{BFZIII, ChevalierPositivity}. However, these positivity criteria are sometimes not very explicit. Explicit tests for positivity have been described in type A \cite{BlochKarp,BossingerLi} and for certain flag varieties of types B and C \cite{BBEG24}. In this article we give an explicit and minimal positivity test for a skew symmetric matrix $A$ in terms of its minors, which mirrors the fact that total positivity on $\SL(n)$ is determined by the positivity of minors.

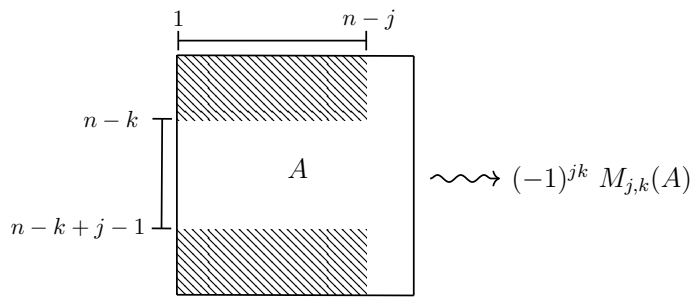
\begin{figure}[ht]
            \centering
            \scalebox{0.881}{
                \tikzset{every picture/.style={line width=0.75pt}}
                \begin{tikzpicture}[x=0.75pt,y=0.75pt,yscale=-1,xscale=1]
                
                \draw    (135,44.4) -- (135,181.6) ;
                \draw    (270.6,44.8) -- (270.6,182) ;
                
                \draw    (135,44.4) --  (270.6,44.8)  ;
                \draw    (135,181.6)-- (270.6,182) ;

                \draw    (135.8,36) -- (243.4,36) ;
                \draw [shift={(243.4,36)}, rotate = 180] [color={rgb, 255:red, 0; green, 0; blue, 0 }  ][line width=0.75]    (0,5.59) -- (0,-5.59)   ;
                \draw [shift={(135.8,36)}, rotate = 180] [color={rgb, 255:red, 0; green, 0; blue, 0 }  ][line width=0.75]    (0,5.59) -- (0,-5.59)   ;
                \draw    (126.6,80.8) -- (126.2,143.6) ;
                \draw [shift={(126.2,143.6)}, rotate = 270.36] [color={rgb, 255:red, 0; green, 0; blue, 0 }  ][line width=0.75]    (0,5.59) -- (0,-5.59)   ;
                \draw [shift={(126.6,80.8)}, rotate = 270.36] [color={rgb, 255:red, 0; green, 0; blue, 0 }  ][line width=0.75]    (0,5.59) -- (0,-5.59)   ;
                
                \path[pattern color=black, pattern=north west lines] (135,44.4) -- (243.4,44.4) -- (243.4,81.6) -- (135,81.6) -- cycle ;
                \path[pattern color=black, pattern=north west lines] (135,144.4) -- (243.4,144.4) -- (243.4,181.6) -- (135,181.6) -- cycle ;
                
                \draw (228.4,16.6) node [anchor=north west][inner sep=0.75pt]  [font=\footnotesize]  {$n-j$};
                \draw (131.2,17.4) node [anchor=north west][inner sep=0.75pt]  [font=\footnotesize]  {$1$};
                \draw (80,74.6) node [anchor=north west][inner sep=0.75pt]  [font=\footnotesize]  {$n-k$};
                \draw (38.4,136.2) node [anchor=north west][inner sep=0.75pt]  [font=\footnotesize]  {$n-k+j-1$};
                \path[draw, ->, decorate, decoration ={snake, amplitude = 1.5}] (280,115) -- (320,115);
                \draw (325,105) node [anchor=north west][inner sep=0.75pt]  [font=\normalsize]  {$( -1)^{jk} \ M_{j,k}(A)$};
                \draw (196.8,102.2) node [anchor=north west][inner sep=0.75pt]  [font=\normalsize]  {$A$};
                \end{tikzpicture}
                }
            \caption{The shading indicates which minor of the matrix $A$ is used to compute $M_{j,k}$.}
            \label{fig:Minor}
        \end{figure}

\begin{definition}\label{def:SpecialMinorsPfaff}
    For any $n \times n$ matrix $A$ we denote by $\Delta_{I}^J(A)$ be the determinant of the submatrix of $A$ in rows $I$ and columns $J$. We denote by $M_{j,k}(A)$ the signed minor:
    \begin{equation}\label{eq:SpecialMinors}
    M_{j,k}(A) = (-1)^{jk} \Delta_{\{1,\ldots,n-k-1,n-k+j, \ldots, n \}}^{\{1,2,\ldots, n-j\}}(A) \qquad\text{for any } 1 \leq j \leq k \leq n-1.    
    \end{equation}
    Note that the minor $M_{j,k}$ is a polynomial of degree $n-j$ in the entries of $A$. It corresponds up to a sign to a left justified minor where the rows are indexed by the complement of an interval, as illustrated by the shaded region in \Cref{fig:Minor}.

\end{definition}

\begin{example}[$n=4$]\label{ex:n=4Minors}
    The minors $M_{j,k}(A)$ for $1 \leq j \leq k \leq 3$ for a $4 \times 4$ skew-symmetric matrix $A=(a_{ij})$ are the following:
    \begin{alignat*}{3}
    M_{1,1}(A) &= a_{12}a_{14}a_{23}-a_{12}a_{13}a_{24}+a_{12}^2a_{34},   & &   \\
    M_{1,2}(A) &= a_{13}^{2}a_{24}-a_{13}a_{14}a_{23}-a_{12}a_{13}a_{34},   \quad M_{2,2}(A) &=&a_{12}a_{14}, \quad  \\
    M_{1,3}(A) &= a_{14}a_{23}^2-a_{13}a_{23}a_{24}+a_{12}a_{23}a_{34},   \quad M_{2,3}(A) &=& a_{13}a_{24} - a_{14}a_{23}, \quad M_{3,3}(A) &= a_{14}.
    \end{alignat*}

\end{example}

We realize these minors via a graphical interpretation of the Marsh-Rietsh parametrization \cite{MR} of $\OGr^{>0}(n,2n)$, using the Lindst\"rom-Gessel-Viennot (LGV) lemma. Our first main result is a positivity test for $\OGr^{>0}(n,2n)$ using the signed minors in \Cref{def:SpecialMinorsPfaff}.

\begin{figure}[H]
        \centering
         \scalebox{0.9}{\begin{tikzpicture}
            \coordinate (l6) at (0,4.5);    
            \coordinate (l7) at (0,4);
            \coordinate (l8) at (0,3.5);
            \coordinate (l9) at (0,3);
            \coordinate (l10) at (0,2.5);        
            \coordinate (l5) at (0,2);        
            \coordinate (l4) at (0,1.5);    
            \coordinate (l3) at (0,1);
            \coordinate (l2) at (0,0.5);
            \coordinate (l1) at (0,0);
        
            \coordinate (r6) at (7,4.5);    
            \coordinate (r7) at (7,4);
            \coordinate (r8) at (7,3.5);
            \coordinate (r9) at (7,3);
            \coordinate (r10) at (7,2.5);
            \coordinate (r5) at (7,2);        
            \coordinate (r4) at (7,1.5);    
            \coordinate (r3) at (7,1);
            \coordinate (r2) at (7,0.5);
            \coordinate (r1) at (7,0);

            \coordinate (v11) at (2.5,0);
            
            \coordinate (v21) at (2,0.5);
            \coordinate (v22) at (2.5,0.5);
            \coordinate (v23) at (4,0.5);
        
            \coordinate (v31) at (1.5,1);
            \coordinate (v32) at (2,1);
            \coordinate (v33) at (3.5,1);
            \coordinate (v34) at (4,1);
            \coordinate (v35) at (5.5,1);
            
            \coordinate (v41) at (1,1.5);
            \coordinate (v42) at (1.5,1.5);
            \coordinate (v43) at (3,1.5);
            \coordinate (v44) at (3.5,1.5);
            \coordinate (v45) at (5,1.5);
            \coordinate (v46) at (5.5,1.5);
            \coordinate (v47) at (6,1.5);
            
            \coordinate (v51) at (0.5,2);
            \coordinate (v52) at (3,2);
            \coordinate (v53) at (4.5,2);
            \coordinate (v54) at (6,2);

            \coordinate (v101) at (1,2.5);    
            \coordinate (v102) at (3,2.5);
            \coordinate (v103) at (5,2.5);
            \coordinate (v104) at (6,2.5);

            \coordinate (v91) at (0.5,3);
            \coordinate (v92) at (1.5,3);
            \coordinate (v93) at (3,3);
            \coordinate (v94) at (3.5,3);
            \coordinate (v95) at (4.5,3);
            \coordinate (v96) at (5.5,3);
            \coordinate (v97) at (6,3);

            \coordinate (v81) at (1.5,3.5);
            \coordinate (v82) at (2,3.5);
            \coordinate (v83) at (3.5,3.5);
            \coordinate (v84) at (4,3.5);
            \coordinate (v85) at (5.5,3.5);

            \coordinate (v71) at (2,4);
            \coordinate (v72) at (2.5,4);
            \coordinate (v73) at (4,4);
        
            \coordinate (v61) at (2.5,4.5);
            
            \draw[black] (-0.5,0) node  [xscale = 0.8, yscale = 0.8] {$1$};
            \draw[black] (-0.5,0.5) node  [xscale = 0.8, yscale = 0.8] {$2$};
            \draw[black] (-0.5,1) node  [xscale = 0.8, yscale = 0.8] {$3$};
            \draw[black] (-0.5,1.5) node  [xscale = 0.8, yscale = 0.8] {$4$};
            \draw[black] (-0.5,2) node  [xscale = 0.8, yscale = 0.8] {$5$};
            \draw[black] (-0.5,2.5) node  [xscale = 0.8, yscale = 0.8] {$10$};
            \draw[black] (-0.5,3) node  [xscale = 0.8, yscale = 0.8] {$9$};
            \draw[black] (-0.5,3.5) node  [xscale = 0.8, yscale = 0.8] {$8$};
            \draw[black] (-0.5,4) node  [xscale = 0.8, yscale = 0.8] {$7$};
            \draw[black] (-0.5,4.5) node  [xscale = 0.8, yscale = 0.8] {$6$};
            
            \draw[black] (7.5,0) node  [xscale = 0.8, yscale = 0.8] {$1$};
            \draw[black] (7.5,0.5) node  [xscale = 0.8, yscale = 0.8] {$2$};
            \draw[black] (7.5,1) node  [xscale = 0.8, yscale = 0.8] {$3$};
            \draw[black] (7.5,1.5) node  [xscale = 0.8, yscale = 0.8] {$4$};
            \draw[black] (7.5,2) node  [xscale = 0.8, yscale = 0.8] {$5$};
            \draw[black] (7.5,2.5) node  [xscale = 0.8, yscale = 0.8] {$10$};
            \draw[black] (7.5,3) node  [xscale = 0.8, yscale = 0.8] {$9$};
            \draw[black] (7.5,3.5) node  [xscale = 0.8, yscale = 0.8] {$8$};
            \draw[black] (7.5,4) node  [xscale = 0.8, yscale = 0.8] {$7$};
            \draw[black] (7.5,4.5) node  [xscale = 0.8, yscale = 0.8] {$6$};

            \draw[draw=dorange, line width=2pt] (l1)  -- (v11);
            \draw[draw=black, line width=1pt] (v11) -- (r1);

            \draw[draw=dviolet, line width=2pt] (l2) -- (v21);
            \draw[draw=black, line width=1pt] (v21) -- (v22);
            \draw[draw=dorange, line width=2pt] (v22) -- (v23);
            \draw[draw=black, line width=1pt] (v23) -- (r2);

            \draw[draw=dgreen, line width=2pt] (l3) -- (v31);
            \draw[draw=black, line width=1pt] (v31) -- (v32);
            \draw[draw=dviolet, line width=2pt] (v32) -- (v33);
            \draw[draw=black, line width=1pt] (v33) -- (v34);
            \draw[draw=dorange, line width=2pt] (v34) -- (r3);

            \draw[draw=blue, line width=2pt] (l4) -- (v41);
            \draw[draw=black, line width=1pt] (v41) -- (v42);
            \draw[draw=dgreen, line width=2pt] (v42) -- (v43);
            \draw[draw=black, line width=1pt] (v43) -- (v44);
            \draw[draw=dviolet, line width=2pt] (v44) -- (r4);
            
            \draw[draw=red, line width=2pt] (l5) -- (v51);
            \draw[draw=black, line width=1pt] (v51) -- (v52);
            \draw[draw=dgreen, line width=2pt] (v52) -- (v53);
            \draw[draw=black, line width=1pt] (v53) -- (r5);

            \draw[draw=black, line width=1pt] (l10) -- (v101);
            \draw[draw=blue, line width=2pt] (v101) -- (v102);
            \draw[draw=black, line width=1pt] (v102) -- (r10);

            \draw[draw=black, line width=1pt] (l9) -- (v91);
            \draw[draw=red, line width=2pt] (v91) -- (v92);
            \draw[draw=black, line width=1pt] (v92) -- (v93);
            \draw[draw=blue, line width=2pt] (v93) -- (v94);
            \draw[draw=black, line width=1pt] (v94) -- (v95);
            \draw[draw=dgreen, line width=2pt] (v95) -- (v96);
            \draw[draw=black, line width=1pt] (v96) -- (r9);

            \draw[draw=black, line width=1pt] (l8) -- (v81);
            \draw[draw=red, line width=2pt] (v81) -- (v82);
            \draw[draw=black, line width=1pt] (v82) -- (v83);
            \draw[draw=blue, line width=2pt] (v83) -- (v84);
            \draw[draw=black, line width=1pt] (v84) -- (v85);
            \draw[draw=dgreen, line width=2pt] (v85) -- (r8);

            \draw[draw=black, line width=1pt] (l7) -- (v71);
            \draw[draw=red, line width=2pt] (v71) -- (v72);
            \draw[draw=black, line width=1pt] (v72) -- (v73);
            \draw[draw=blue, line width=2pt] (v73) -- (r7);

            \draw[draw=black, line width=1pt] (l6) -- (v61);
            \draw[draw=red, line width=2pt] (v61) -- (r6);

            \draw[draw=blue, line width=2pt, ->]  (v41) .. controls (1.25,1.75) and (1.25,2.25) .. (v101);
            \draw[draw=red, line width=2pt, ->]  (v51) .. controls (0.75,2.25) and (0.75,2.75) .. (v91);
            
            \draw[draw=dgreen, line width=2pt, ->]  (v31) -- (v42);
            \draw[draw=red, line width=2pt, ->]  (v92) -- (v81);
            
            \draw[draw=dviolet, line width=2pt, ->]  (v21) -- (v32);
            \draw[draw=red, line width=2pt, ->]  (v82) -- (v71);
            
            \draw[draw=dorange, line width=2pt, ->]  (v11) -- (v22);
            \draw[draw=red, line width=2pt, ->]  (v72) -- (v61);
        
            \draw[draw=dgreen, line width=2pt, ->]  (v43) -- (v52);
            \draw[draw=blue, line width=2pt, ->]  (v102) -- (v93);
        
            \draw[draw=dviolet, line width=2pt, ->]  (v33) -- (v44);
            \draw[draw=blue, line width=2pt, ->]  (v94) -- (v83);
        
            \draw[draw=dorange, line width=2pt, ->]  (v23) -- (v34);
            \draw[draw=blue, line width=2pt, ->]  (v84) -- (v73);
            
            \draw[line width=2pt, ->]  (v45) .. controls (5.25,1.75) and (5.25,2.25) .. (v103);
            \draw[draw=dgreen, line width=2pt, ->]  (v53) .. controls (4.75,2.25) and (4.75,2.75) .. (v95); 
        
            \draw[line width=2pt, ->]  (v35) -- (v46);
            \draw[draw=dgreen, line width=2pt, ->]  (v96) -- (v85) ;
        
            \draw[line width=2pt, ->]  (v47) -- (v54);
            \draw[line width=2pt, ->]  (v104) -- (v97);
        \end{tikzpicture}}
    \caption{ The collection of non-intersecting paths in the LGV diagram corresponding to the minor $M_{2,1}(A)$ for $n = 5$.}
    \label{fig:PathCollectionExample}
\end{figure}
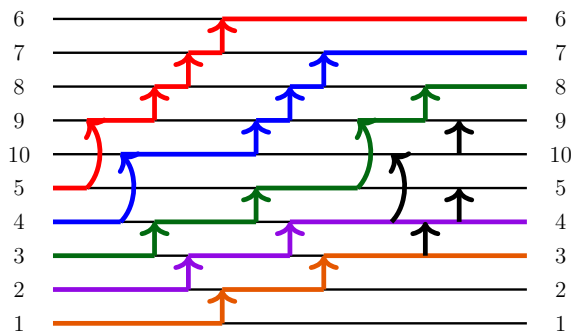

\begin{theorem}
    \label{thm:Main}
  A skew-symmetric matrix $A \in \SS_n$ is totally positive if and only if
  \[
    M_{j,k}(A) > 0  \quad \text{for any } 1 \leq j \leq k \leq n-1.
  \]
  This test is minimal in the sense that it uses the fewest possible number of inequalities.
\end{theorem}

The set $\SS_n^{\geq 0}$ of totally nonnegative skew-symmetric matrices is the Euclidean closure of $\SS_n^{>0}$. While the minors $M_{j,k}$ are non-negative on $\SS_n^{\geq 0}$, there exist skew-symmetric matrices $A \not \in \SS_n^{\geq 0}$ with $M_{j,k}(A) = 0$ for all $1 \leq j \leq k \leq n-1$. So the minors $M_{j,k}(A)$ are not enough to test for the nonnegativity of $A$. Nonetheless, together with the semigroup property of $\SO^{>0}(2n)$ we are able to give a nonnegativity test in the following~form.

\begin{theorem}\label{thm:Nonnegative}
    Fix $X \in \OGr(n,2n)$. Then, for any smooth $1$-parameter family $Z(\epsilon)$ in $\SO^{>0}(2n)$ such that $Z(\epsilon) \xrightarrow[\epsilon \to 0] {} \Id_{2n} $ and $X(\epsilon) \coloneqq X \cdot Z(\epsilon)$, the following are equivalent.
    \begin{enumerate}[wide=40pt, leftmargin = 58pt]
        \item \label{nonnegativeitem1} $X$ is totally nonnegative.
        \item \label{nonnegativeitem2} $X(\epsilon)$ is totally positive for all $\epsilon>0$ sufficiently small.
        \item \label{nonnegativeitem3} For all $1\leq j\leq k\leq n-1$, the leading coefficient in the Taylor expansion of $M_{j,k}(B(\epsilon))$ is positive, where $B(\epsilon)$ is defined by $X(\epsilon) = \rowspan \big([\Id_n|B(\epsilon)]\big)$.
    \end{enumerate}
    Moreover, the family $Z(\epsilon)$ can be chosen so that $M_{j,k}(B(\epsilon))$ are polynomials in $\epsilon$.
\end{theorem}

As for flag varieties, the set $\SS_{n}^{\geq 0}$ decomposes into a disjoint union of semi-algebraic sets called \emph{positive Richardson cells} as follows:
\begin{equation}
    \SS_n^{\geq 0} = \bigsqcup \mathcal{R}^{>0}_{v, w},
\end{equation} 
where the union is over all minimal coset representatives $w$ in the parabolic quotient $W^{[n-1]}$ of the Weyl group of $\SO(2n)$, and all $v \leq w$ in Bruhat order. See \Cref{subsec:RichardsonDeodhar} for more details. Our next result determines the Richardson cell that contains a given $A \in \SS_n^{\geq 0}$. A constructive version of this theorem in stated in \Cref{thm:RichardsonRestated}.
 
\begin{theorem}\label{thm:Richardson}
    Let $A \in \SS^{\geq 0}_n$ and $\mathscr{M}_A$ be the realizable rank $n$ matroid on $[2n]$ associated to $[\Id_n|A]$. Then, the Richardson cell containing $A$ can be determined from $\mathscr{M}_A$. 
\end{theorem}

Given a skew-symmetric matrix $A \in \SS_n$, its principal minors are perfect squares whose square root is a polynomial in the entries of $A$. These polynomials are called the \emph{Pfaffians} of $A$. As described here, there is a sign ambiguity for Pfaffians. However, in \Cref{sec:5}, we give a more intrinsic definition that fixes the sign. Given $I \subset [n]$ we denote by $\Pf_I(A)$ the Pfaffian of $A$ corresponding to the principal minor $\det(A_I^I)=\Pf_I(A)^2$. We take the convention that $\Pf_{\emptyset}(A)=1$ and $\Pf(A):=\Pf_{[n]}(A)$; also note that if $I$ has odd size, $\Pf_I(A)=0$. Similar to positive definite symmetric matrices, whose principal minors are positive, one could alternatively consider defining positive skew-symmetric matrices in terms of their Pfaffians. Remarkably, it turns out that the Pfaffians do have a fixed sign on the $\SS_n^{>0}$: 

\begin{theorem}\label{thm:PfaffianSign}
    For any $A \in \SS_n^{>0}$, and $I \subset [n]$  of even size, we have
        \begin{equation}\label{eq:signPattern}
            \sgn(\Pf_I(A)) = (-1)^{\sum_{i \in I} i - \binom{|I|+1}{2}} =: \sgn(I, [n]).
        \end{equation}
    If $I=\{i_1<\cdots < i_{|I|}\}$ and $[n]\setminus I=\{j_1<\cdots < j_{n-|I|}\}$, this is the sign of the permutation $i_1, \dots, i_{|I|}, j_{1} \dots, j_{n - |I|}$ in one-line notation.
\end{theorem}
\noindent
However, we note that there also exist skew-symmetric matrices that are not totally positive, or even totally nonnegative, whose Pfaffians have this sign pattern.

\subsection{Outline} We begin by collecting some necessary background and a few preliminary results in \Cref{sec:2}. In \Cref{sec:3}, we prove the positivity criterion for skew-symmetric matrices in \Cref{thm:Main}. We deal with the set of nonnegative skew-symmetric matrices $\SS_n^{\geq 0}$ and prove \Cref{thm:Nonnegative} and \Cref{thm:Richardson} in \Cref{sec:4}. In \Cref{sec:5} we show that the Pfaffians of a skew-symmetric matrix $A \in \SS_n^{\geq 0}$ have a remarkable sign pattern and briefly discuss Pfaffian positivity. Finally, in \Cref{sec:6}, we discuss future directions and open questions. 

\smallskip

Some computations in this article were carried on using the computer algebra system Macaulay2 \cite{M2}. The relevant Macaulay2 script we used is available at 

\begin{center}
\href{https://mathrepo.mis.mpg.de/PositiveSkewMatrices}{\tt https://mathrepo.mis.mpg.de/PositiveSkewMatrices}.     
\end{center}

\smallskip

\subsection{Acknowledgment}
We thank Bernd Sturmfels for suggesting this problem, Chris Eur for helpful Macaulay2 code, Steven Karp for pointing us to \cite[Proposition 8.17]{Lusztig1} and \cite[Theorem 3.4]{Lusztig2} used in the proof of \Cref{lem:semiGroup}, and Grant Barkley for helpful discussions around generalized minors. YEM was partially supported by Deutsche Forschungsgemeinschaft (DFG, German Research Foundation) SFB-TRR 195 “Symbolic Tools in Mathematics and their Application.”

\section{Preliminaries}\label{sec:2}

\subsection{A pinch of Lie theory}\label{subsec:LieTheory}
Let $G$ be a split reductive algebraic group over $\RR$. A \emph{pinning} for $G$ is a choice of a maximal torus $T$, a Borel subgroup $B$ with opposite Borel group $B^{-}$, and a group homomorphism $\phi_i: \SL_2 \to G$ associated to each simple root. We will focus on the case $G=\SO(2n)$, the special orthogonal group with respect to the bilinear form $q$. Let $\mathfrak{so}_{2n}$ denote its Lie algebra, the vector space of $2n\times 2n$ matrices $X$ such that $X^T Q + QX=0$ with the usual Lie bracket. 

\smallskip

We fix the maximal torus to be the subgroup of $\SO(2n)$  of diagonal matrices
\begin{equation}\label{eq:maxTorus}
T := \big\{\mathrm{diag}(t_1,\ldots,t_n,t^{-1}_1,\ldots,t^{-1}_n): \ t_1,\ldots, t_n \in \RR^* \big \}.    
\end{equation}
We fix the Borel subgroup $B$ to be the upper triangular matrices in $\SO(2n)$. The opposite Borel subgroup is the subgroup of lower triangular matrices in $\SO(2n)$. We take, as simple roots, $\{\varepsilon_1-\varepsilon_2, \varepsilon_2-\varepsilon_3, \ldots, \varepsilon_{n-1}-\varepsilon_{n}, \varepsilon_{n-1}+\varepsilon_n\}$ where $\{\varepsilon_i\}_{i=1}^{n}$ is the standard basis of $\RR^n$. The group homomorphisms $\big(\phi_i: \SL(2) \to \SO(2n)\big)_{i=1}^{n}$ for our pinning~are
\begin{equation}\label{eq:pinning}
\resizebox{0.85\textwidth}{!}{$
    \begin{aligned}
        & \phi_i\begin{pmatrix}
                a&b\\
                c&d
            \end{pmatrix}=\begin{blockarray}{cccccccccc}
            \begin{block}{c(ccccccccc)}
              & 1 &  &  & & & & & &\\
              & & \ddots &  & &  & & & & \\
              i& &  & a & b &  & & & & \\
              i+1& &  &c  & d &  &  & & & \\
              & &  &  &  & \ddots  & & & &  \\
              n+i&&&&&&d&-c&&\\
              n+i+1&&&&&&-b&a&&\\
             &&&&&&&&\ddots&\\
             &&&&&&&&&1\\
            \end{block} 
            \end{blockarray}\hspace{5pt} \qquad \text{for } 1 \leq i \leq n-1,\\
        &\text{and}\\
        &  \phi_n\begin{pmatrix}
                a&b\\
                c&d
            \end{pmatrix}=\begin{blockarray}{cccccccccc}
        & &  &n-1 & n & & & & 2n-1 &2n\\
        \begin{block}{c(ccccccccc)}
          & 1 &  &  &  & & & & &\\
          & & \ddots &   &  & & & & &\\
          n-1& &  & a &    & & & && b\\
          n& &  &  & a &   & & &-b &\\
          &&&&&1&&&&\\
          &&&&&&\ddots&&&\\
         2n-1&&&&-c&&&&d&\\
         2n&&&c&&&&&&d\\
        \end{block} 
        \end{blockarray}\hspace{5pt}.
\end{aligned}$} 
\end{equation}
The Weyl group of $\operatorname{SO}(2n)$ is the group $W$ of permutations $\sigma \in S_{2n}$ satisfying
\begin{enumerate}[wide=60pt]
    \item \ $\quad \sigma(n+i)=n+\sigma(i)$ modulo $2n$  for all $i \in [n]$,   
    \item \ $\quad \big| \{i \in [n]: \sigma(i)>n\} \big|$ is even.
\end{enumerate}
This is a Coxeter group; it is generated by
\begin{equation}\label{eq:WeylGroupGens}
\resizebox{0.90\textwidth}{!}{
    $s_i=(i, i+1)(n+i, n+i+1) \quad \text{for } i=1,\ldots n-1 \quad \text{and} \quad s_n=(n,2n-1)(n-1,2n)
    $},   
\end{equation}
We note that these generators satisfy the following braid relations
\begin{align*}
    &s_i^2={\rm id}, \ s_{j}s_{j+1}s_j=s_{j+1}s_js_{j+1} \ \text{ and } \ s_{n-2}s_ns_{n-2}=s_{n-2}s_ns_{n-2} \quad \text{for any } i \in [n], j \in [n-1],\\
    &s_is_j = s_js_i \quad \text{whenever } 1\leq i \neq j \leq n \text{ and } \{i,j\} \notin \big\{\{k,k+1\}_{k\in [n-1]}\big\} \cup  \big\{ \{n-2,n\} \big\}.
\end{align*}

Given $w \in W$ the \emph{length of $w$} is $\ell(w) := \min\{k \in \NN: w=s_{i_1}\cdots s_{i_k}$ where $i_j \in [n]\}$. A  \emph{reduced expression of $w$} is a particular choice of (ordered) simple transpositions $s_{i_j}$ such that $w=s_{i_1}\cdots s_{i_{\ell(w)}}$. We denote by $w_0$ the element of maximal length in $W$. 

\medskip

In \cite[Theorem 11.3]{MR} Marsh and Rietsch give a parametrization of the totally positive part of the complete flag variety for any split reductive group $G$ and choice of pinning in terms of a reduced expression for the longest element in the Weyl group. Lusztig's characterization of positivity for flag varieties \cite{Lusztig2} guarantees that we can project this parametrization from the complete flag variety of $G$ down to any partial flag variety.

\begin{remark}\label{rem:MRParamConvention}
    We use a different convention than Marsh and Rietsch by viewing flags as row spans of matrices rather than column spans. This allows us to present examples more easily. 
    Concretely, a matrix $M$ represents a positive flag in our convention if and only if its transpose $M^T$ represents a positive flag in~\cite{MR}.
\end{remark}

\smallskip

In this section, we rely on some basic facts about the Weyl group of $\SO(2n)$, which can be found in \cite{HumphreysCoxeter}. For $t \in \RR$ and $1\leq i \leq n$ we denote by $x_i(t)$ the matrices
\begin{equation}\label{eq:x_i(t)}
x_i(t) := \phi_i\begin{pmatrix}
        1&t\\
        0&1
    \end{pmatrix}.    
\end{equation}

Projecting the Marsh-Rietsch parametrization \cite[Section 11]{MR} for $G=\SO(2n)$, we obtain
\begin{equation} \label{eq:parameterizationCompleteFlag}
    \OGr(n,2n)^{>0} = \Big\{ \rowspan \big(\pi_n \big(x_{j_1}(t_1)\cdots x_{j_{\hat{N}}}(t_{\hat{N}}) \big)\big): t_1,\ldots, t_{\hat{N}} \in \RR_{>0} \Big\},
\end{equation}
where $\pi_n:\operatorname{SO}(2n) \to \OGr(n,2n)$ is the map that takes a $2n \times 2n$ orthogonal matrix to the row span of its first $n$ rows, and $w_0 = s_{j_1}s_{j_2} \cdots s_{j_{\hat{N}}}$ is the longest element of the Weyl group, and as such, $\hat{N} = 2 \binom{n}{2}$.

We will find it useful to fix a reduced expression $\underline{w}_0$ for $w_0$, that is, a minimal sequence of generators which multiplied together yield $w_0$. For reasons that will be evident shortly, we also give a name to the final $\binom{n}{2}$ generators of $\underline{w}_0$, as follows:

    \begin{equation}\label{eq:w_0ReducedExpr}
      \resizebox{.915\hsize}{!}{$
    \begin{aligned}
    \underline{w}_0^{[n-1]} &= s_n(s_{n-2}\cdots s_1)(s_{n-1}\cdots s_2)s_n(s_{n-2}\cdots s_3)(s_{n-1}\cdots s_4)\cdots s_n && \text{for $n$ even,}\\
    \underline{w_0}^{[n-1]} &=  s_n(s_{n-2}\cdots s_1)(s_{n-1}\cdots s_2)s_n(s_{n-2}\cdots s_3)(s_{n-1}\cdots s_4)\cdots s_n s_{n-2} s_{n-1}&& \text{for $n$ odd,}
    \end{aligned}$}
    \end{equation}
    \begin{equation}\label{eq:w_0LongReducedExpr}
    \underline{w}_0 = s_1(s_2s_1) \dots (s_{n-1}s_{n-2}\dots s_1) \ \underline{w}_0^{[n-1]}.
    \end{equation}

Observe that, for $Z\in \SO(2n)$, $i\in [n-1]$, and $t\in \mathbb{R}$, $\pi_n(Z)=\pi_n(x_i(t)Z)$. Thus, instead of using an expression for $w_0$ in the Marsh-Rietsch parametrization of $\OGr(n,2n)^{>0}$, it suffices to use any representative of the right coset $w_0W_{[n-1]}$, where $W_{[n-1]}\coloneqq \langle s_1,\ldots,s_{n-1}\rangle$. We will use a minimum length coset representative. Such a representative should have length $N\coloneqq \binom{n}{2}$. Note that, by equation \eqref{eq:w_0LongReducedExpr}, $\underline{w}_0^{[n-1]}$ is an expression for a right coset representative of $w_0$ and, by equation \eqref{eq:w_0ReducedExpr}, it is of length $N$. We have

\begin{equation}\label{eq:posorthogonalgrassmannian}
    \OGr(n,2n)^{>0} = \Big\{ \rowspan \big(\pi_n \big(x_{i_1}(t_1)\cdots x_{i_N}(t_{N}) \big)\big): t_1,\ldots, t_{N} \in \RR_{>0} \Big\},
\end{equation}
where $\underline{w}_0^{[n-1]}=s_{i_1}s_{i_2}\cdots s_{i_N}$. We denote by $X = X(t_1, \dots, t_N) $ the $n \times 2n$ matrix in the first $n$ rows of $x_{i_1}(t_1)\cdots x_{i_N}(t_N)$. Since the matrices \eqref{eq:x_i(t)} are all unipotent, the leftmost $n\times n$ block matrix of $X$ is unipotent. We may thus row reduce $X$ to the form $[\Id_n|A]$, where $A=A(t_1,\dots,t_N)$ is skew-symmetric matrix with polynomial entries in the $t_i$.

\begin{definition}\label{def:LusztigPositive}
    A real skew-symmetric $n \times n$ matrix is \textit{totally positive} if it is of the form $A(t_1, \dots, t_N)$ with $t_1,\dots, t_N \in \RR_{> 0}$.
\end{definition}

\begin{example}[$n=4$]\label{ex:n=4MRparam} In this case we have $ N = \binom{4}{2} = 6$ parameters. The expressions for $w_0$ and $w_0^{[3]}$ are $\underline{w}_0 =s_{1}(s_{2}s_{1})(s_{3}s_{2}s_{1})s_{4}(s_{2}s_{1})(s_{3}s_{2})s_{4}$ and $\underline{w}_0^{[3]} =  s_{4}(s_{2}s_{1})(s_{3}s_{2})s_{4}$, meaning

\[
    i_1 = 4 , \quad i_2 = 2 \quad i_3 = 1 \quad i_4 = 3 \quad i_5 = 2 \quad i_6 = 4.
\]
The matrices $X$ and $A$ are given by
\[
X =  \begin{pmatrix}
       1&t_{3}&t_{3}t_{5}&0&0&0&0&t_{3}t_{5}t_{6}\\
      0&1&t_{2}+t_{5}&t_{2}t_{4}&0&0&-t_{2}t_{4}t_{6}&t_{2}t_{6}+t_{5}t_{6}\\
      0&0&1&t_{4}&0&t_{1}t_{4}t_{5}&-t_{1}t_{4}-t_{4}t_{6}&t_{1}+t_{6}\\
      0&0&0&1&-t_{1}t_{2}t_{3}&t_{1}t_{2}+t_{1}t_{5}&-t_{1}-t_{6}&0

\end{pmatrix}
\]
and
\[
A = \begin{pmatrix}
        0&t_{1}t_{2}t_{3}t_{4}t_{5}&-t_{1}t_{2}t_{3}t_{4}&t_{1}t_{2}t_{3}\\
      -t_{1}t_{2}t_{3}t_{4}t_{5}&0&t_{1}t_{2}t_{4}&-t_{1}t_{2}-t_{1}t_{5}\\
      t_{1}t_{2}t_{3}t_{4}&-t_{1}t_{2}t_{4}&0&t_{1}+t_{6}\\
      -t_{1}t_{2}t_{3}&t_{1}t_{2}+t_{1}t_{5}&-t_{1}-t_{6}&0

\end{pmatrix}.
\]
In this parametrization, the minors $M_{j,k}(A)$, which we computed in Example \ref{ex:n=4Minors}, are
  \begin{alignat*}{3}
    M_{1,1}(A) &= t_{1}^{2}t_{2}^{2}t_{3}^{2}t_{4}^{2}t_{5}^{2}t_{6},   & &   \\
    M_{1,2}(A) &= t_{1}^{2}t_{2}^{2}t_{3}^{2}t_{4}^{2}t_{5}t_{6},   \quad M_{2,2}(A) &=& \ t_{1}^{2}t_{2}^{2}t_{3}^{2}t_{4}t_{5}, \quad  \\
    M_{1,3}(A) &= t_{1}^{2}t_{2}^{2}t_{3}t_{4}^{2}t_{5}t_{6},   \quad M_{2,3}(A) &=& \ t_{1}^{2}t_{2}t_{3}t_{4}t_{5}, \quad M_{3,3}(A) &= t_{1}t_{2}t_{3}.
    \end{alignat*}
    Note that these minors are all monomials in the parameters $t_1, \ldots, t_6$ and that if all of these minors are positive, then all the $t_i$'s are positive. This essentially proves Theorem \ref{thm:Main} in the case $n=4$ and motivates our general approach. Moreover, since the exponent vectors of the $M_{j,k}$ form a $\ZZ$-basis for $\ZZ^6$, each $t_i$ is a Laurent monomial in the minors $M_{j,k}(A)$.
\end{example}

\subsection{The Lindstr\"om-Gessel-Viennot diagram}

In this section, we introduce the Lindstr\"om-Gessel-Viennot (LGV) diagram which will be useful in the remainder of this article. Loosely speaking, this diagram encodes the combinatorics that govern the multiplication of matrices $x_{i_1}(t_1) \cdots x_{i_N}(t_N)$ and hence also of $X = X(t_1,\dotsm t_N)$ and $A = A(t_1,\dotsm t_N)$. This construction is similar to \cite[Section 3.2]{boretsky}, in which the matrices $x_i(t)$ are seen as adjacency matrices of a weighted directed graph.

\begin{definition}\label{def:LGVdiagram}
    We define the \textit{LGV diagram} as a directed graph consisting of $2n$ horizontal edges called \textit{strands}, with $2\binom{n}{2}$ weighted vertical arrows connecting them in the following way. Let $i_1,\ldots, i_N$ be such that $\underline{w}_0^{[n-1]}=s_{i_1}s_{i_2}\cdots s_{i_N}$.
    \begin{enumerate}
        \item Strands are always directed rightwards. We call the ends of the strands \textit{vertices}. Label the vertices on both sides in order $1,\ldots, n, 2n,\ldots, n+1$ starting from the bottom. Vertices on the left are called \textit{source vertices} and those on the right are called \textit{sink vertices}.
        \smallskip
        \item For each $1\leq l\leq N$, we add arrows to the strands, with arrows corresponding to smaller values of $l$ always appearing to the left of those corresponding to larger values of $l$. To obtain the arrows corresponding to $l$, look at the non-diagonal entries of $x_{i_l}(t_l)$. For every non-zero entry in position $(j,k)$ we draw an arrow from strand $j$ to strand $k$ with weight given by the value of the entry, which is either $t_l$ or $-t_l$.
    \end{enumerate}
    We stress that the order in which the arrows from different $x_{i_l}(t_l)$ matters, reflecting the noncommutativity of these matrices. However, the multiple arrows corresponding to a single matrix $x_{i_l}(t_l)$ can be drawn in any order.
\end{definition}

\begin{example}[$n=4$]\label{ex:n=4LGV}
   In this case, $(i_1, i_2, i_3, i_4, i_5, i_6) = (4,2,1,3,2,4)$ and the matrices $x_1, x_2, x_3, x_4$ are, respectively,
    \[
    \resizebox{1\textwidth}{!}{$
    \begin{bmatrix}
        1 & t &&&&&&\\
         & 1 &&&&&&\\
         &&1&&&&&\\
         &&&1&&&&\\
         &&&&1&&&\\
         &&&&-t&1&&\\
         &&&&&&1&\\
         &&&&&&&1\\
    \end{bmatrix},\quad
    \begin{bmatrix}
        1 &  &&&&&&\\
         & 1 & t&&&&&\\
         &&1&&&&&\\
         &&&1&&&&\\
         &&&&1&&&\\
         &&&&&1&&\\
         &&&&&-t&1&\\
         &&&&&&&1\\
    \end{bmatrix}, \quad
    \begin{bmatrix}
         1&&&&&&&\\
         &1&&&&&&\\
         &&1&t&&&&\\
         &&&1&&&&\\
         &&&&1&&&\\
         &&&&&1&&\\
         &&&&&&1&\\
         &&&&&&-t&1\\
    \end{bmatrix} 
    \quad \text{and} \quad 
    \begin{bmatrix}
        1 &  &&&&&&\\
         & 1 &&&&&&\\
         &&1&&&&&t\\
         &&&1&&&-t&\\
         &&&&1&&&\\
         &&&&&1&&\\
         &&&&&&1&\\
         &&&&&&&1\\
    \end{bmatrix}.$}
    \]
    Since $\underline{w}_0^{[3]} = s_{4}(s_{2}s_{1})(s_{3}s_{2})s_{4}$, the LGV diagram for $n=4$ is as follows. As the horizontal strands are always directed rightwards, we do not explicitly indicate it in our drawings.

    \begin{center}
    \scalebox{1}{\begin{tikzpicture}  
    \coordinate (l7) at (0,4);
    \coordinate (l8) at (0,3.5);
    \coordinate (l9) at (0,3);
    \coordinate (l10) at (0,2.5);        
    \coordinate (l5) at (0,2);        
    \coordinate (l4) at (0,1.5);    
    \coordinate (l3) at (0,1);
    \coordinate (l2) at (0,0.5);
       
    \coordinate (r7) at (5.5,4);
    \coordinate (r8) at (5.5,3.5);
    \coordinate (r9) at (5.5,3);
    \coordinate (r10) at (5.5,2.5);
    \coordinate (r5) at (5.5,2);        
    \coordinate (r4) at (5.5,1.5);    
    \coordinate (r3) at (5.5,1);
    \coordinate (r2) at (5.5,0.5);

    \coordinate (v21) at (2,0.5);

    \coordinate (v31) at (1.5,1);
    \coordinate (v32) at (2,1);
    \coordinate (v33) at (3.5,1);
    
    \coordinate (v41) at (1,1.5);
    \coordinate (v42) at (1.5,1.5);
    \coordinate (v43) at (3,1.5);
    \coordinate (v44) at (3.5,1.5);
    \coordinate (v45) at (5,1.5);
    
    \coordinate (v51) at (0.5,2);
    \coordinate (v52) at (3,2);
    \coordinate (v53) at (4.5,2);

    \coordinate (v101) at (1,2.5);    
    \coordinate (v102) at (3,2.5);
    \coordinate (v103) at (5,2.5);

    \coordinate (v91) at (0.5,3);
    \coordinate (v92) at (1.5,3);
    \coordinate (v93) at (3,3);
    \coordinate (v94) at (3.5,3);
    \coordinate (v95) at (4.5,3);

    \coordinate (v81) at (1.5,3.5);
    \coordinate (v82) at (2,3.5);
    \coordinate (v83) at (3.5,3.5);
    \coordinate (v84) at (4,3.5);

    \coordinate (v71) at (2,4);

    \draw[black] (-0.5,0.5) node  [xscale = 0.8, yscale = 0.8] {$1$};
    \draw[black] (-0.5,1) node  [xscale = 0.8, yscale = 0.8] {$2$};
    \draw[black] (-0.5,1.5) node  [xscale = 0.8, yscale = 0.8] {$3$};
    \draw[black] (-0.5,2) node  [xscale = 0.8, yscale = 0.8] {$4$};
    \draw[black] (-0.5,2.5) node  [xscale = 0.8, yscale = 0.8] {$8$};
    \draw[black] (-0.5,3) node  [xscale = 0.8, yscale = 0.8] {$7$};
    \draw[black] (-0.5,3.5) node  [xscale = 0.8, yscale = 0.8] {$6$};
    \draw[black] (-0.5,4) node  [xscale = 0.8, yscale = 0.8] {$5$};
    
    \draw[black] (6,0.5) node  [xscale = 0.8, yscale = 0.8] {$1$};
    \draw[black] (6,1) node  [xscale = 0.8, yscale = 0.8] {$2$};
    \draw[black] (6,1.5) node  [xscale = 0.8, yscale = 0.8] {$3$};
    \draw[black] (6,2) node  [xscale = 0.8, yscale = 0.8] {$4$};
    \draw[black] (6,2.5) node  [xscale = 0.8, yscale = 0.8] {$8$};
    \draw[black] (6,3) node  [xscale = 0.8, yscale = 0.8] {$7$};
    \draw[black] (6,3.5) node  [xscale = 0.8, yscale = 0.8] {$6$};
    \draw[black] (6,4) node  [xscale = 0.8, yscale = 0.8] {$5$};

    \draw[draw=black, line width=1pt] (l2) -- (r2);
    \draw[draw=black, line width=1pt] (l3) -- (r3);
    \draw[draw=black, line width=1pt] (l4) -- (r4);
    \draw[draw=black, line width=1pt] (l5) -- (r5);
    \draw[draw=black, line width=1pt] (l10) -- (r10);
    \draw[draw=black, line width=1pt] (l9) -- (r9);
    \draw[draw=black, line width=1pt] (l8) -- (r8);
    \draw[draw=black, line width=1pt] (l7) -- (r7);

    \draw[draw=black, line width=1pt, ->]  (v41) .. controls (1.25,1.75) and (1.25,2.25) .. (v101) node[midway, below left] {\tiny $t_{1}$};
    
    \draw[draw=black, line width=1pt, ->]  (v51) .. controls (0.75,2.25) and (0.75,2.75) .. (v91) node[midway, below left] {\tiny $-t_{1}$};
    
    \draw[draw=black, line width=1pt, ->]  (v31) -- (v42) node[below right] {\tiny $t_{2}$};
    \draw[draw=black, line width=1pt, ->]  (v92) -- (v81) node[below right] {\tiny $-t_{2}$};
    
    \draw[draw=black, line width=1pt, ->]  (v21) -- (v32) node[below right] {\tiny $t_{3}$};
    \draw[draw=black, line width=1pt, ->]  (v82) -- (v71) node[below right] {\tiny $-t_{3}$};

    \draw[draw=black, line width=1pt, ->]  (v43) -- (v52) node[below right] {\tiny $t_{4}$};
    \draw[draw=black, line width=1pt, ->]  (v102) -- (v93) node[below right] {\tiny $-t_{4}$};

    \draw[draw=black, line width=1pt, ->]  (v33) -- (v44) node[below right] {\tiny $t_{5}$};
    \draw[draw=black, line width=1pt, ->]  (v94) -- (v83) node[below right] {\tiny $-t_{5}$};
    
    \draw[draw=black, line width=1pt, ->]  (v45) .. controls (5.25,1.75) and (5.25,2.25) .. (v103)node[midway, below left] {\tiny $t_{6}$};
    \draw[draw=black, line width=1pt, ->]  (v53) .. controls (4.75,2.25) and (4.75,2.75) .. (v95) node[midway, below left] {\tiny $-t_{6}$}; 

\end{tikzpicture}}
    \end{center}
\end{example}

\begin{definition}
    A \textit{path} in the LGV diagram is a path in the directed graph from a source vertex to a sink vertex. We will identify a path $\mathcal{P}$ with the set of vertical arrows it uses. Observe that this set uniquely defines $\mathcal{P}$ 
    
    \smallskip
    
    Given $I,J \subset [2n]$ of the same cardinality, we define a \textit{path collection} from $I$ to $J$ as $|I|$ paths from the source vertices labeled by elements of $I$ to the sink vertices labeled by elements of $J$. We say a path collection is \textit{non-intersecting} if the paths do not intersect when drawn in the LGV diagram. We say that a non-intersecting path collection from $I$ to $J$ is \textit{unique} if it is the only non-intersecting path collection with source vertices $I$ and sink vertices $J$.

    \smallskip
    
    A \textit{left greedy path} in the LGV diagram originating from $v$ is a path that starts from source vertex $v$ and uses each vertical arrow it encounters until it reaches a sink vertex. These are paths that turn left whenever possible. We denote this path by ${\rm LGP}(v)$. A \textit{left greedy path collection} is a path collection where all the paths are left greedy.

    \smallskip

    We say that a path  $\mathcal{P}$ in a non-intersecting path collection $P$ is \textit{right greedy} if it never uses vertical arrows unless that is the only way to avoid intersecting paths originating below it in $P$. These are paths that turn right whenever possible, subject to $P$ being a non-intersecting path collection.

\end{definition}

\begin{figure}[H]
    \centering
     \scalebox{0.9}{\begin{tikzpicture}
            \coordinate (l6) at (0,4.5);    
            \coordinate (l7) at (0,4);
            \coordinate (l8) at (0,3.5);
            \coordinate (l9) at (0,3);
            \coordinate (l10) at (0,2.5);        
            \coordinate (l5) at (0,2);        
            \coordinate (l4) at (0,1.5);    
            \coordinate (l3) at (0,1);
            \coordinate (l2) at (0,0.5);
            \coordinate (l1) at (0,0);
        
            \coordinate (r6) at (7,4.5);    
            \coordinate (r7) at (7,4);
            \coordinate (r8) at (7,3.5);
            \coordinate (r9) at (7,3);
            \coordinate (r10) at (7,2.5);
            \coordinate (r5) at (7,2);        
            \coordinate (r4) at (7,1.5);    
            \coordinate (r3) at (7,1);
            \coordinate (r2) at (7,0.5);
            \coordinate (r1) at (7,0);

            \coordinate (v11) at (2.5,0);
            
            \coordinate (v21) at (2,0.5);
            \coordinate (v22) at (2.5,0.5);
            \coordinate (v23) at (4,0.5);
        
            \coordinate (v31) at (1.5,1);
            \coordinate (v32) at (2,1);
            \coordinate (v33) at (3.5,1);
            \coordinate (v34) at (4,1);
            \coordinate (v35) at (5.5,1);
            
            \coordinate (v41) at (1,1.5);
            \coordinate (v42) at (1.5,1.5);
            \coordinate (v43) at (3,1.5);
            \coordinate (v44) at (3.5,1.5);
            \coordinate (v45) at (5,1.5);
            \coordinate (v46) at (5.5,1.5);
            \coordinate (v47) at (6,1.5);
            
            \coordinate (v51) at (0.5,2);
            \coordinate (v52) at (3,2);
            \coordinate (v53) at (4.5,2);
            \coordinate (v54) at (6,2);

            \coordinate (v101) at (1,2.5);    
            \coordinate (v102) at (3,2.5);
            \coordinate (v103) at (5,2.5);
            \coordinate (v104) at (6,2.5);

            \coordinate (v91) at (0.5,3);
            \coordinate (v92) at (1.5,3);
            \coordinate (v93) at (3,3);
            \coordinate (v94) at (3.5,3);
            \coordinate (v95) at (4.5,3);
            \coordinate (v96) at (5.5,3);
            \coordinate (v97) at (6,3);

            \coordinate (v81) at (1.5,3.5);
            \coordinate (v82) at (2,3.5);
            \coordinate (v83) at (3.5,3.5);
            \coordinate (v84) at (4,3.5);
            \coordinate (v85) at (5.5,3.5);

            \coordinate (v71) at (2,4);
            \coordinate (v72) at (2.5,4);
            \coordinate (v73) at (4,4);
        
            \coordinate (v61) at (2.5,4.5);
            
            \draw[black] (-0.5,0) node  [xscale = 0.8, yscale = 0.8] {$1$};
            \draw[black] (-0.5,0.5) node  [xscale = 0.8, yscale = 0.8] {$2$};
            \draw[black] (-0.5,1) node  [xscale = 0.8, yscale = 0.8] {$3$};
            \draw[black] (-0.5,1.5) node  [xscale = 0.8, yscale = 0.8] {$4$};
            \draw[black] (-0.5,2) node  [xscale = 0.8, yscale = 0.8] {$5$};
            \draw[black] (-0.5,2.5) node  [xscale = 0.8, yscale = 0.8] {$10$};
            \draw[black] (-0.5,3) node  [xscale = 0.8, yscale = 0.8] {$9$};
            \draw[black] (-0.5,3.5) node  [xscale = 0.8, yscale = 0.8] {$8$};
            \draw[black] (-0.5,4) node  [xscale = 0.8, yscale = 0.8] {$7$};
            \draw[black] (-0.5,4.5) node  [xscale = 0.8, yscale = 0.8] {$6$};
            
            \draw[black] (7.5,0) node  [xscale = 0.8, yscale = 0.8] {$1$};
            \draw[black] (7.5,0.5) node  [xscale = 0.8, yscale = 0.8] {$2$};
            \draw[black] (7.5,1) node  [xscale = 0.8, yscale = 0.8] {$3$};
            \draw[black] (7.5,1.5) node  [xscale = 0.8, yscale = 0.8] {$4$};
            \draw[black] (7.5,2) node  [xscale = 0.8, yscale = 0.8] {$5$};
            \draw[black] (7.5,2.5) node  [xscale = 0.8, yscale = 0.8] {$10$};
            \draw[black] (7.5,3) node  [xscale = 0.8, yscale = 0.8] {$9$};
            \draw[black] (7.5,3.5) node  [xscale = 0.8, yscale = 0.8] {$8$};
            \draw[black] (7.5,4) node  [xscale = 0.8, yscale = 0.8] {$7$};
            \draw[black] (7.5,4.5) node  [xscale = 0.8, yscale = 0.8] {$6$};

            \draw[draw=dorange, line width=2pt] (l1)  -- (v11);
            \draw[draw=black, line width=1pt] (v11) -- (r1);

            \draw[draw=dviolet, line width=2pt] (l2) -- (v21);
            \draw[draw=black, line width=1pt] (v21) -- (v22);
            \draw[draw=dorange, line width=2pt] (v22) -- (v23);
            \draw[draw=black, line width=1pt] (v23) -- (r2);

            \draw[draw=dgreen, line width=2pt] (l3) -- (v31);
            \draw[draw=black, line width=1pt] (v31) -- (v32);
            \draw[draw=dviolet, line width=2pt] (v32) -- (v33);
            \draw[draw=black, line width=1pt] (v33) -- (v34);
            \draw[draw=dorange, line width=2pt] (v34) -- (r3);

            \draw[draw=blue, line width=2pt] (l4) -- (v41);
            \draw[draw=black, line width=1pt] (v41) -- (v42);
            \draw[draw=dgreen, line width=2pt] (v42) -- (v43);
            \draw[draw=black, line width=1pt] (v43) -- (v44);
            \draw[draw=dviolet, line width=2pt] (v44) -- (r4);
            
            \draw[draw=red, line width=2pt] (l5) -- (v51);
            \draw[draw=black, line width=1pt] (v51) -- (v52);
            \draw[draw=dgreen, line width=2pt] (v52) -- (v53);
            \draw[draw=black, line width=1pt] (v53) -- (r5);

            \draw[draw=black, line width=1pt] (l10) -- (v101);
            \draw[draw=blue, line width=2pt] (v101) -- (v102);
            \draw[draw=black, line width=1pt] (v102) -- (r10);

            \draw[draw=black, line width=1pt] (l9) -- (v91);
            \draw[draw=red, line width=2pt] (v91) -- (v92);
            \draw[draw=black, line width=1pt] (v92) -- (v93);
            \draw[draw=blue, line width=2pt] (v93) -- (v94);
            \draw[draw=black, line width=1pt] (v94) -- (v95);
            \draw[draw=dgreen, line width=2pt] (v95) -- (v96);
            \draw[draw=black, line width=1pt] (v96) -- (r9);

            \draw[draw=black, line width=1pt] (l8) -- (v81);
            \draw[draw=red, line width=2pt] (v81) -- (v82);
            \draw[draw=black, line width=1pt] (v82) -- (v83);
            \draw[draw=blue, line width=2pt] (v83) -- (v84);
            \draw[draw=black, line width=1pt] (v84) -- (v85);
            \draw[draw=dgreen, line width=2pt] (v85) -- (r8);

            \draw[draw=black, line width=1pt] (l7) -- (v71);
            \draw[draw=red, line width=2pt] (v71) -- (v72);
            \draw[draw=black, line width=1pt] (v72) -- (v73);
            \draw[draw=blue, line width=2pt] (v73) -- (r7);

            \draw[draw=black, line width=1pt] (l6) -- (v61);
            \draw[draw=red, line width=2pt] (v61) -- (r6);

            \draw[draw=blue, line width=2pt, ->]  (v41) .. controls (1.25,1.75) and (1.25,2.25) .. (v101) node[midway, below left] {\scriptsize \textcolor{blue}{$t_{1}$}};
            
            \draw[draw=red, line width=2pt, ->]  (v51) .. controls (0.75,2.25) and (0.75,2.75) .. (v91) node[midway, below left] {\scriptsize \textcolor{red}{$-t_{1}$}};
            
            \draw[draw=dgreen, line width=2pt, ->]  (v31) -- (v42) 
             node[midway, right=0.5mm] {\scriptsize \textcolor{dgreen}{$t_{2}$}};

            \draw[draw=red, line width=2pt, ->]  (v92) -- (v81)
            node[midway, right=0.5mm] {\scriptsize \textcolor{red}{$-t_{2}$}};

            \draw[draw=dviolet, line width=2pt, ->]  (v21) -- (v32)
            node[midway, right=0.5mm] {\scriptsize \textcolor{dviolet}{$t_{3}$}};

            \draw[draw=red, line width=2pt, ->]  (v82) -- (v71)
            node[midway, right=0.5mm] {\scriptsize \textcolor{red}{$-t_{3}$}};
            
            \draw[draw=dorange, line width=2pt, ->]  (v11) -- (v22)
            node[midway, right=0.5mm] {\scriptsize \textcolor{orange}{$t_{4}$}};

            \draw[draw=red, line width=2pt, ->]  (v72) -- (v61)
            node[midway, right=0.5mm] {\scriptsize \textcolor{red}{$-t_{4}$}};

            \draw[draw=dgreen, line width=2pt, ->]  (v43) -- (v52)
            node[midway, right=0.5mm] {\scriptsize \textcolor{dgreen}{$t_{5}$}};
            
            \draw[draw=blue, line width=2pt, ->]  (v102) -- (v93)
            node[midway, right=0.5mm] {\scriptsize \textcolor{blue}{$-t_{5}$}};
        
            \draw[draw=dviolet, line width=2pt, ->]  (v33) -- (v44)
            node[midway, right=0.5mm] {\scriptsize \textcolor{dviolet}{$t_{6}$}};
            
            \draw[draw=blue, line width=2pt, ->]  (v94) -- (v83)
            node[midway, right=0.5mm] {\scriptsize \textcolor{blue}{$-t_{6}$}};
        
            \draw[draw=dorange, line width=2pt, ->]  (v23) -- (v34)
            node[midway, right=0.5mm] {\scriptsize \textcolor{orange}{$t_{7}$}};
            
            \draw[draw=blue, line width=2pt, ->]  (v84) -- (v73)
            node[midway, right=0.5mm] {\scriptsize \textcolor{blue}{$-t_{7}$}};
            
            \draw[line width=2pt, ->]  (v45) .. controls (5.25,1.75) and (5.25,2.25) .. (v103)
            node[midway, below left] {\scriptsize {$t_{8}$}};
            
            \draw[draw=dgreen, line width=2pt, ->]  (v53) .. controls (4.75,2.25) and (4.75,2.75) .. (v95)
            node[midway, below left] {\scriptsize \textcolor{dgreen}{$-t_{8}$}}; 
        
            \draw[line width=2pt, ->]  (v35) -- (v46)
            node[midway, right] {\scriptsize {$t_{9}$}};
            
            \draw[draw=dgreen, line width=2pt, ->]  (v96) -- (v85) 
            node[midway, right] {\scriptsize \textcolor{dgreen}{$-t_{9}$}};
        
            \draw[line width=2pt, ->]  (v47) -- (v54)
            node[midway, right] {\scriptsize {$t_{10}$}};
            
            \draw[line width=2pt, ->]  (v104) -- (v97)
            node[midway, right] {\scriptsize {$t_{10}$}};
        \end{tikzpicture}}
    \caption{The unique path collection from $I = \{1,2,3,4,5\}$ to $\{3,4,6,7,8\}$. }
    \label{fig:examplePathCollection2}
\end{figure}
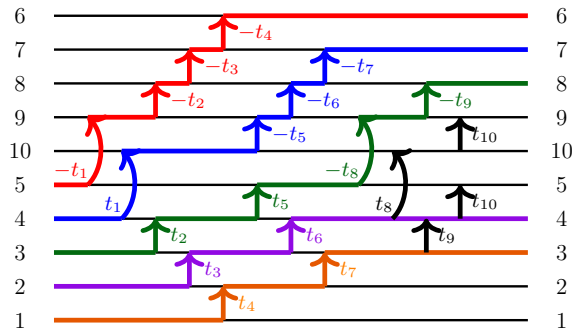

\begin{example}\label{example:n5pathcollection}
    Figure \ref{fig:examplePathCollection2} depicts a non-intersecting path collection in the LGV diagram for $n=5$ with $5$ paths from $I=\{1,2,3,4,5\}$ to $J=\{3,4,8,7,6\}$. The red, blue and green paths from $5$ to $6$, from $4$ to $7$, and from $3$ to $8$ respectively, are examples of left greedy paths. The orange and purple paths from $1$ to $3$ and from $2$ to $4$, respectively, are right greedy in this path collection. 
\end{example}

\begin{proposition}\label{LGVlemma}
    Let $I \in \binom{[2n]}{n}$ and $\mathcal{P}_I$ be the set of non-intersecting path collections from $\{1,\ldots,n\}$ to $I$ in the LGV diagram. Then, the maximal minor of $X$ in the columns indexed by $I$ is given by
    \[
        \Delta^I(X) :=\det\left(X_I\right)= \sum_{P=(P_1,\ldots,P_n) \in \mathcal{P}_I} \sgn(\pi_P)\prod_{i=1}^n \prod_{e \in P_i} w(e),
    \]
    where $\pi_P$ is the permutation that sends $i$ to the end point of $P_i$, and $w(e) \in \{\pm t_1, \dots, \pm t_N\}$ is the weight of the edge $e$.
\end{proposition}
\begin{proof}
    Recall that $X$ is the first $n$ rows of $Z:= x_{i_1}(t_1) \cdots x_{i_N}(t_N)$. Hence, the maximal minors of $X$ are $n\times n$ minors of $Z$ where we always pick the first $n$ rows. 
    Note that we can view each matrix $Y\coloneqq x_{i_l}(t_l)$ as the adjacency matrix for a weighted directed graph $G_l$ with $2n$ source vertices, $2n$ sink vertices and with source $j$ connected to sink $k$ by a directed edge of weight $(Y)_{jk}$ if $(Y)_{jk}\neq 0$. A classical result from graph theory tells us that the entry $(Z)_{jk}$ is the sum of the weights of paths from source $j$ to sink $k$ in the concatenation of the graphs $G_l$, where for $1\leq l\leq N-1$ we identify the sinks of $G_l$ with the corresponding sources of $G_{l+1}$. Up to adding in horizontal strands, which don't contribute to the weights of paths, this concatenation is precisely the LGV diagram. Thus, the result follows from a simple corollary \cite[Corollary 2.21]{boretsky} of the LGV lemma \cite{GV,Lindstrom}.
\end{proof}

In the LGV diagram, each matrix $x_{i_l}$ contributes two arrows. There are two different types of arrows, depending on whether or not $i_l = n$. If $i_l \neq n$ then the arrows corresponding to $x_{i_l}(t_l)$ are $i_l \to i_l +1$ and $n+i_l-1 \to n+i_l$. If $i_l = n$, the arrows are $n-1 \to 2n$ and $n \to 2n-1$, which jump over strands $n$ and $2n$ respectively. Now, we introduce notation to encode each arrow in the graph and their corresponding weights.

\begin{itemize}
    \item We start with arrows that do not skip strands. These will be denoted $a_i^{(j)}$, with the lower index describing which strands of the LGV diagram that the arrow spans and the upper index describing where it lies horizontally. For $i<n$, $a_i^{(j)}$ is the arrow $i \to i+1$ coming from the $s_i$ in the $j$th set of parentheses in expression \eqref{eq:w_0ReducedExpr}. For $i>n$, $a_i^{(j)}$ is the arrow $i \to i-1$ coming from the $s_{n-i-1}$ in the $j$th set of parentheses in expression \eqref{eq:w_0ReducedExpr}.
    
    \item All other arrows skip strands. These are denoted $b_i^{(j)}$, where the upper and lower indices again describe the horizontal and vertical position of the arrow in the LGV diagram, respectively. For $i=n$ or $i=n-1$, $b_i^{(j)}$ is the arrow $n \to 2n-1$  or from $n-1\to 2n$, respectively, coming from the $j^{\text{th}}$ time $s_n$ appears in expression \eqref{eq:w_0ReducedExpr}. 
\end{itemize}

Each of these arrows is weighted by $\pm t_l$ for some $l$. We define a function $\lambda$ which assigns to an arrow $a$ with weight $t_l$ the index of its weight $\lambda(a)=l$. Note that for $i<n$ and any suitable $j$, $\lambda\left(a_i^{(j)}\right) = \lambda\left(a_{n+i+1}^{(j)}\right)$ and $\lambda\left(b_{n-1}^{(j)}\right) = \lambda\left(b_n^{(j)}\right)$.

\begin{example}
    Paths in the LGV diagram, can be explicitly described in terms of the arrows $a_{i}^{(j)}, b_{i}^{(j)}$. For example, the path collection in \Cref{example:n5pathcollection} is given by the five paths:
    \[
    \textcolor{dorange}{a_1^1a_2^2}, \quad \textcolor{dviolet}{a_2^1a_3^2}, \quad \textcolor{dgreen}{a_3^1a_4^2b_5^2a_9^3}, \quad \textcolor{blue}{b_4^1a_{10}^2a_9^2a_8^2}, \quad \text{and} \quad \textcolor{red}{b_5^1a_{9}^1a_8^1a_7^1}.
    \]
    The weight indices of the arrows in this LGV diagram are as follows.
    \smallskip
    \begin{center}
        \begin{tabular}{c|c|c|c|c|c|c|c|c|c|c}
             \textbf{Arrow $\mathbf{a}$} & $b_4^1, b_5^1$ & $a_1^1, a_9^1$ & $a_2^1, a_8^1$ & $a_3^1, a_7^1$ & $a_4^2, a_{10}^2$ & $a_3^2, a_7^2$ & $a_2^2, a_8^2$ & $b_4^2, b_5^2$ & $a_3^3, a_7^3$ & $a_4^4, a_{10}^4$\\
             \hline
             $\lambda(\mathbf{a})$& $1$ & $2$ & $3$ & $4$ & $5$ & $6$ & $7$ & $8$ & $9$ & $10$
        \end{tabular}
    \end{center}
\end{example}

\begin{remark}
    The LGV diagram can be explicitly described using expression \eqref{eq:w_0ReducedExpr}. The LGV diagram consists of $2n$ horizontal strands, with a sequence of arrows added from left to right, one pair corresponding to each $s_i$ in \eqref{eq:w_0ReducedExpr}. (We have been drawing pairs of arrows directly on top of one another, but there is no change in anything we have described if we draw the higher arrow slightly to the right of the lower arrow.) For $n$ even the sequence of arrows is as follows, with parentheses matching those in \eqref{eq:w_0ReducedExpr}:
    \begin{equation}\label{eq:curlyW0}
        \mathcal{W}_0 \coloneqq \begin{array}{l}
             b_n^{(1)}b_{n-1}^{(1)} (a_{n-2}^{(1)}a_{2n-1}^{(1)} \ \cdots \ a_{1}^{(1)}a_{n+2}^{(1)})(a_{n-1}^{(2)}a_{2n}^{(2)} \ \cdots \ a_{2}^{(2)}a_{n+3}^{(2)})\\
             b_n^{(2)}b_{n-1}^{(2)}(a_{n-2}^{(3)}a_{2n-1}^{(3)} \ \cdots \ a_{3}^{(3)}a_{n+4}^{(3)})(a_{n-1}^{(4)}a_{2n}^{(4)}\ \cdots \ a_{4}^{(4)}a_{n+5}^{(4)}) \\
             \multicolumn{1}{c}{\vdots}\\
             b_n^{(l)}b_{n-1}^{(l)} (a_{n-2}^{(2l-1)}a_{2n-1}^{(2l-1)} \ \cdots \ a_{2l-1}^{(2l-1)}a_{n+2l}^{(2l-1)}) (a_{n-1}^{(2l)}a_{2n}^{(2l)} \ \cdots  \ a_{2l}^{(2l)}a_{n+2l+1}^{(2l)})\\
             \multicolumn{1}{c}{\vdots}\\
             b_n^{(n/2)}b_{n-1}^{(n/2)}.
        \end{array}
    \end{equation}
   
For $n$ odd, we replace the last line with $b_n^{(\lfloor n/2\rfloor)}b_{n-1}^{(\lfloor n/2\rfloor)}(a_{n-2}^{(n-2)}a_{2n-1}^{(n-2)}) (a_{n-1}^{(n-1)} a_{2n}^{(n-1)})$.

\medskip

Note that $\mathcal{W}_0$ contains either $(n-2)$ or $(n-1)$ sets of parentheses, for $n$ even or odd respectively. Each set of parentheses defines a ``diagonal lines'' in the LGV diagram. For instance, see the arrows $\textcolor{dgreen}{a_3^{(1)}}\textcolor{dviolet}{a_2^{(1)}}\textcolor{dorange}{a_1^{(1)}}$ and $\textcolor{red}{a_9^{(1)}}\textcolor{red}{a_8^{(1)}}\textcolor{red}{a_7^{(1)}}$ in \Cref{fig:examplePathCollection2}, which all appear in the first set of parentheses of $\mathcal{W}_0$. We will refer to the $m^{\text{th}}$ pair of parentheses in $\mathcal{W}_0$ as \textit{parentheses pair $m$}. Explicitly, parentheses pair $m$ contains all arrows $a_i^{(2m-1)}$ and $a_i^{(2m)}$.
\end{remark}

\section{Positivity criterion} \label{sec:3}

\subsection{Monomial Minors}

The minors in \Cref{def:SpecialMinorsPfaff} correspond to particular path collections in the LGV diagram. First, note that the maximal minors of $X = X(t_1,\ldots, t_N)$ are, up to sign, equal to the minors of $A = A(t_1,\ldots,t_N)$. We defer the specific sign analysis of $M_{j,k}$ to \Cref{subsec:signAnalysis}, but we have the following: For $1\leq j\leq k\leq n-1$,
\[
    M_{j,k} := M_{j,k}(A) = \pm \ \Delta^{n-k,\ldots,n-k+j-1,n+1,n+2,\ldots,2n-j}(X), \quad \text{ for } 1 \leq j \leq k \leq n-1.
\] 
By virtue of \Cref{LGVlemma}, to compute $M_{j,k}$, we must look at non-intersecting path collections from $\{1,\ldots, n\}$ to $\{n-k,\ldots,n-k+j-1,n+1,n+2,\ldots,2n-j\}$.

\begin{lemma}\label{lem:LGVdiagramLeftGreedy}
    The left greedy path collection from source vertices $[n]$ in the LGV diagram is a unique non-intersecting path collection that uses all arrows in the LGV diagram. Moreover, if ${\rm LGP}(i)$ terminates at sink vertex $j$, then no other path originating from source vertex $i$ terminates above strand $j$. 
\end{lemma}
\begin{proof}
    We will prove these statements in more generality in \Cref{prop:pathcollectionproperties} and \Cref{lem:greedyisextremeuv}.
\end{proof}

\begin{proposition}\label{prop:UniquePathColl}
    For any $1 \leq j \leq k \leq n-1$, there is a unique non-intersecting path collection from the source vertices $\{1,\ldots, n\}$ to the sink vertices $\{n-k,\ldots,n-k+j-1,n+1,n+2,\ldots,2n-j\}$. Explicitly it is given by the following $n$ paths:
    \[
    1 \to n-k, \quad \cdots, \quad j \to n-k+j-1, \quad j+1 \to 2n-j, \quad \cdots, \quad n \to n+1.
    \]
    In particular, the minors $M_{j,k}$ are monomials in the variables $t_i$.
\end{proposition}

\smallskip

Before giving a proof, let us first give some intuition on the path collections in \Cref{prop:UniquePathColl}. For $1\leq j \leq k \leq n-1$, the corresponding path collection can be obtained as follows:
\begin{itemize}
    \item Start with the path collection consisting of the path from $1$ to $k$ and the left greedy paths starting at $2,\ldots, n$.
    \item Modify sequentially the paths starting at $2,\ldots, n-k+j-1$ and make them right greedy. That is, make path starting at $2$ right greedy, and then the one starting at $3$, and so on.
\end{itemize}
Then, the index $k$ indicates where the path starting at $1$ ends while the index $j$ tells us which paths are right greedy and which are left greedy.

\begin{example}\label{ex:somePathsCollectionsInLGVDiagram}
    In the case of $n=4,5$ there are respectively $6$ and $10$ minors $M_{j,k}$. We illustrate some examples of how $M_{1,1}$, $M_{1,2}$ and $M_{2,2}$ arise from unique non-intersecting path collections. For $k=j=1$ we obtain the following path collections which corresponds to the minor $M_{1,1}=(t_1\cdots t_5)^2t_{6}$ for $n=4$ and $M_{1,1}=(t_1\cdots t_9)^2t_{10}$ for $n=5$.
   
    \begin{center}
    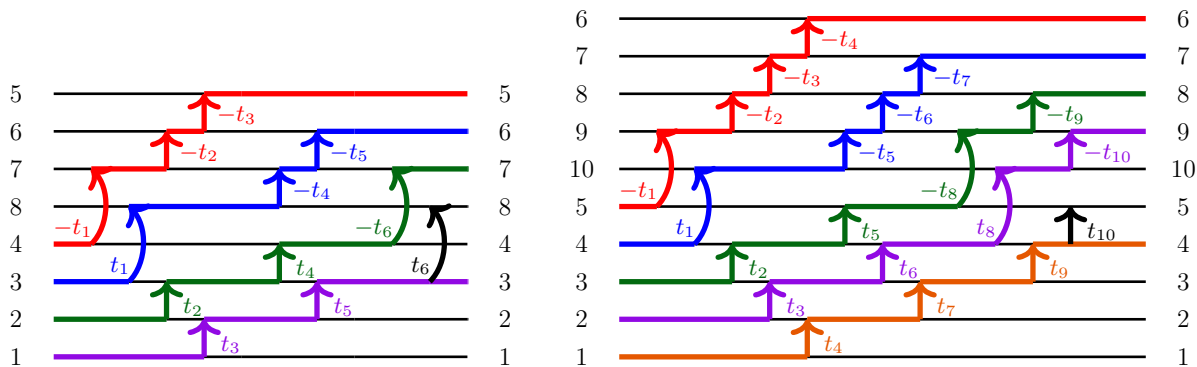
\begin{figure}[ht]
    \begin{tikzpicture}
    
    \coordinate (l7) at (0,4);
    \coordinate (l8) at (0,3.5);
    \coordinate (l9) at (0,3);
    \coordinate (l10) at (0,2.5);        
    \coordinate (l5) at (0,2);        
    \coordinate (l4) at (0,1.5);    
    \coordinate (l3) at (0,1);
    \coordinate (l2) at (0,0.5);

    \coordinate (r7) at (5.5,4);
    \coordinate (r8) at (5.5,3.5);
    \coordinate (r9) at (5.5,3);
    \coordinate (r10) at (5.5,2.5);
    \coordinate (r5) at (5.5,2);        
    \coordinate (r4) at (5.5,1.5);    
    \coordinate (r3) at (5.5,1);
    \coordinate (r2) at (5.5,0.5);

    \coordinate (v21) at (2,0.5);
    \coordinate (v22) at (2.5,0.5);
    \coordinate (v23) at (4,0.5);

    \coordinate (v31) at (1.5,1);
    \coordinate (v32) at (2,1);
    \coordinate (v33) at (3.5,1);
    \coordinate (v34) at (4,1);
    \coordinate (v35) at (5.5,1);
    
    \coordinate (v41) at (1,1.5);
    \coordinate (v42) at (1.5,1.5);
    \coordinate (v43) at (3,1.5);
    \coordinate (v44) at (3.5,1.5);
    \coordinate (v45) at (5,1.5);
    \coordinate (v46) at (5.5,1.5);
    \coordinate (v47) at (5.5,1.5);
    
    \coordinate (v51) at (0.5,2);
    \coordinate (v52) at (3,2);
    \coordinate (v53) at (4.5,2);
    \coordinate (v54) at (5.5,2);

    \coordinate (v101) at (1,2.5);    
    \coordinate (v102) at (3,2.5);
    \coordinate (v103) at (5,2.5);
    \coordinate (v104) at (5.5,2.5);

    \coordinate (v91) at (0.5,3);
    \coordinate (v92) at (1.5,3);
    \coordinate (v93) at (3,3);
    \coordinate (v94) at (3.5,3);
    \coordinate (v95) at (4.5,3);
    \coordinate (v96) at (5.5,3);
    \coordinate (v97) at (5.5,3);

    \coordinate (v81) at (1.5,3.5);
    \coordinate (v82) at (2,3.5);
    \coordinate (v83) at (3.5,3.5);
    \coordinate (v84) at (4,3.5);
    \coordinate (v85) at (5.5,3.5);

    \coordinate (v71) at (2,4);
    \coordinate (v72) at (2.5,4);
    \coordinate (v73) at (4,4);

    \draw[black] (-0.5,0.5) node  [xscale = 0.8, yscale = 0.8] {$1$};
    \draw[black] (-0.5,1) node  [xscale = 0.8, yscale = 0.8] {$2$};
    \draw[black] (-0.5,1.5) node  [xscale = 0.8, yscale = 0.8] {$3$};
    \draw[black] (-0.5,2) node  [xscale = 0.8, yscale = 0.8] {$4$};
    \draw[black] (-0.5,2.5) node  [xscale = 0.8, yscale = 0.8] {$8$};
    \draw[black] (-0.5,3) node  [xscale = 0.8, yscale = 0.8] {$7$};
    \draw[black] (-0.5,3.5) node  [xscale = 0.8, yscale = 0.8] {$6$};
    \draw[black] (-0.5,4) node  [xscale = 0.8, yscale = 0.8] {$5$};

    \draw[black] (6,0.5) node  [xscale = 0.8, yscale = 0.8] {$1$};
    \draw[black] (6,1) node  [xscale = 0.8, yscale = 0.8] {$2$};
    \draw[black] (6,1.5) node  [xscale = 0.8, yscale = 0.8] {$3$};
    \draw[black] (6,2) node  [xscale = 0.8, yscale = 0.8] {$4$};
    \draw[black] (6,2.5) node  [xscale = 0.8, yscale = 0.8] {$8$};
    \draw[black] (6,3) node  [xscale = 0.8, yscale = 0.8] {$7$};
    \draw[black] (6,3.5) node  [xscale = 0.8, yscale = 0.8] {$6$};
    \draw[black] (6,4) node  [xscale = 0.8, yscale = 0.8] {$5$};

    \draw[draw=dviolet, line width=2pt] (l2) -- (v21);
    \draw[draw=black, line width=1pt] (v21) -- (v22);
    \draw[draw=black, line width=1pt] (v22) -- (v23);
    \draw[draw=black, line width=1pt] (v23) -- (r2);

    \draw[draw=dgreen, line width=2pt] (l3) -- (v31);
    \draw[draw=black, line width=1pt] (v31) -- (v32);
    \draw[draw=dviolet, line width=2pt] (v32) -- (v33);
    \draw[draw=black, line width=1pt] (v33) -- (v34);
    \draw[draw=black, line width=1pt] (v34) -- (v35);
    \draw[draw=black, line width=1pt] (v35) -- (r3);

    \draw[draw=blue, line width=2pt] (l4) -- (v41);
    \draw[draw=black, line width=1pt] (v41) -- (v42);
    \draw[draw=dgreen, line width=2pt] (v42) -- (v43);
    \draw[draw=black, line width=1pt] (v43) -- (v44);
    \draw[draw=dviolet, line width=2pt] (v44) -- (v45);
    \draw[draw=dviolet, line width=2pt] (v45) -- (v46);
    \draw[draw=dviolet, line width=2pt] (v46) -- (r4);

    \draw[draw=red, line width=2pt] (l5) -- (v51);
    \draw[draw=black, line width=1pt] (v51) -- (v52);
    \draw[draw=dgreen, line width=2pt] (v52) -- (v53);
    \draw[draw=black, line width=1pt] (v53) -- (r5);

    \draw[draw=black, line width=1pt] (l10) -- (v101);
    \draw[draw=blue, line width=2pt] (v101) -- (v102);
    \draw[draw=black, line width=1pt] (v102) -- (v103);
    \draw[draw=black, line width=1pt] (v103) -- (v104);
    \draw[draw=black, line width=1pt] (v104) -- (r10);

    \draw[draw=black, line width=1pt] (l9) -- (v91);
    \draw[draw=red, line width=2pt] (v91) -- (v92);
    \draw[draw=black, line width=1pt] (v92) -- (v93);
    \draw[draw=blue, line width=2pt] (v93) -- (v94);
    \draw[draw=black, line width=1pt] (v94) -- (v95);
    \draw[draw=dgreen, line width=2pt] (v95) -- (v96);
    \draw[draw=dgreen, line width=2pt] (v96) -- (v97);
    \draw[draw=dgreen, line width=2pt] (v97) -- (r9);

    \draw[draw=black, line width=1pt] (l8) -- (v81);
    \draw[draw=red, line width=2pt] (v81) -- (v82);
    \draw[draw=black, line width=1pt] (v82) -- (v83);
    \draw[draw=blue, line width=2pt] (v83) -- (v84);
    \draw[draw=blue, line width=2pt] (v84) -- (v85);
    \draw[draw=blue, line width=2pt] (v85) -- (r8);

    \draw[draw=black, line width=1pt] (l7) -- (v71);
    \draw[draw=red, line width=2pt] (v71) -- (v72);
    \draw[draw=red, line width=2pt] (v72) -- (v73);
    \draw[draw=red, line width=2pt] (v73) -- (r7);

    \draw[draw=blue, line width=2pt, ->]  (v41) .. controls (1.25,1.75) and (1.25,2.25) .. (v101) node[midway, below left] {\scriptsize \textcolor{blue}{$t_{1}$}};
            
    \draw[draw=red, line width=2pt, ->]  (v51) .. controls (0.75,2.25) and (0.75,2.75) .. (v91) node[midway, below left] {\scriptsize \textcolor{red}{$-t_{1}$}};
    
    \draw[draw=dgreen, line width=2pt, ->]  (v31) -- (v42) node[below right=0.5mm] {\scriptsize \textcolor{dgreen}{$t_{2}$}};
    \draw[draw=red, line width=2pt, ->]  (v92) -- (v81) node[below right] {\scriptsize \textcolor{red}{$-t_{2}$}};
    
    \draw[draw=dviolet, line width=2pt, ->]  (v21) -- (v32) node[below right=0.5mm] {\scriptsize \textcolor{dviolet}{$t_{3}$}};
    \draw[draw=red, line width=2pt, ->]  (v82) -- (v71) node[below right] {\scriptsize \textcolor{red}{$-t_{3}$}};

    \draw[draw=dgreen, line width=2pt, ->]  (v43) -- (v52) node[below right=0.5mm] {\scriptsize \textcolor{dgreen}{$t_{4}$}};
    \draw[draw=blue, line width=2pt, ->]  (v102) -- (v93) node[below right] {\scriptsize \textcolor{blue}{$-t_{4}$}};

    \draw[draw=dviolet, line width=2pt, ->]  (v33) -- (v44) node[below right=0.5mm] {\scriptsize \textcolor{dviolet}{$t_5$}};
    \draw[draw=blue, line width=2pt, ->]  (v94) -- (v83) node[below right] {\scriptsize \textcolor{blue}{$-t_{5}$}};

    \draw[draw=black, line width=2pt, ->]  (v45) .. controls (5.25,1.75) and (5.25,2.25) .. (v103)
            node[midway, below left] {\scriptsize $t_{6}$};
            
    \draw[draw=dgreen, line width=2pt, ->]  (v53) .. controls (4.75,2.25) and (4.75,2.75) .. (v95)
            node[midway, below left] {\scriptsize \textcolor{dgreen}{$-t_{6}$}};

\end{tikzpicture}
\quad
\begin{tikzpicture}
    \coordinate (l6) at (0,4.5);    
    \coordinate (l7) at (0,4);
    \coordinate (l8) at (0,3.5);
    \coordinate (l9) at (0,3);
    \coordinate (l10) at (0,2.5);        
    \coordinate (l5) at (0,2);        
    \coordinate (l4) at (0,1.5);    
    \coordinate (l3) at (0,1);
    \coordinate (l2) at (0,0.5);
    \coordinate (l1) at (0,0);

    \coordinate (r6) at (7,4.5);    
    \coordinate (r7) at (7,4);
    \coordinate (r8) at (7,3.5);
    \coordinate (r9) at (7,3);
    \coordinate (r10) at (7,2.5);
    \coordinate (r5) at (7,2);        
    \coordinate (r4) at (7,1.5);    
    \coordinate (r3) at (7,1);
    \coordinate (r2) at (7,0.5);
    \coordinate (r1) at (7,0);

    \coordinate (v11) at (2.5,0);
    
    \coordinate (v21) at (2,0.5);
    \coordinate (v22) at (2.5,0.5);
    \coordinate (v23) at (4,0.5);

    \coordinate (v31) at (1.5,1);
    \coordinate (v32) at (2,1);
    \coordinate (v33) at (3.5,1);
    \coordinate (v34) at (4,1);
    \coordinate (v35) at (5.5,1);
    
    \coordinate (v41) at (1,1.5);
    \coordinate (v42) at (1.5,1.5);
    \coordinate (v43) at (3,1.5);
    \coordinate (v44) at (3.5,1.5);
    \coordinate (v45) at (5,1.5);
    \coordinate (v46) at (5.5,1.5);
    \coordinate (v47) at (6,1.5);
    
    \coordinate (v51) at (0.5,2);
    \coordinate (v52) at (3,2);
    \coordinate (v53) at (4.5,2);
    \coordinate (v54) at (6,2);

    \coordinate (v101) at (1,2.5);    
    \coordinate (v102) at (3,2.5);
    \coordinate (v103) at (5,2.5);
    \coordinate (v104) at (6,2.5);

    \coordinate (v91) at (0.5,3);
    \coordinate (v92) at (1.5,3);
    \coordinate (v93) at (3,3);
    \coordinate (v94) at (3.5,3);
    \coordinate (v95) at (4.5,3);
    \coordinate (v96) at (5.5,3);
    \coordinate (v97) at (6,3);

    \coordinate (v81) at (1.5,3.5);
    \coordinate (v82) at (2,3.5);
    \coordinate (v83) at (3.5,3.5);
    \coordinate (v84) at (4,3.5);
    \coordinate (v85) at (5.5,3.5);

    \coordinate (v71) at (2,4);
    \coordinate (v72) at (2.5,4);
    \coordinate (v73) at (4,4);

    \coordinate (v61) at (2.5,4.5);
    
    \draw[black] (-0.5,0) node  [xscale = 0.8, yscale = 0.8] {$1$};
    \draw[black] (-0.5,0.5) node  [xscale = 0.8, yscale = 0.8] {$2$};
    \draw[black] (-0.5,1) node  [xscale = 0.8, yscale = 0.8] {$3$};
    \draw[black] (-0.5,1.5) node  [xscale = 0.8, yscale = 0.8] {$4$};
    \draw[black] (-0.5,2) node  [xscale = 0.8, yscale = 0.8] {$5$};
    \draw[black] (-0.5,2.5) node  [xscale = 0.8, yscale = 0.8] {$10$};
    \draw[black] (-0.5,3) node  [xscale = 0.8, yscale = 0.8] {$9$};
    \draw[black] (-0.5,3.5) node  [xscale = 0.8, yscale = 0.8] {$8$};
    \draw[black] (-0.5,4) node  [xscale = 0.8, yscale = 0.8] {$7$};
    \draw[black] (-0.5,4.5) node  [xscale = 0.8, yscale = 0.8] {$6$};
    
    \draw[black] (7.5,0) node  [xscale = 0.8, yscale = 0.8] {$1$};
    \draw[black] (7.5,0.5) node  [xscale = 0.8, yscale = 0.8] {$2$};
    \draw[black] (7.5,1) node  [xscale = 0.8, yscale = 0.8] {$3$};
    \draw[black] (7.5,1.5) node  [xscale = 0.8, yscale = 0.8] {$4$};
    \draw[black] (7.5,2) node  [xscale = 0.8, yscale = 0.8] {$5$};
    \draw[black] (7.5,2.5) node  [xscale = 0.8, yscale = 0.8] {$10$};
    \draw[black] (7.5,3) node  [xscale = 0.8, yscale = 0.8] {$9$};
    \draw[black] (7.5,3.5) node  [xscale = 0.8, yscale = 0.8] {$8$};
    \draw[black] (7.5,4) node  [xscale = 0.8, yscale = 0.8] {$7$};
    \draw[black] (7.5,4.5) node  [xscale = 0.8, yscale = 0.8] {$6$};

    \draw[draw=dorange, line width=2pt] (l1)  -- (v11);
    \draw[draw=black, line width=1pt] (v11) -- (r1);

    \draw[draw=dviolet, line width=2pt] (l2) -- (v21);
    \draw[draw=black, line width=1pt] (v21) -- (v22);
    \draw[draw=dorange, line width=2pt] (v22) -- (v23);
    \draw[draw=black, line width=1pt] (v23) -- (r2);

    \draw[draw=dgreen, line width=2pt] (l3) -- (v31);
    \draw[draw=black, line width=1pt] (v31) -- (v32);
    \draw[draw=dviolet, line width=2pt] (v32) -- (v33);
    \draw[draw=black, line width=1pt] (v33) -- (v34);
    \draw[draw=dorange, line width=2pt] (v34) -- (v35);
    \draw[draw=black, line width=1pt] (v35) -- (r3);

    \draw[draw=blue, line width=2pt] (l4) -- (v41);
    \draw[draw=black, line width=1pt] (v41) -- (v42);
    \draw[draw=dgreen, line width=2pt] (v42) -- (v43);
    \draw[draw=black, line width=1pt] (v43) -- (v44);
    \draw[draw=dviolet, line width=2pt] (v44) -- (v45);
    \draw[draw=black, line width=1pt] (v45) -- (v46);
    \draw[draw=dorange, line width=2pt] (v46) -- (r4);

    \draw[draw=red, line width=2pt] (l5) -- (v51);
    \draw[draw=black, line width=1pt] (v51) -- (v52);
    \draw[draw=dgreen, line width=2pt] (v52) -- (v53);
    \draw[draw=black, line width=1pt] (v53) -- (r5);

    \draw[draw=black, line width=1pt] (l10) -- (v101);
    \draw[draw=blue, line width=2pt] (v101) -- (v102);
    \draw[draw=black, line width=1pt] (v102) -- (v103);
    \draw[draw=dviolet, line width=2pt] (v103) -- (v104);
    \draw[draw=black, line width=1pt] (v104) -- (r10);

    \draw[draw=black, line width=1pt] (l9) -- (v91);
    \draw[draw=red, line width=2pt] (v91) -- (v92);
    \draw[draw=black, line width=1pt] (v92) -- (v93);
    \draw[draw=blue, line width=2pt] (v93) -- (v94);
    \draw[draw=black, line width=1pt] (v94) -- (v95);
    \draw[draw=dgreen, line width=2pt] (v95) -- (v96);
    \draw[draw=black, line width=1pt] (v96) -- (v97);
    \draw[draw=dviolet, line width=2pt] (v97) -- (r9);

    \draw[draw=black, line width=1pt] (l8) -- (v81);
    \draw[draw=red, line width=2pt] (v81) -- (v82);
    \draw[draw=black, line width=1pt] (v82) -- (v83);
    \draw[draw=blue, line width=2pt] (v83) -- (v84);
    \draw[draw=black, line width=1pt] (v84) -- (v85);
    \draw[draw=dgreen, line width=2pt] (v85) -- (r8);

    \draw[draw=black, line width=1pt] (l7) -- (v71);
    \draw[draw=red, line width=2pt] (v71) -- (v72);
    \draw[draw=black, line width=1pt] (v72) -- (v73);
    \draw[draw=blue, line width=2pt] (v73) -- (r7);

    \draw[draw=black, line width=1pt] (l6) -- (v61);
    \draw[draw=red, line width=2pt] (v61) -- (r6);

    \draw[draw=blue, line width=2pt, ->]  (v41) .. controls (1.25,1.75) and (1.25,2.25) .. (v101) node[midway, below left] {\scriptsize \textcolor{blue}{$t_{1}$}};
            
    \draw[draw=red, line width=2pt, ->]  (v51) .. controls (0.75,2.25) and (0.75,2.75) .. (v91) node[midway, below left] {\scriptsize \textcolor{red}{$-t_{1}$}};
    
    \draw[draw=dgreen, line width=2pt, ->]  (v31) -- (v42) node[below right=0.5mm] {\scriptsize \textcolor{dgreen}{$t_{2}$}};
    \draw[draw=red, line width=2pt, ->]  (v92) -- (v81) node[below right] {\scriptsize \textcolor{red}{$-t_{2}$}};
    
    \draw[draw=dviolet, line width=2pt, ->]  (v21) -- (v32) node[below right=0.5mm] {\scriptsize \textcolor{dviolet}{$t_{3}$}};
    \draw[draw=red, line width=2pt, ->]  (v82) -- (v71) node[below right] {\scriptsize \textcolor{red}{$-t_{3}$}};
    
    \draw[draw=dorange, line width=2pt, ->]  (v11) -- (v22) node[below right=0.5mm] {\scriptsize \textcolor{dorange}{$t_{4}$}};
    \draw[draw=red, line width=2pt, ->]  (v72) -- (v61) node[below right] {\scriptsize \textcolor{red}{$-t_{4}$}};

    \draw[draw=dgreen, line width=2pt, ->]  (v43) -- (v52) node[below right=0.5mm] {\scriptsize \textcolor{dgreen}{$t_{5}$}};
    \draw[draw=blue, line width=2pt, ->]  (v102) -- (v93) node[below right] {\scriptsize \textcolor{blue}{$-t_{5}$}};

    \draw[draw=dviolet, line width=2pt, ->]  (v33) -- (v44) node[below right=0.5mm] {\scriptsize \textcolor{dviolet}{$t_{6}$}};
    \draw[draw=blue, line width=2pt, ->]  (v94) -- (v83) node[below right] {\scriptsize \textcolor{blue}{$-t_{6}$}};

    \draw[draw=dorange, line width=2pt, ->]  (v23) -- (v34) node[below right=0.5mm] {\scriptsize \textcolor{dorange}{$t_{7}$}};
    \draw[draw=blue, line width=2pt, ->]  (v84) -- (v73) node[below right] {\scriptsize \textcolor{blue}{$-t_{7}$}};

    \draw[draw=dviolet, line width=2pt, ->]  (v45) .. controls (5.25,1.75) and (5.25,2.25) .. (v103)
            node[midway, below left] {\scriptsize \textcolor{dviolet}{$t_{8}$}};
            
    \draw[draw=dgreen, line width=2pt, ->]  (v53) .. controls (4.75,2.25) and (4.75,2.75) .. (v95)
            node[midway, below left] {\scriptsize \textcolor{dgreen}{$-t_{8}$}}; 
        
    \draw[draw=dorange, line width=2pt, ->]  (v35) -- (v46) node[below right=0.5mm] {\scriptsize \textcolor{dorange}{$t_{9}$}};
    \draw[draw=dgreen, line width=2pt, ->]  (v96) -- (v85) node[below right] {\scriptsize \textcolor{dgreen}{$-t_{9}$}};

    \draw[line width=2pt, ->]  (v47) -- (v54)node[below right=0.5mm] {\scriptsize  $t_{10}$};
    \draw[draw=dviolet, line width=2pt, ->]  (v104) -- (v97) node[below right] {\scriptsize \textcolor{dviolet}{$-t_{10}$}};
\end{tikzpicture}
\caption{Illustration of the path collection corresponding to $M_{1,1}$.}
\end{figure}
\vspace{-1.5em}
\end{center}
For $k=2,j=1$ we modify the path starting at $1$ and obtain the following path collections corresponding to $M_{1,2}=(t_1t_2t_3t_4)^2t_5t_{6}$ for $n=4$ and $M_{1,2}=(t_1\cdots t_8)^2t_9t_{10}$ for $n=5$.
    \begin{center}
    \begin{figure}[ht]
    \begin{tikzpicture}
    
    \coordinate (l7) at (0,4);
    \coordinate (l8) at (0,3.5);
    \coordinate (l9) at (0,3);
    \coordinate (l10) at (0,2.5);        
    \coordinate (l5) at (0,2);        
    \coordinate (l4) at (0,1.5);    
    \coordinate (l3) at (0,1);
    \coordinate (l2) at (0,0.5);

    \coordinate (r7) at (5.5,4);
    \coordinate (r8) at (5.5,3.5);
    \coordinate (r9) at (5.5,3);
    \coordinate (r10) at (5.5,2.5);
    \coordinate (r5) at (5.5,2);        
    \coordinate (r4) at (5.5,1.5);    
    \coordinate (r3) at (5.5,1);
    \coordinate (r2) at (5.5,0.5);

    \coordinate (v21) at (2,0.5);
    \coordinate (v22) at (2.5,0.5);
    \coordinate (v23) at (4,0.5);

    \coordinate (v31) at (1.5,1);
    \coordinate (v32) at (2,1);
    \coordinate (v33) at (3.5,1);
    \coordinate (v34) at (4,1);
    \coordinate (v35) at (5.5,1);
    
    \coordinate (v41) at (1,1.5);
    \coordinate (v42) at (1.5,1.5);
    \coordinate (v43) at (3,1.5);
    \coordinate (v44) at (3.5,1.5);
    \coordinate (v45) at (5,1.5);
    \coordinate (v46) at (5.5,1.5);
    \coordinate (v47) at (5.5,1.5);
    
    \coordinate (v51) at (0.5,2);
    \coordinate (v52) at (3,2);
    \coordinate (v53) at (4.5,2);
    \coordinate (v54) at (5.5,2);

    \coordinate (v101) at (1,2.5);    
    \coordinate (v102) at (3,2.5);
    \coordinate (v103) at (5,2.5);
    \coordinate (v104) at (5.5,2.5);

    \coordinate (v91) at (0.5,3);
    \coordinate (v92) at (1.5,3);
    \coordinate (v93) at (3,3);
    \coordinate (v94) at (3.5,3);
    \coordinate (v95) at (4.5,3);
    \coordinate (v96) at (5.5,3);
    \coordinate (v97) at (5.5,3);

    \coordinate (v81) at (1.5,3.5);
    \coordinate (v82) at (2,3.5);
    \coordinate (v83) at (3.5,3.5);
    \coordinate (v84) at (4,3.5);
    \coordinate (v85) at (5.5,3.5);

    \coordinate (v71) at (2,4);
    \coordinate (v72) at (2.5,4);
    \coordinate (v73) at (4,4);

    \draw[black] (-0.5,0.5) node  [xscale = 0.8, yscale = 0.8] {$1$};
    \draw[black] (-0.5,1) node  [xscale = 0.8, yscale = 0.8] {$2$};
    \draw[black] (-0.5,1.5) node  [xscale = 0.8, yscale = 0.8] {$3$};
    \draw[black] (-0.5,2) node  [xscale = 0.8, yscale = 0.8] {$4$};
    \draw[black] (-0.5,2.5) node  [xscale = 0.8, yscale = 0.8] {$8$};
    \draw[black] (-0.5,3) node  [xscale = 0.8, yscale = 0.8] {$7$};
    \draw[black] (-0.5,3.5) node  [xscale = 0.8, yscale = 0.8] {$6$};
    \draw[black] (-0.5,4) node  [xscale = 0.8, yscale = 0.8] {$5$};

    \draw[black] (6,0.5) node  [xscale = 0.8, yscale = 0.8] {$1$};
    \draw[black] (6,1) node  [xscale = 0.8, yscale = 0.8] {$2$};
    \draw[black] (6,1.5) node  [xscale = 0.8, yscale = 0.8] {$3$};
    \draw[black] (6,2) node  [xscale = 0.8, yscale = 0.8] {$4$};
    \draw[black] (6,2.5) node  [xscale = 0.8, yscale = 0.8] {$8$};
    \draw[black] (6,3) node  [xscale = 0.8, yscale = 0.8] {$7$};
    \draw[black] (6,3.5) node  [xscale = 0.8, yscale = 0.8] {$6$};
    \draw[black] (6,4) node  [xscale = 0.8, yscale = 0.8] {$5$};

    \draw[draw=dviolet, line width=2pt] (l2) -- (v21);
    \draw[draw=black, line width=1pt] (v21) -- (v22);
    \draw[draw=black, line width=1pt] (v22) -- (v23);
    \draw[draw=black, line width=1pt] (v23) -- (r2);

    \draw[draw=dgreen, line width=2pt] (l3) -- (v31);
    \draw[draw=black, line width=1pt] (v31) -- (v32);
    \draw[draw=dviolet, line width=2pt] (v32) -- (r3);

    \draw[draw=blue, line width=2pt] (l4) -- (v41);
    \draw[draw=black, line width=1pt] (v41) -- (v42);
    \draw[draw=dgreen, line width=2pt] (v42) -- (v43);
    \draw[draw=black, line width=1pt] (v43) -- (v44);
    \draw[draw=black, line width=1pt] (v44) -- (v45);
    \draw[draw=black, line width=1pt] (v45) -- (v46);
    \draw[line width=2pt] (v46) -- (r4);

    \draw[draw=red, line width=2pt] (l5) -- (v51);
    \draw[draw=black, line width=1pt] (v51) -- (v52);
    \draw[draw=dgreen, line width=2pt] (v52) -- (v53);
    \draw[draw=black, line width=1pt] (v53) -- (r5);

    \draw[draw=black, line width=1pt] (l10) -- (v101);
    \draw[draw=blue, line width=2pt] (v101) -- (v102);
    \draw[draw=black, line width=1pt] (v102) -- (v103);
    \draw[draw=black, line width=1pt] (v103) -- (v104);
    \draw[draw=black, line width=1pt] (v104) -- (r10);

    \draw[draw=black, line width=1pt] (l9) -- (v91);
    \draw[draw=red, line width=2pt] (v91) -- (v92);
    \draw[draw=black, line width=1pt] (v92) -- (v93);
    \draw[draw=blue, line width=2pt] (v93) -- (v94);
    \draw[draw=black, line width=1pt] (v94) -- (v95);
    \draw[draw=dgreen, line width=2pt] (v95) -- (v96);
    \draw[draw=dgreen, line width=2pt] (v96) -- (v97);
    \draw[draw=dgreen, line width=2pt] (v97) -- (r9);

    \draw[draw=black, line width=1pt] (l8) -- (v81);
    \draw[draw=red, line width=2pt] (v81) -- (v82);
    \draw[draw=black, line width=1pt] (v82) -- (v83);
    \draw[draw=blue, line width=2pt] (v83) -- (v84);
    \draw[draw=blue, line width=2pt] (v84) -- (v85);
    \draw[draw=blue, line width=2pt] (v85) -- (r8);

    \draw[draw=black, line width=1pt] (l7) -- (v71);
    \draw[draw=red, line width=2pt] (v71) -- (v72);
    \draw[draw=red, line width=2pt] (v72) -- (v73);
    \draw[draw=red, line width=2pt] (v73) -- (r7);

    \draw[draw=blue, line width=2pt, ->]  (v41) .. controls (1.25,1.75) and (1.25,2.25) .. (v101) node[midway, below left] {\scriptsize \textcolor{blue}{$t_{1}$}};
            
    \draw[draw=red, line width=2pt, ->]  (v51) .. controls (0.75,2.25) and (0.75,2.75) .. (v91) node[midway, below left] {\scriptsize \textcolor{red}{$-t_{1}$}};
    
    \draw[draw=dgreen, line width=2pt, ->]  (v31) -- (v42) node[below right=0.5mm] {\scriptsize \textcolor{dgreen}{$t_{2}$}};
    \draw[draw=red, line width=2pt, ->]  (v92) -- (v81) node[below right] {\scriptsize \textcolor{red}{$-t_{2}$}};
    
    \draw[draw=dviolet, line width=2pt, ->]  (v21) -- (v32) node[below right=0.5mm] {\scriptsize \textcolor{dviolet}{$t_{3}$}};
    \draw[draw=red, line width=2pt, ->]  (v82) -- (v71) node[below right] {\scriptsize \textcolor{red}{$-t_{3}$}};

    \draw[draw=dgreen, line width=2pt, ->]  (v43) -- (v52) node[below right=0.5mm] {\scriptsize \textcolor{dgreen}{$t_{4}$}};
    \draw[draw=blue, line width=2pt, ->]  (v102) -- (v93) node[below right] {\scriptsize \textcolor{blue}{$-t_{4}$}};

    \draw[line width=2pt, ->]  (v33) -- (v44) node[below right=0.5mm] {\scriptsize $t_5$};
    \draw[draw=blue, line width=2pt, ->]  (v94) -- (v83) node[below right] {\scriptsize \textcolor{blue}{$-t_{5}$}};

    \draw[draw=black, line width=2pt, ->]  (v45) .. controls (5.25,1.75) and (5.25,2.25) .. (v103)
            node[midway, below left] {\scriptsize $t_{6}$};
            
    \draw[draw=dgreen, line width=2pt, ->]  (v53) .. controls (4.75,2.25) and (4.75,2.75) .. (v95)
            node[midway, below left] {\scriptsize \textcolor{dgreen}{$-t_{6}$}};

\end{tikzpicture}
\quad
        \begin{tikzpicture}
    \coordinate (l6) at (0,4.5);    
    \coordinate (l7) at (0,4);
    \coordinate (l8) at (0,3.5);
    \coordinate (l9) at (0,3);
    \coordinate (l10) at (0,2.5);        
    \coordinate (l5) at (0,2);        
    \coordinate (l4) at (0,1.5);    
    \coordinate (l3) at (0,1);
    \coordinate (l2) at (0,0.5);
    \coordinate (l1) at (0,0);

    \coordinate (r6) at (7,4.5);    
    \coordinate (r7) at (7,4);
    \coordinate (r8) at (7,3.5);
    \coordinate (r9) at (7,3);
    \coordinate (r10) at (7,2.5);
    \coordinate (r5) at (7,2);        
    \coordinate (r4) at (7,1.5);    
    \coordinate (r3) at (7,1);
    \coordinate (r2) at (7,0.5);
    \coordinate (r1) at (7,0);

    \coordinate (v11) at (2.5,0);
    
    \coordinate (v21) at (2,0.5);
    \coordinate (v22) at (2.5,0.5);
    \coordinate (v23) at (4,0.5);

    \coordinate (v31) at (1.5,1);
    \coordinate (v32) at (2,1);
    \coordinate (v33) at (3.5,1);
    \coordinate (v34) at (4,1);
    \coordinate (v35) at (5.5,1);
    
    \coordinate (v41) at (1,1.5);
    \coordinate (v42) at (1.5,1.5);
    \coordinate (v43) at (3,1.5);
    \coordinate (v44) at (3.5,1.5);
    \coordinate (v45) at (5,1.5);
    \coordinate (v46) at (5.5,1.5);
    \coordinate (v47) at (6,1.5);
    
    \coordinate (v51) at (0.5,2);
    \coordinate (v52) at (3,2);
    \coordinate (v53) at (4.5,2);
    \coordinate (v54) at (6,2);

    \coordinate (v101) at (1,2.5);    
    \coordinate (v102) at (3,2.5);
    \coordinate (v103) at (5,2.5);
    \coordinate (v104) at (6,2.5);

    \coordinate (v91) at (0.5,3);
    \coordinate (v92) at (1.5,3);
    \coordinate (v93) at (3,3);
    \coordinate (v94) at (3.5,3);
    \coordinate (v95) at (4.5,3);
    \coordinate (v96) at (5.5,3);
    \coordinate (v97) at (6,3);
    
    \coordinate (v81) at (1.5,3.5);
    \coordinate (v82) at (2,3.5);
    \coordinate (v83) at (3.5,3.5);
    \coordinate (v84) at (4,3.5);
    \coordinate (v85) at (5.5,3.5);

    \coordinate (v71) at (2,4);
    \coordinate (v72) at (2.5,4);
    \coordinate (v73) at (4,4);

    \coordinate (v61) at (2.5,4.5);
    
    \draw[black] (-0.5,0) node  [xscale = 0.8, yscale = 0.8] {$1$};
    \draw[black] (-0.5,0.5) node  [xscale = 0.8, yscale = 0.8] {$2$};
    \draw[black] (-0.5,1) node  [xscale = 0.8, yscale = 0.8] {$3$};
    \draw[black] (-0.5,1.5) node  [xscale = 0.8, yscale = 0.8] {$4$};
    \draw[black] (-0.5,2) node  [xscale = 0.8, yscale = 0.8] {$5$};
    \draw[black] (-0.5,2.5) node  [xscale = 0.8, yscale = 0.8] {$10$};
    \draw[black] (-0.5,3) node  [xscale = 0.8, yscale = 0.8] {$9$};
    \draw[black] (-0.5,3.5) node  [xscale = 0.8, yscale = 0.8] {$8$};
    \draw[black] (-0.5,4) node  [xscale = 0.8, yscale = 0.8] {$7$};
    \draw[black] (-0.5,4.5) node  [xscale = 0.8, yscale = 0.8] {$6$};
    
    \draw[black] (7.5,0) node  [xscale = 0.8, yscale = 0.8] {$1$};
    \draw[black] (7.5,0.5) node  [xscale = 0.8, yscale = 0.8] {$2$};
    \draw[black] (7.5,1) node  [xscale = 0.8, yscale = 0.8] {$3$};
    \draw[black] (7.5,1.5) node  [xscale = 0.8, yscale = 0.8] {$4$};
    \draw[black] (7.5,2) node  [xscale = 0.8, yscale = 0.8] {$5$};
    \draw[black] (7.5,2.5) node  [xscale = 0.8, yscale = 0.8] {$10$};
    \draw[black] (7.5,3) node  [xscale = 0.8, yscale = 0.8] {$9$};
    \draw[black] (7.5,3.5) node  [xscale = 0.8, yscale = 0.8] {$8$};
    \draw[black] (7.5,4) node  [xscale = 0.8, yscale = 0.8] {$7$};
    \draw[black] (7.5,4.5) node  [xscale = 0.8, yscale = 0.8] {$6$};

    \draw[draw=dorange, line width=2pt] (l1)  -- (v11);
    \draw[draw=black, line width=1pt] (v11) -- (r1);

    \draw[draw=dviolet, line width=2pt] (l2) -- (v21);
    \draw[draw=black, line width=1pt] (v21) -- (v22);
    \draw[draw=dorange, line width=2pt] (v22) -- (v23);
    \draw[draw=black, line width=1pt] (v23) -- (r2);

    \draw[draw=dgreen, line width=2pt] (l3) -- (v31);
    \draw[draw=black, line width=1pt] (v31) -- (v32);
    \draw[draw=dviolet, line width=2pt] (v32) -- (v33);
    \draw[draw=black, line width=1pt] (v33) -- (v34);
    \draw[draw=dorange, line width=2pt] (v34) -- (r3);

    \draw[draw=blue, line width=2pt] (l4) -- (v41);
    \draw[draw=black, line width=1pt] (v41) -- (v42);
    \draw[draw=dgreen, line width=2pt] (v42) -- (v43);
    \draw[draw=black, line width=1pt] (v43) -- (v44);
    \draw[draw=dviolet, line width=2pt] (v44) -- (v45);
    \draw[draw=black, line width=1pt] (v45) -- (r4);
    
    \draw[draw=red, line width=2pt] (l5) -- (v51);
    \draw[draw=black, line width=1pt] (v51) -- (v52);
    \draw[draw=dgreen, line width=2pt] (v52) -- (v53);
    \draw[draw=black, line width=1pt] (v53) -- (r5);

    \draw[draw=black, line width=1pt] (l10) -- (v101);
    \draw[draw=blue, line width=2pt] (v101) -- (v102);
    \draw[draw=black, line width=1pt] (v102) -- (v103);
    \draw[draw=dviolet, line width=2pt] (v103) -- (v104);
    \draw[draw=black, line width=1pt] (v104) -- (r10);

    \draw[draw=black, line width=1pt] (l9) -- (v91);
    \draw[draw=red, line width=2pt] (v91) -- (v92);
    \draw[draw=black, line width=1pt] (v92) -- (v93);
    \draw[draw=blue, line width=2pt] (v93) -- (v94);
    \draw[draw=black, line width=1pt] (v94) -- (v95);
    \draw[draw=dgreen, line width=2pt] (v95) -- (v96);
    \draw[draw=black, line width=1pt] (v96) -- (v97);
    \draw[draw=dviolet, line width=2pt] (v97) -- (r9);

    \draw[draw=black, line width=1pt] (l8) -- (v81);
    \draw[draw=red, line width=2pt] (v81) -- (v82);
    \draw[draw=black, line width=1pt] (v82) -- (v83);
    \draw[draw=blue, line width=2pt] (v83) -- (v84);
    \draw[draw=black, line width=1pt] (v84) -- (v85);
    \draw[draw=dgreen, line width=2pt] (v85) -- (r8);

    \draw[draw=black, line width=1pt] (l7) -- (v71);
    \draw[draw=red, line width=2pt] (v71) -- (v72);
    \draw[draw=black, line width=1pt] (v72) -- (v73);
    \draw[draw=blue, line width=2pt] (v73) -- (r7);

    \draw[draw=black, line width=1pt] (l6) -- (v61);
    \draw[draw=red, line width=2pt] (v61) -- (r6);

    \draw[draw=blue, line width=2pt, ->]  (v41) .. controls (1.25,1.75) and (1.25,2.25) .. (v101) node[midway, below left] {\scriptsize \textcolor{blue}{$t_{1}$}};
            
    \draw[draw=red, line width=2pt, ->]  (v51) .. controls (0.75,2.25) and (0.75,2.75) .. (v91) node[midway, below left] {\scriptsize \textcolor{red}{$-t_{1}$}};
    
    \draw[draw=dgreen, line width=2pt, ->]  (v31) -- (v42) node[below right=0.5mm] {\scriptsize \textcolor{dgreen}{$t_{2}$}};
    \draw[draw=red, line width=2pt, ->]  (v92) -- (v81) node[below right] {\scriptsize \textcolor{red}{$-t_{2}$}};
    
    \draw[draw=dviolet, line width=2pt, ->]  (v21) -- (v32) node[below right=0.5mm] {\scriptsize \textcolor{dviolet}{$t_{3}$}};
    \draw[draw=red, line width=2pt, ->]  (v82) -- (v71) node[below right] {\scriptsize \textcolor{red}{$-t_{3}$}};
    
    \draw[draw=dorange, line width=2pt, ->]  (v11) -- (v22) node[below right=0.5mm] {\scriptsize \textcolor{dorange}{$t_{4}$}};
    \draw[draw=red, line width=2pt, ->]  (v72) -- (v61) node[below right] {\scriptsize \textcolor{red}{$-t_{4}$}};

    \draw[draw=dgreen, line width=2pt, ->]  (v43) -- (v52) node[below right=0.5mm] {\scriptsize \textcolor{dgreen}{$t_{5}$}};
    \draw[draw=blue, line width=2pt, ->]  (v102) -- (v93) node[below right] {\scriptsize \textcolor{blue}{$-t_{5}$}};

    \draw[draw=dviolet, line width=2pt, ->]  (v33) -- (v44) node[below right=0.5mm] {\scriptsize \textcolor{dviolet}{$t_6$}};
    \draw[draw=blue, line width=2pt, ->]  (v94) -- (v83) node[below right] {\scriptsize \textcolor{blue}{$-t_{6}$}};

    \draw[draw=dorange, line width=2pt, ->]  (v23) -- (v34) node[below right=0.5mm] {\scriptsize \textcolor{dorange}{$t_{7}$}};
    \draw[draw=blue, line width=2pt, ->]  (v84) -- (v73) node[below right] {\scriptsize \textcolor{blue}{$-t_{7}$}};
    
    \draw[draw=dviolet, line width=2pt, ->]  (v45) .. controls (5.25,1.75) and (5.25,2.25) .. (v103)
            node[midway, below left] {\scriptsize \textcolor{dviolet}{$t_{8}$}};
            
    \draw[draw=dgreen, line width=2pt, ->]  (v53) .. controls (4.75,2.25) and (4.75,2.75) .. (v95)
            node[midway, below left] {\scriptsize \textcolor{dgreen}{$-t_{8}$}}; 

    \draw[line width=2pt, ->]  (v35) -- (v46) node[below right=0.5mm] {\scriptsize $t_{9}$};
    \draw[draw=dgreen, line width=2pt, ->]  (v96) -- (v85) node[below right] {\scriptsize \textcolor{dgreen}{$-t_{9}$}};

    \draw[line width=2pt, ->]  (v47) -- (v54)node[below right=0.5mm] {\scriptsize  $t_{10}$};
    \draw[draw=dviolet, line width=2pt, ->]  (v104) -- (v97) node[below right] {\scriptsize \textcolor{dviolet}{$-t_{10}$}};
\end{tikzpicture}
\caption{Illustration of the path collection corresponding to $M_{1,2}$.}
\end{figure}
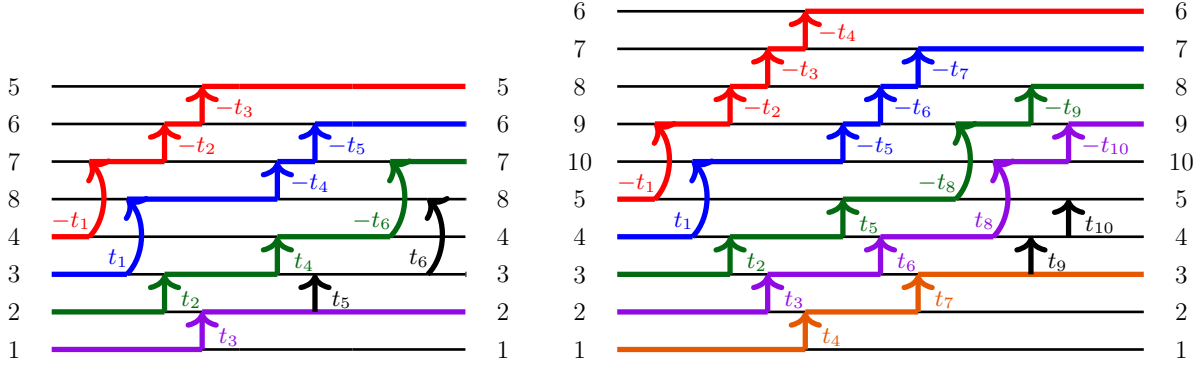
\vspace{-1.5em}
\end{center}
    The path collection in example \ref{example:n5pathcollection} is obtained from the previous one by changing the path starting at $2$ from left greedy to right greedy, and it corresponds to $M_{2,2}=(t_1\cdots t_7)^2t_8t_9$ for $n=5$.
\end{example}

\begin{proof}[Proof of \Cref{prop:UniquePathColl}]
    For each $1\leq j \leq k \leq n-1$ we need to argue that the path collection consisting of a path from $1 \to n-k$, non-intersecting right greedy paths from sources $2,\ldots,j$ and left greedy paths from sources $j+1,\ldots,n$ is unique and has the claimed sink set. By \Cref{lem:LGVdiagramLeftGreedy} and \eqref{eq:curlyW0}, the left greedy path collection originating from $[n]$ is explicitly given by:
    \[
    \begin{array}{ll}
        1 \to n &\text{ if } n  \text{ is odd}\\
        1 \to 2n &\text{ if } n \text{ is even}
    \end{array}, \quad 2 \to 2n-1, \quad 3 \to 2n-2, \quad \cdots, \quad n \to n+1.
    \]
    It also follows from \Cref{lem:LGVdiagramLeftGreedy} that the left greedy paths from $j+1,\ldots,n$ are unique. Moreover, combining \Cref{lem:LGVdiagramLeftGreedy} with the explicit description of the paths in the left greedy path collection above, we see that there is no path collection originating from any source set contained in $[n]$ other than $j+1,\ldots,n$ with the same sink set.  Then, we just need to check the result for the subcollection corresponding to the paths with source vertices $1,\ldots,j$. 
    \smallskip
    We start by showing there exists a unique path in the LGV diagram from $1$ to $n-k$. First, we describe the left greedy path ${\rm LGP}(1)$ starting at $1$. From $\mathcal{W}_0$ we see that the only arrow with a tail in strand $1$ is $a_1^{(1)}$, so it must be the first arrow in ${\rm LGP}(1)$. The unique arrow to the right of $a_1^{(1)}$ which starts on strand $2$ is $a_2^{(2)}$. Similarly, for $i\in[n-3]$, $a_{i+1}^{(i+1)}$ is the unique arrow originating on strand $i+1$ to the right of the head of $a_i^{(i)}$. Thus, ${\rm LGP}(1)$ begins with $a_1^{(1)}a_2^{(2)}a_3^{(3)}\cdots a_{n-2}^{(n-2)}$. There is then a unique arrow originating on strand $n-1$ to the right of the head of $a_{n-2}^{(n-2)}$. For $n$ odd, it is $a_{n-1}^{(n-1)}$ and for $n$ even, $b_{n-1}^{(n/2)}$. Accordingly, for $1\leq k \leq n-1$ the unique path from $1$ to $n-k$ follows ${\rm LGP}(1)$ until it reaches strand $n-k$, and then terminates horizontally. Explicitly, it is given by the arrows $a_1^{(1)}a_2^{(2)}\cdots a_{n-k-1}^{(n-k-1)}$.
    \smallskip
    When we change ${\rm LGV}(2)$ to the right greedy path originating at $2$, it only uses the arrows forced by the fact that the bottom $n-k$ strands are blocked at some point by the path $1 \to n-k$. That is, the path $2 \to n-k+1$ is given by $a_2^{(1)}a_3^{(2)}\cdots a_{n-k-2}^{(n-k-1)}$. Similarly, the path starting at $3$ will only use the arrows forced by the fact that the path $2\to n-k+1$ at some point blocks the bottom $n-k+1$ strands, and so on. Continuing this argument, the path $j \to n-k+j-1$ is given by $a_{j}^{(1)}a_{j+1}^{(2)}\cdots a_{n-k+j-2}^{(n-k-1)}$. In particular, these paths have the desired sinks and each arrow they use is forced. Moreover, note that since these paths all terminate at sinks in $[n]$, they do not use any arrows that skip over other arrows. Thus, the path originating at $i$ must terminate below the path originating at $i+1$ for $1\leq i < j$. Thus, the path collection described is the unique non-intersecting path collection from $\{1,\ldots, j\}$ to $\{n-k, \ldots, n-k+j-1\}$. 
\end{proof}

\begin{lemma}\label{lem:determiningtheparameters}
    Let $X=X(t_1,\ldots, t_N)$ be as in \Cref{subsec:LieTheory} for some parameters $t_i \in \RR^*$. For $1 \leq j \leq k \leq n-1$, let us denote by $\#(j,k)$ the index
    \[
    \#(j,k) := N - \left(\binom{k}{2} + (j-1)\right).
    \]
    Then, the signs of the minors $\big\{M_{j',k'} \colon (j',k')\leq (j,k) \textnormal{ in reverse lexicographic order}\big\}$ determine the sign of $t_{\#(j,k)}$. Moreover, the indeterminates $t_i$ can be written as Laurent monomials in the minors $M_{j,k}$.
\end{lemma}
\begin{proof}
    We will prove the result by reverse induction on $\#(j,k)$. Equivalently, we induct on $(j,k)$ in reverse lexicographic order. By \Cref{prop:UniquePathColl}, to determine $M_{j,k}$ we need to look at the weights appearing in the path collection given by 
    \[
    1 \to n-k, \quad \cdots, \quad j \to n-k+j-1, \quad j+1 \to 2n-j, \quad \cdots ,\quad n \to n+1.
    \]
    We will denote this path collection as $(n-k,\ldots,n-k+j-1,2n-j,\ldots, n+1)$. We know from \Cref{prop:UniquePathColl} that each $M_{j,k}$ is a monomial in the $t_i$ with coefficient $\pm1$. We write $M_{j,k}=\epsilon_{j,k}\mathcal{M}_{j,k}$, where $\epsilon_{j,k}\in \{1,-1\}$ and $\mathcal{M}_{j,k}$ is a monomial in the $t_i$ with coefficient $1$. We will show in \Cref{lem:Mjkpositive} that $\epsilon_{j,k}=1$, but for the purposes of this lemma, it suffices to work with $\mathcal{M}_{j,k}$ instead of $M_{j,k}$. We denote the minor corresponding to the left greedy path collection via \Cref{LGVlemma} as $M=\epsilon \mathcal{M}$, where $\mathcal{M}$ is a monomial in the $t_i$ with coefficient $1$. By \Cref{lem:LGVdiagramLeftGreedy}, the left greedy path collection originating from $[n]$ unique and uses all arrows in the LGV diagram, so $\mathcal{M}=(t_1\cdots t_N)^2$.

    For $k=1, j=1$, the unique path collection $(n-1,2n-1,\ldots, n+1)$ is almost identical to the left greedy path collection originating from $[n]$, only the path originating at $1$ is modified. Recall from the proof of \Cref{prop:UniquePathColl} that for $k=1$ the unique path $1 \to n-1$ uses all the arrows other than the last one in ${\rm LGP}(1)$. Thus, in the path collection corresponding to $M_{1,1}$, the only arrow which is not used is $a_{n-1}^{(n-1)}$ if $n$ is odd and $b_{n-1}^{(n/2)}$ if $n$ is even. Since $\lambda(a_{n-1}^{(n-1)})=N$ for $n$ odd and $\lambda(b_{n-1}^{(n/2)})=N$ for $n$ even, $\mathcal{M}_{1,1}=(\prod_{i=1}^{N-1}t_i^2)t_N$. Also, since $\mathcal{M}$ is a perfect square, $\mathcal{M}_{1,1}$ determines the sign of $t_N=\frac{\mathcal{M}}{\mathcal{M}_{1,1}}$. We use the following lemma to complete the proof for $k=1, j=1$.

    \begin{lemma}
        The minor $M$ can be written as a Laurent monomial in $\{M_{1,1},\ldots,M_{n-1,n-1}\}$.
    \end{lemma}
    \begin{proof}
        We continue to work with $\mathcal{M}$ and $\mathcal{M}_{j,k}$, dealing with the sign analysis later. Note that $\mathcal{M} = \Pf(A)^2$ if $n$ is even and $\mathcal{M}=\Pf_{[n-1]}(A)^2$ if $n$ is odd. Thus, it suffices to determine $|\Pf(A)|$ or $|\Pf_{[n-1]}(A)|$, for $n$ even or odd, respectively, as a Laurent monomial in $\{\mathcal{M}_{1,1},\ldots, \mathcal{M}_{n-1,n-1}\}$. Using \cite[Theorem 1]{Weyman}, we can write the minors of a skew-symmetric matrix in terms of its Pfaffians. Applying this to the (unsigned) minors $\mathcal{M}_{j,k}$ of $A$ with $j=k$ yields
        \begin{align*}
            &\mathcal{M}_{k,k} =  |\Pf_{[n-k-1]}(A) \ \Pf_{1,\ldots,n-k,n}(A)| &&\text{if } n-k \text{ is odd}\\
            &\mathcal{M}_{k,k} = | \Pf_{[n-k]}(A) \ \Pf_{1,\ldots,n-k-1,n}(A)| &&\text{if } n-k \text{ is even.}
        \end{align*}
        Note that there is a common term appearing in both $\mathcal{M}_{k,k}$ and $\mathcal{M}_{k-1,k-1}$. Explicitly,
          \[
            \frac{\mathcal{M}_{k-1,k-1}}{\mathcal{M}_{k,k}} =\begin{cases}
                \left|\frac{\Pf_{[n-k+1]}(A)}{\Pf_{[n-k-1]}(A)}\right|,  &\text{if } n-k \text{ is odd}\\[8pt]
                \left|\frac{\Pf_{1,\ldots,n-k+1,n}(A)}{\Pf_{1,\ldots,n-k-1,n}(A)}\right|,&\text{if } n-k \text{ is even.}
            \end{cases}
          \]
          This implies that 
          \[
            \begin{cases}
                |\Pf_{[n-k+1]}(A)|=\frac{|\Pf_{[n-k-1]}(A)| \ \mathcal{M}_{k-1,k-1}}{\mathcal{M}_{k,k}},&\text{if } n-k \text{ is odd}\\[8pt]
                 |\Pf_{1,\ldots, n-k+1,n}(A)|=\frac{|\Pf_{1,\ldots, n-k-1,n}(A)| \ \mathcal{M}_{k-1,k-1}}{\mathcal{M}_{k,k}}, &\text{if } n-k \text{ is even.}
            \end{cases}
          \]
          
        Observe that $\mathcal{M}_{n-1,n-1}=|\Pf_{1,n}(A)|$. By reverse induction on $k$, we can write $|\Pf_{1,\cdots n-k,n}|$ if $n-k$ is odd and $|\Pf_{[n-k]}|$ if $n-k$ is even as a Laurent monomial in $\{\mathcal{M}_{k,k},\ldots, \mathcal{M}_{n-1,n-1}\}$. In particular, for $k=1$ we obtain expressions for $\Pf_{[n-1]}(A)$ for $n$ odd and $\Pf(A)$ for $n$ even as Laurent monomials in the $\mathcal{M}_{k,k}$, as desired.
    \end{proof}

    Assume we have determined the signs of $t_N,\ldots, t_{\#(j,k)+1}$ and we have written them as Laurent monomials in $\mathcal{M}_{j',k'}$ for $(j',k')\leq (j,k)$.
    
    \smallskip
    
    For $j=1$, the minor $M_{1,k}$ corresponds to the path collection using the unique path from $1 \to n-k$ and the left greedy path collection originating from $2,\ldots, n$. This differs from the path collection corresponding to $M_{1,1}$ only in the path starting at $1$. For an illustration of this fact, we refer the reader to \Cref{ex:somePathsCollectionsInLGVDiagram}. Hence, $\frac{\mathcal{M}_{1,1}}{\mathcal{M}_{1,k}}= t_{i_1}\cdots t_{i_{k-1}}$ where $i_1=\lambda(a_{n-k}^{(n-k)})<i_2=\lambda(a_{n-k+1}^{(n-k+1)}) \ldots <  i_{k-1}=\lambda(a_{n-2}^{(n-2)})<i_{k}=N $ are the indices of the labels of the last $k$ arrows in ${\rm LGP}(1)$. In other words, since $\mathcal{M}_{1,1}=(\prod_{i=1}^{N-1}t_i^2)t_N$, then all the variables $t_i, 1\leq i \leq N$ other than $t_{i_1},\cdots t_{i_{k-1}}, t_{i_k}=t_N$ appear squared in $\mathcal{M}_{1,k}$. We now prove that $i_1=\#(1,k)$, which allows us to express $t_{\#(1,k)}$ as a Laurent monomial in $\mathcal{M}_{1,1},\mathcal{M}_{1,k}$ and $\{t_i\colon i> \#(1,k)\}$. By induction, this proves for $j=1$ that the variables can written as Laurent monomials in the unsigned minors $\mathcal{M}_{j',k'}$.
    
    Recall the sequence $\mathcal{W}_0$ of arrows appearing in the LGV graph and that, up to sign, $t_i$ is the weight of each of a pair of consecutive arrows in $\mathcal{W}_0$. To compute $i_1$, observe that $a_{n-k}^{(n-k)}$ is one of the last two arrows in some parentheses pair. 
    
    For $n-k$ even, $a_{n-k}^{(n-k)}$ lies in parentheses pair $\frac{n-k}{2}$. Each parentheses pair is preceded by a pair of $b$ arrows. As each set of parentheses in parentheses pair $l$ contain $n-2l$ pairs of arrows, we have
  \begin{align}
  \begin{split}\label{eq:lambdaofalpha}
      \lambda\left(a_{n-k}^{(n-k)}\right) 
      &= \sum_{l=1}^{\frac{n-k}{2}} \left[ 1+2(n-2l)\right] \\
      &= \frac{n-k}{2}(1+2n)-2\frac{n-k}{2}\left(\frac{n-k}{2}+1\right)\\
      &=\left(\frac{n-k}{2}\right)(n+k-1)\\
      &=\frac{1}{2}(n(n-1)-k(k-1))=N-\binom{k}{2}.
      \end{split}
  \end{align}

  A similar calculation, taking into account that there is one extra set of parentheses of length $k-1$ and an extra pair of $b$ arrows, gives the same result for $n-k$ odd. Thus, $i_1=N-\binom{k}{2}=\#(k,1)$.
  
  \smallskip
  
  For $1<j\leq k$, the minor $M_{j,k}$ corresponds to the path collection $(n-k,\ldots, n-k+j-1,2n-j,\ldots,n+1)$. This is given by taking the path collection corresponding to $M_{j-1,k}$ and modifying the path starting at $j$. Concretely, we are removing arrows from ${\rm LGP}(j)$ in order to make it right greedy. Hence, $\frac{\mathcal{M}_{j-1,k}}{\mathcal{M}_{j,k}}= t_{i_1}\cdots t_{i_l}$ where $i_1\leq \ldots \leq  i_l$ are the indices of the labels of the arrows in ${\rm LGP}(j)$ but not in the path $j \to n-k+j-1$. 
  \smallskip
  Moreover, we can explicitly identify the first arrow, $\alpha_j$, along the path ${\rm LGP}(j)$ that is not on the path $j \to n-k+j-1$. Observe that by definition of a left greedy path, $\alpha_j$ is the first arrow starting on strand $n-k+j-1$ to the right of the head of $a_{n-k+j-2}^{(n-k-1)}$. As we can identify where $a_{n-k+j-2}^{(n-k-1)}$ lies in $\mathcal{W}_0$, we can identify $\alpha_j$:.
  \begin{itemize}
      \item If $j<k$, then $n-k+j-1<n-1$ and the only arrows originating on strand $n-k+j-1$ are of the form $a_{n-k+j-1}^{(l)}$ for some $l$. If we also ask for this arrow to be to the right of the head of $a_{n-k+j-2}^{(n-k-1)}$, we must have $l\geq n-k-1$. Since the parentheses in \eqref{eq:curlyW0} have decreasing lower indices, we obtain $\alpha_j=a_{n-k+j-1}^{(n-k)}$.
      \item If $j=k$, then $n-k+j-1=n-1$ and additionally to the $a_{n-1}^{(l)}$ arrows originating on strand $n-1$, we also have $b_{n-1}^{(l')}$ arrows. We must consider two cases.
      \begin{itemize}
          \item  For $n-k-1$ odd, the arrow $a_{n-2}^{(n-k-1)}$ is in the first set of parentheses in parentheses pair $\frac{n-k}{2}$ of $\mathcal{W}_0$. Then, the next available arrow originating on strand $n-1$ is $\alpha_j=a_{n-1}^{(n-k)}$, which is in the second set of parentheses in the same parentheses pair. 
          \item  For $n-k-1$ even, $a_{n-2}^{(n-k-1)}$ is in the second parentheses of parentheses pair $m=\frac{n-k-1}{2}$ in $\mathcal{W}_0$. Then, there is a jump arrow on strand $n-1$ before the next parentheses pair, implying $\alpha_j = b_{n-1}^{(m+1)}$.
      \end{itemize}
  \end{itemize}
  Now, we want to compute $\lambda(\alpha_j)$. When $\alpha_j=a_{n-k+j-1}^{(n-k)}$ it follows from \eqref{eq:curlyW0} that $\lambda(\alpha_j)=\lambda(a_{n-k}^{(n-k)})-(j-1)$. When $\alpha_j=b_{n-1}^{(m+1)}$, we claim the same is true. For $j=k=1$, we must have $n$ even and this follows from our analysis of the $j=k=1$ case earlier in this proof. For $j=k>1$, note that $\alpha_{j-1}=a_{n-2}^{(n-k)}$ and we have $\lambda(\alpha_{j-1})=\lambda(a_{n-k}^{(n-k)})-(j-2)$. As $\alpha_{j-1}=a_{n-2}^{(n-k)}$ is the first arrow in parentheses pair $m+1=\frac{n-k+1}{2}$ of $\mathcal{W}_0$, it follows that $\lambda(\alpha_j)=\lambda(\alpha_{j-1})-1=\lambda(a_{n-k}^{(n-k)}-(j-1)$, as claimed. In conclusion, using \eqref{eq:lambdaofalpha}, $i_1=\lambda(\alpha_j)=N-\binom{k}{2}-(j-1)=\#(j,k)$, concluding the proof.
\end{proof}

\subsection{Sign analysis}\label{subsec:signAnalysis}

We now show that the parameters $(t_i, \ 1 \leq i \leq N)$ are positive if and only if the minors $(M_{j,k}, \ 1 \leq j \leq k \leq n-1)$ are positive. We begin with the following lemma.

\begin{lemma}\label{lem:signXtoA}
    For $1\leq j \leq k \leq n-1$, we have 
    \[
    \Delta^{n-k,\ldots,n-k+j-1,n+1,n+2,\ldots,2n-j}(X)=(-1)^{j(n-1-k)}\Delta_{\{1,\ldots,n-k-1,n-k+j, \ldots, n\}}^{\{1,2,\ldots, n-j\}}(A).
    \]
\end{lemma}
\begin{proof}
    This follows immediately from Laplace expansion.
\end{proof}

\begin{lemma}\label{lem:Mjkpositive}
    If $t_1,\ldots, t_N >0$ then the minors $M_{j,k}$ are positive for all $1\leq j \leq k \leq n-1$.
\end{lemma}
\begin{proof}

    Concretely we prove that $M_{j,k}$ are monomials on the $t_i$ with coefficient $1$. That is, using notation from \Cref{lem:determiningtheparameters}, for all $1\leq j\leq k\leq n-1$, $\epsilon_{j,k}=1$ where $M_{j,k}=\epsilon_{j,k}\mathcal{M}_{j,k}$.

    \smallskip
    
    Note that we can divide the LGV diagram into two parts: from strand $1$ to $n$ (\textit{lower part}) and from strand $2n$ to $n+1$ (\textit{upper part}). These are connected by the ``jump arrows'' $b_i^{(j)}$. The lower part of the LGV diagram has all weights positive and the upper part has all weights negative. Therefore, to determine the sign of the minors $M_{j,k}$, we need to count the number of arrows in the upper part of the diagram and the number of jump arrows of negative weight used by the corresponding path collection; we also must take into account the sign from the LGV lemma, \Cref{LGVlemma}, the sign from \Cref{lem:signXtoA}, and the sign appearing in the definition of $M_{j,k}$.

    \smallskip
    
    The path collection corresponding to $M_{j,k}$ is $(n-k,\ldots, n-k+j-1,2n-j,\ldots,n+1)$. Since the paths starting at $1,\ldots, j$ remain in the lower part of the diagram, their weights are positive. The paths starting at $j+1,\ldots,n$ are left greedy paths which use some arrows in the lower part of the diagram then take one jump arrow and continue in the upper part. Explicitly, they use all the arrows in the upper part of the graph coming from the first $n-j$ sets of parentheses in \eqref{eq:curlyW0}. Each path uses all the arrows in the upper part of the graph coming from a single set of parentheses. For an illustration of this fact, we refer the reader to \Cref{ex:somePathsCollectionsInLGVDiagram}. If $n-j$ (the number of left greedy paths) is even then we can pair paths with consecutive sources: 
    \[
        \big\{\LGP(n), \LGP(n-1)\big\}, \big\{\LGP(n-2),\LGP(n-3)\big\},\ldots, \big\{\LGP(n-j+1),\LGP(n-j)\big\}.
    \]
    Each pair of paths then uses all the arrows in the upper part of the graph coming from a parentheses pair in \eqref{eq:curlyW0}. Each set of parentheses in a parentheses pair has the same number of arrows coming from the upper part of the graph. Thus, these arrows contribute an even number of factors of $-1$. Then, the sign given by the arrows depends on the parity of the number of jump arrows of negative weight used. Each parentheses pair is preceded by a pair of jump arrows, one each of positive and negative weight. So, there are exactly $\frac{n-j}{2}$ jump arrows of negative weight used by the path collection contributing a sign of $(-1)^{\frac{n-j}{2}}$. The permutation that sorts $(n-k,\ldots, n-k+j-1,2n-j,\ldots,n+1)$ can be written as the product of $\frac{n-j}{2}$ transpositions corresponding to exchanging $n+a$ and $2n-j-a+1$ for $1\leq a\leq \frac{n-j}{2}$. Thus, the sign contributed by \Cref{LGVlemma} is $(-1)^{\frac{n-j}{2}}$. If $n-j$ is odd, a similar argument shows that the sign contributed by arrows of negative weight is $(-1)^{\frac{n-j-1}{2}+n+1}$ and the sign contributed by \Cref{LGVlemma} is $(-1)^{\frac{n-j-1}{2}}$. In either case, the sign contributed by arrows of negative weight and by \Cref{LGVlemma} together is $(-1)^{j(n-1)}$. Multiplying this by the sign $(-1)^{j(n-1-k)}$ contributed by \Cref{lem:signXtoA}, we obtain a sign of $(-1)^{jk}$. This matches the sign that appears in \Cref{def:SpecialMinorsPfaff}, the definition of $M_{j,k}$.
\end{proof}

\begin{lemma}\label{lem:Mjkdeterminematrix}
 Let $A$ be a $n \times n$ skew-symmetric matrix such that $M_{j,k}(A) \neq 0$ for all $1\leq j \leq k\leq n-1$. Then, these minors fully determine the matrix $A$.
\end{lemma}

\begin{proof}
    Let $U$ be the set of matrices $A$ in $\SS_n$ such that $M_{j,k}(A) \neq 0$ for any $1 \leq j \leq k \leq n-1$. Similar to \Cref{def:LusztigPositive}, we denote by $V$ the following set
    \[
        V := \Big\{ A(t_1, \dots, t_N) \colon t_1,\dots,t_N \in \RR^\ast \Big\},
    \]
    where we recall that $A(t_1,\dots,t_N)$ is the $n\times n$ skew-symmetric matrix with polynomial entries in the parameters $t_i$ obtained after row-reducing in \eqref{eq:posorthogonalgrassmannian}. 
    We note that the sets $U$ and $V$ are both Zariski dense in $\OGr(n,2n)$. By \Cref{lem:Mjkpositive}, we have $V \subset U$. Let $\varphi \colon U  \to V$ be a rational map which takes a skew-symmetric matrix $A$ in $U$, determines the $t_i$'s in the Marsh-Rietsch parametrization, as described in the proof of \Cref{lem:determiningtheparameters}, and then, using these $t_i$'s, outputs an element in $V$ via the map \eqref{eq:posorthogonalgrassmannian}. Note that the function $\varphi$ extends to a rational map $\tilde{\varphi} \colon \OGr(n,2n) \dashrightarrow \OGr(n,2n)$. Restricting $\tilde{\varphi}$ to $V$ gives the identity and, since $V$ is Zariski dense in $\OGr(n,2n)$, we deduce that $\tilde{\varphi} = \rm{id}$. In particular, $\varphi$ as a function applied to $A\in U$ depends only on the minors $M_{j,k}(A)$ and outputs $A$. Thus, the minors $M_{j,k}(A)$ fully determine the matrix $A$.
    \end{proof}

    \begin{proof}[\bf Proof of \Cref{thm:Main}]
        By \Cref{lem:determiningtheparameters} and \Cref{lem:Mjkpositive}, for a point $A(t_1, \dots, t_N)$ with $t_i \in \RR^\ast$, the parameters $t_i$ are positive if and only if all the $M_{j,k}(A)$ are. The first statement of the result then follows from \Cref{lem:Mjkdeterminematrix}. For the second statement, observe that any positivity test of $\OGr^{>0}(n,2n)$ using regular functions $f_1, \dots, f_r$ can be turned into a positivity test for the positive orthant $(\RR_{>0})^{N}$ using the Marsh-Rietsh parametrization. The minimal number of inequalities that cut out the orthant $(\RR_{>0})^N$ in $\RR^N$ is $N$, so we deduce that $r \geq N$. The positivity test in \Cref{thm:Main} uses exactly $N = \binom{n}{2}$ inequalities so it is indeed a minimal positivity test.
    \end{proof}

\section{Nonnegativity test and boundary cells}\label{sec:4}

\subsection{Nonnegativity test}

Recall that the totally nonnegative orthogonal Grassmannian $\OGr^{\geq 0}(n,2n)$ is the Euclidean closure of $\OGr^{>0}(n,2n)$. We now prove \Cref{thm:Nonnegative}, which provides a nonnegativity criterion for skew-symmetric matrices. To do so, we use \Cref{lem:semiGroup}, which exploits the semigroup property of $\SO^{>0}(2n)$.

\smallskip

Recall that $\SO(2n)$ acts on $\OGr(n,2n)$ via matrix multiplication, i.e. given $M \in \SO(2n)$ and $X=\rowspan (N) \in \OGr(n,2n)$ we have $X \cdot M = \rowspan (NM)$.

\begin{lemma}[]\label{lem:semiGroup}
    For any $Z \in \SO^{ > 0}(2n)$ and $X\in \OGr^{\geq 0}(n,2n)$, we have $X \cdot Z \in \OGr^{>0}(n,2n)$.
\end{lemma}
\begin{proof}
    In Proposition 3.2(c) and Proposition 8.17 in \cite{Lusztig1} this statement is proven for the complete flag. By Theorem 3.4 in \cite{Lusztig2} the result follows for partial flags by projecting.
\end{proof}

\begin{proof}[\bf Proof of \Cref{thm:Nonnegative}]
    The equivalence between (\ref{nonnegativeitem2}) and (\ref{nonnegativeitem3}) follows from \Cref{thm:Main}.
    
    To show the equivalence of (\ref{nonnegativeitem1}) and (\ref{nonnegativeitem2}), let $Z(\epsilon) \in \SO^{>0}(2n)$ be a smooth $1$-parameter family that converges to $\Id_{2n}$ as $\epsilon \to 0$. Fix $X \in \OGr^{\geq 0}(n,2n)$. Then, by virtue of \Cref{lem:semiGroup}, $X \cdot Z(\epsilon)$ is totally positive. By our positivity criterion, \Cref{thm:Main}, we deduce that the minors $M_{j,k}(B(\epsilon))$ are positive for any $\epsilon > 0$.  Conversely, suppose that the leading term in the Taylor expansion of all $M_{j,k}(B(\epsilon))$ is positive for small enough $\epsilon$. Then using \Cref{thm:Main} we deduce that $X(\epsilon) \in \OGr(n,2n)^{> 0}$ for $\epsilon$ small enough. Since $Z(\epsilon) \xrightarrow[]{\epsilon \to 0} \Id_{2n}$, we conclude that $X(\epsilon)  \xrightarrow[]{\epsilon \to 0} X$. This means that $X \in \overline{\OGr^{>0}(n,2n)} = \OGr^{\geq 0}(n,2n)$. Thus, (\ref{nonnegativeitem1}) and (\ref{nonnegativeitem2}) are equivalent.
    \smallskip

    We now show that we can choose $Z(\epsilon)$ to have polynomial entries in $\epsilon$. The parametrization of $\OGr^{>0}(n,2n)$ in equation \eqref{eq:parameterizationCompleteFlag} is a projection of the parametrization of the complete flag variety given in \cite{MR}. Thus, for $\epsilon >0$ the following matrix representing a point in the positive complete flag variety for $\SO(2n)$:
    \[
        z(\epsilon) = x_{j_1}(\epsilon)\cdots x_{j_{\hat{N}}}(\epsilon).
    \]
    We then have $ Z(\epsilon) := z(\epsilon)^Tz(\epsilon) \in \SO(2n)^{>0}$. To see why, recall that the Marsh-Rietsch parametrization of the totally positive complete flag variety is obtained by representing a complete flag with a matrix in the totally positive unipotent radical $U_{>0}^{+}$ of $\SO(2n)$. Using \cite[Section 11]{MR}, we observe that $z(\epsilon) \in U_{>0}^{+}$ and $z(\epsilon)^T \in U_{>0}^{-}$. It follows from the definition of $\SO(2n)^{>0}$, see \cite[2.12]{Lusztig1}, that $Z(\epsilon)$ is totally positive. Finally, note that as $\epsilon \to 0$ each $x_i(\epsilon)$ converges to the identity, so $Z(\epsilon) \xrightarrow[\epsilon \to 0]{} \Id_{2n}$. Moreover, since the entries of each $x_{j_i}(\epsilon)$ are polynomial in $\epsilon$, the same is true of the entries and minors of $Z(\epsilon)$. 
\end{proof}

\begin{example}[$n=4$]\label{ex:NonNegn4}
   The matrix $Z(\epsilon)$ from the parametrization in \cite[Section 11]{MR} is
    \[
    \resizebox{1\textwidth}{!}{$\begin{bmatrix}
       1&4\,\epsilon&10\,\epsilon^{2}&7\,\epsilon^{3}&-\epsilon^{6}&5\,\epsilon^{5}&-11\,\epsilon^{4}&13\,\epsilon^{3}\\
       4\,\epsilon&16\,\epsilon^{2}+1&40\,\epsilon^{3}+4\,\epsilon&28\,\epsilon^{4}+4\,\epsilon^{2}&-4\,\epsilon^{7}-\epsilon^{5}&20\,\epsilon^{6}+4\,\epsilon^{4}&-44\,\epsilon^{5}-7\,\epsilon^{3}&52\,\epsilon^{4}+6\,\epsilon^{2}\\
       10\,\epsilon^{2}&40\,\epsilon^{3}+4\,\epsilon&100\,\epsilon^{4}+16\,\epsilon^{2}+1&70\,\epsilon^{5}+16\,\epsilon^{3}+2\,\epsilon&-10\,\epsilon^{8}-4\,\epsilon^{6}-\epsilon^{4}&50\,\epsilon^{7}+16\,\epsilon^{5}+3\,\epsilon^{3}&-110\,\epsilon^{6}-28\,\epsilon^{4}-4\,\epsilon^{2}&130\,\epsilon^{5}+24\,\epsilon^{3}+2\,\epsilon\\
       7\,\epsilon^{3}&28\,\epsilon^{4}+4\,\epsilon^{2}&70\,\epsilon^{5}+16\,\epsilon^{3}+2\,\epsilon&49\,\epsilon^{6}+16\,\epsilon^{4}+4\,\epsilon^{2}+1&-7\,\epsilon^{9}-4\,\epsilon^{7}-2\,\epsilon^{5}-\epsilon^{3}&35\,\epsilon^{8}+16\,\epsilon^{6}+6\,\epsilon^{4}+2\,\epsilon^{2}&-77\,\epsilon^{7}-28\,\epsilon^{5}-8\,\epsilon^{3}-2\,\epsilon&91\,\epsilon^{6}+24\,\epsilon^{4}+4\,\epsilon^{2}\\
       -\epsilon^{6}&-4\,\epsilon^{7}-\epsilon^{5}&-10\,\epsilon^{8}-4\,\epsilon^{6}-\epsilon^{4}&-7\,\epsilon^{9}-4\,\epsilon^{7}-2\,\epsilon^{5}-\epsilon^{3}&\epsilon^{12}+\epsilon^{10}+\epsilon^{8}+10\,\epsilon^{6}+36\,\epsilon^{4}+16\,\epsilon^{2}+1&-5\,\epsilon^{11}-4\,\epsilon^{9}-3\,\epsilon^{7}-14\,\epsilon^{5}-24\,\epsilon^{3}-4\,\epsilon&11\,\epsilon^{10}+7\,\epsilon^{8}+4\,\epsilon^{6}+8\,\epsilon^{4}+6\,\epsilon^{2}&-13\,\epsilon^{9}-6\,\epsilon^{7}-2\,\epsilon^{5}-3\,\epsilon^{3}\\
       5\,\epsilon^{5}&20\,\epsilon^{6}+4\,\epsilon^{4}&50\,\epsilon^{7}+16\,\epsilon^{5}+3\,\epsilon^{3}&35\,\epsilon^{8}+16\,\epsilon^{6}+6\,\epsilon^{4}+2\,\epsilon^{2}&-5\,\epsilon^{11}-4\,\epsilon^{9}-3\,\epsilon^{7}-14\,\epsilon^{5}-24\,\epsilon^{3}-4\,\epsilon&25\,\epsilon^{10}+16\,\epsilon^{8}+9\,\epsilon^{6}+20\,\epsilon^{4}+16\,\epsilon^{2}+1&-55\,\epsilon^{9}-28\,\epsilon^{7}-12\,\epsilon^{5}-12\,\epsilon^{3}-4\,\epsilon&65\,\epsilon^{8}+24\,\epsilon^{6}+6\,\epsilon^{4}+4\,\epsilon^{2}\\
       -11\,\epsilon^{4}&-44\,\epsilon^{5}-7\,\epsilon^{3}&-110\,\epsilon^{6}-28\,\epsilon^{4}-4\,\epsilon^{2}&-77\,\epsilon^{7}-28\,\epsilon^{5}-8\,\epsilon^{3}-2\,\epsilon&11\,\epsilon^{10}+7\,\epsilon^{8}+4\,\epsilon^{6}+8\,\epsilon^{4}+6\,\epsilon^{2}&-55\,\epsilon^{9}-28\,\epsilon^{7}-12\,\epsilon^{5}-12\,\epsilon^{3}-4\,\epsilon&121\,\epsilon^{8}+49\,\epsilon^{6}+16\,\epsilon^{4}+8\,\epsilon^{2}+1&-143\,\epsilon^{7}-42\,\epsilon^{5}-8\,\epsilon^{3}-2\,\epsilon\\
       13\,\epsilon^{3}&52\,\epsilon^{4}+6\,\epsilon^{2}&130\,\epsilon^{5}+24\,\epsilon^{3}+2\,\epsilon&91\,\epsilon^{6}+24\,\epsilon^{4}+4\,\epsilon^{2}&-13\,\epsilon^{9}-6\,\epsilon^{7}-2\,\epsilon^{5}-3\,\epsilon^{3}&65\,\epsilon^{8}+24\,\epsilon^{6}+6\,\epsilon^{4}+4\,\epsilon^{2}&-143\,\epsilon^{7}-42\,\epsilon^{5}-8\,\epsilon^{3}-2\,\epsilon&169\,\epsilon^{6}+36\,\epsilon^{4}+4\,\epsilon^{2}+1
       \end{bmatrix}$.}
    \]
    We use this $1$-parameter family and \Cref{thm:Nonnegative} to check if the following matrix is in $\SS_4^{\geq 0}$:
    \[
    A=\begin{bmatrix}
       0&0&0&2\\
       0&0&0&0\\
       0&0&0&-2\\
       -2&0&2&0
       \end{bmatrix}.
    \]
       All the minors $M_{j,k}(A)$ are nonnegative. However, 
       \begin{alignat*}{3}
        & M_{1,1}(B(\epsilon))= 80\epsilon^{5} + o(\epsilon^{5}), \quad M_{1,2}(B(\epsilon)) = 40\epsilon^{4} + o(\epsilon^{4}) , \quad M_{1,3}(B(\epsilon))= 16\epsilon^{2} + o(\epsilon^{3}) , \\
        & M_{2,2}(B(\epsilon)) = 80\epsilon^{5} + o(\epsilon^{5}) , \quad M_{2,3}(B(\epsilon)) =\mathbf{-16}\,\epsilon^{2} + o(\epsilon^{2}), \\
        & M_{3,3}(B(\epsilon)) = 2-8\,\epsilon + o(\epsilon).
        \end{alignat*}
        In particular the leading coefficient of $M_{2,3}(B(\epsilon))$ is $-16$, which implies $A \notin \SS_4^{\geq 0}$. We provide Macaulay2 \cite{M2} code to carry out this calculation in \cite{M2Code}.
\end{example}

\subsection{Expressions in the Weyl group}

In this subsection we go into more details about concepts we already introduced in \Cref{subsec:LieTheory}. Let $G$ be a reductive algebraic group and $T$ any maximal torus. Let $W=N_G(T)/T$ be its Weyl group. We will assume basic familiarity with the definition of Weyl groups and we refer the reader to \cite{HumphreysCoxeter} for further details. Let $s_1, \dots, s_n$ be simple reflections that generate $W$. We will reserve underlined letters for expressions (i.e. words) in the simple reflections. Given an expression $\underline{w}=s_{i_1}\cdots s_{i_l}$, we denote by the non-underlined letter $w$ the element $s_{i_1}\cdots s_{i_l}\in W$. We say an expression $\underline{w}=s_{i_1}\cdots s_{i_l}$ for $w\in W$ has length $\ell(\underline{w})=l$. The expression is \emph{reduced} if it is an expression for $w$ of minimal length and we write $\ell(w)$ for the length of a reduced expression for $w$. As before, we denote by $w_0$ for the element of largest length in $W$.

\begin{definition}
    For any subset $J\subset[n]$, we define the \textit{parabolic subgroup} $W_J$ of $W$ by $W_J\coloneqq\langle s_i\mid i\in J\rangle$. We write $W^J \coloneqq W_J \backslash W$ for the set of right cosets of $W_J$. For each coset $A\in W^J$, there is a unique $w\in W$ of minimal length with $A = W_J \cdot w$ called \textit{minimal coset representative}. Given $w\in W$, we denote its minimal coset representative by $w^J$.
\end{definition}

For a reduced expression $\underline{v}=s_{i_1}\cdots s_{i_p}$ of an element $v\in W$, a \emph{subexpression} $\underline{u}$ of $\underline{v}$ is a choice of either $1$ or $s_{i_j}$ for each $j\in[p]$. We will record this by writing $\underline{u}$ as a string whose $j^{\textnormal{th}}$ entry is either $1$ or $s_{i_j}$. We may interpret the subexpression $\underline{u}$ as an expression for some $u\in W$ by removing the $1$s. We will say $\underline{u}$ is reduced if this expression is reduced. Let $\underline{v}$ be an expression for $v\in W$ and $\underline{u}$ be a subexpression of $\underline{v}$. For $k\geq 0$, we will let $\underline{u}_{(k)}$ be the subexpression of $\underline{v}$ which is identical to $\underline{u}$ in its last $k$ entries and is all $1$s beforehand.

\begin{example}\label{ex:expressions}

    Let $n=3$. We work in the Weyl group of type D generated by the simple transpositions in \eqref{eq:WeylGroupGens}. Consider $\underline{v}=s_1s_2s_3s_1s_2$, an expression for $v=462$ (in one-line notation for signed permutations). An example of a subexpression is $\underline{u}=s_11s_3s_11$, which is an expression for $u=624$. We then have 
    \[
        \underline{u}_{(0)}=\underline{u}_{(1)}=11111,\quad  \underline{u}_{(2)}=111s_11,\quad \underline{u}_{(3)}=\underline{u}_{(4)}=11s_3s_11, \quad \text{and} \quad \underline{u}_{(5)}=\underline{u}.
    \] 
    In terms of Weyl group elements we get
    \[
    u_{(0)} = u_{(1)}=123, \quad u_{(2)} = 213,  \quad  u_{(3)} = u_{(4)}=615, \quad \text{and} \quad  u_{(5)} = u.
    \]
\end{example}
The \emph{Bruhat order} $(W,<)$ is defined by $u<v$ if there is a subexpression for $u$ in some, equivalently any, expression for $v$.
\begin{definition}
    We say that a subexpression $\underline{u}$ of $\underline{v}$ is \textit{distinguished}, denoted $\underline{u} \preceq \underline{v}$, if $u_{(j)}\leq s_{i_{n-j+1}}u_{(j-1)}$ for all $j\in[p]$. That is, if multiplying $u_{(j-1)}$ on the left by $s_{i_{n-j+1}}$ decreases the length of $u_{(j-1)}$, then $\underline{u}$ must contain $s_{i_{n-j+1}}$. 
\end{definition}

\begin{example}
    Continuing \Cref{ex:expressions}, we have $\underline{u}\prec\underline{v}$. However, $\underline{w}=111s_11\npreceq\underline{v}$ because $w_{(5)}=s_1\nleq s_1w_{(4)} = {\rm id}$. This example highlights that, colloquially, distinguished subexpressions are leftmost subexpressions.
\end{example}

\begin{theorem}{\cite[Lemma 3.5]{MR}}
Let $\underline{v}$ be a reduced expression for $v\in W$ and let $u<v$. Then there exists a unique subexpression $\underline{u}$ for $u$ which is reduced and distinguished in $\underline{v}$. We call this subexpression the \textbf{reduced distinguished subexpression} for $u$ in $\underline{v}$. 
    
\end{theorem}

\begin{remark}
    As in \cref{rem:MRParamConvention}, our conventions and notation for distinguished subexpressions differ from \cite{MR} since we use row spans. Accordingly, we work with right cosets of $W$ and our definition of distinguished subexpressions is mirrored from \cite{MR}.
\end{remark}

\subsection{Richardson and Deodhar decompositions}\label{subsec:RichardsonDeodhar}

The boundary structure of the totally nonnegative part of partial flag varieties is combinatorially rich and closely related to the theory of positroids in type A. In this section we describe this boundary structure for arbitrary flag varieties and then specialize to the case of the orthogonal Grassmannian. Our presentation follows \cite{KodamaAndWilliams}.

\smallskip

Recall that $G$ is a split reductive algebraic group over $\mathbb{R}$. Fix a pinning for $G$. The reader unfamiliar with the general theory of split reductive algebraic groups over $\RR$ can continue to work with $G=\SO(2n)$, with the pinning introduced in \Cref{sec:2}. The complete flag variety $G/B$ can be identified with the variety of Borel subgroups $\mathcal{B}$ via conjugation: $gB \in  G/B \leftrightarrow g\cdot B:= gBg^{-1} \in \mathcal{B}$. We can also lift the Weyl group elements $w \in W=N_G(T)/T$ to $G$ using the pinning. For each simple reflection $s_i \in W$, define 
\begin{equation}\label{eq:SiDot}
    \dot{s}_i := \varphi_i \begin{pmatrix}
    0 & -1 \\ 1 & 0
    \end{pmatrix}. 
\end{equation}
For a reduced expression $\underline{w}=s_{i_1}\cdots s_{i_k}$, define $\dot{w}=\dot{s}_{i_1}\cdots \dot{s}_{i_k}$. Then, the Bruhat decomposition of $G$ descends to the complete flag variety. The \emph{Richardson cells} of $G/B$ are intersections of opposite Bruhat cells. Concretely, given $v,w \in W$ we define
\[
    \Rcal_{v,w}:= (B\dot{w} \cdot B) \cap (B^{-}\dot{v} \cdot B).
\]
For $\underline{v},\underline{w}$ expressions of $v,w \in W$, Deodhar \cite{Deodhar85} defined the \emph{Deodhar component} $\Rcal_{\underline{v},\underline{w}}$. They refine the Richardson cells, explicitly for a fixed reduced expression $\underline{w}$ we have
\[
    \Rcal_{v,w} = \bigsqcup_{\underline{v} \prec \underline{w}} \Rcal_{\underline{v},\underline{w}},
\]
where the disjoint union is over all distinguished subexpressions $\underline{v}$ of $\underline{w}$. For each Richardson cell $\Rcal_{v,w}$ there is a unique Deodhar component of maximal dimension, namely when we take $\underline{v}$ to be the unique reduced distinguished subexpression for $v$ in $\underline{w}$, which we denote $\underline{v}^{+}$. We denote by $\mathcal{R}^{>0}_{v,w}$ the positive part $\mathcal{R}_{v,w}\cap (G/B)^{\geq0}$ of the Richardson cell $\mathcal{R}_{v,w}$ and by $\mathcal{R}^{>0}_{\underline{v},\underline{w}}$ the positive part $\mathcal{R}_{\underline{v},\underline{w}}\cap (G/B)^{\geq0}$ of the Deodhar component $\mathcal{R}_{\underline{v},\underline{w}}$.

In each Richardson cell, the only Deodhar component that intersects $(G/B)^{\geq 0}$ is the Deodhar component of maximal dimension; that is, by \cite[Theorem 11.3]{MR},
\[
\Rcal_{v,w}^{>0}=\Rcal^{>0}_{\underline{v}^{+},\underline{w}}.
\]

When $w = w_0$ and $v = {\rm id}$, the component $\Rcal_{id, w_0}=(G/B)^{>0}$, which is dense in $(G/B)^{\geq 0}$. As explained in \cite[Section 3]{KodamaAndWilliams}, partial flag varieties also admit Richardson (resp.\ Deodhar) decompositions. Each Richardson cell (resp.\ Deodhar component) of a partial flag variety is the projection of a Richardson cell (resp.\ Deodhar component) of the complete flag variety. Moreover, restricting Richardson cells and Deodhar components to $(G/B)^{\geq 0}$, we obtain a cell decomposition of the nonnegative part of the flag variety into positive Richardson cells $\mathcal{R}_{v,w}^{>0}$.

\smallskip

We now focus on the concrete case of $\OGr(n,2n)$ and its corresponding parabolic quotient $W^{[n-1]}$. We have a Richardson cell for each pair $v, w^{[n-1]} \in W$, where $w^{[n-1]}$ is a minimal coset representative and $v \leq w^{[n-1]}$. For a fixed expression $\underline{w}^{[n-1]}$ for $w^{[n-1]}$, we have a Deodhar component for each distinguished subexpression $\underline{v}$ of $\underline{w}^{[n-1]}$. Slightly abusing notation, we will also denote by $\Rcal_{v,w^{[n-1]}}$, $\Rcal_{\underline{v},\underline{w}^{[n-1]}}$, and $\Rcal_{v,w^{[n-1]}}^{>0}=\Rcal_{\underline{v}^+,\underline{w}^{[n-1]}}^{>0}$ the Richardson cells, Deodhar components and positive Richardson cells, respectively, in $\OGr(n,2n)$. 

\begin{definition} \label{defn:Jsets}
    For a subexpression $\underline{u}$ of a length $p$ expression $\underline{v}$, we define 
        \begin{align*}
            J_{\underline{u}}^{+}&:=\{k\in[p]\mid u_{(k)}> u_{(k-1)}\},\\
            J_{\underline{u}}^{\circ}&:=\{k\in[p]\mid u_{(k)}=u_{(k-1)}\},\\
            J_{\underline{u}}^{-}&:=\{k\in[p]\mid u_{(k)} < u_{(k-1)}\}.
       \end{align*}
\end{definition}

As we did for the totally positive part in \Cref{sec:2}, we can use Marsh and Rietsch's work to parametrize any Deodhar component in terms of the corresponding pair of subexpressions \cite[Section 5]{MR}. Let $\underline{w}^{[n-1]}=s_{i_1}\cdots s_{i_k}$ be a reduced expression. Then, for any distinguished subexpression $\underline{v}\preceq \underline{w}^{[n-1]}$, we have
\[
    \Rcal_{\underline{v},\ \underline{w}^{[n-1]}} = \left\{ \pi_n(g_1\cdots g_p) : \quad  \begin{array}{ll}
        g_k= x_{i_k}(t_k) \text{ with } t_k \in \RR^* &\text{if } k\in J_{\underline{v}^{+}}^{+}\\
        g_k= s_{i_k}^T &\text{if } k\in J_{\underline{v}^{+}}^{\circ} \\
        g_k= x_{i_k}(a_k)^Ts_{i_k} \text{ with } a_k \in \RR &\text{if } k\in J_{\underline{v}^{+}}^{-} 
    \end{array}\right\}. 
\]
The restriction to $(G/B)^{\geq 0}$ is empty unless $\underline{v}=\underline{v}^+$ is reduced distinguished, in which case
\begin{equation}\label{eq:Deodharparam}
    \Rcal^{>0}_{\underline{v}^{+},\ \underline{w}^{[n-1]}} = \left\{ \pi_n(g_1\cdots g_p) : \quad  \begin{array}{ll}
        g_k= x_{i_k}(t_k) \text{ with } t_k \in \RR_{>0} &\text{if } k\in J_{\underline{v}^{+}}^{+}\\
        g_k= s_{i_k}^T &\text{if } k\in J_{\underline{v}^{+}}^{\circ}
    \end{array}\right\}.
\end{equation}

\subsection{LGV Diagrams for Boundary Components}

We have worked extensively with the LGV diagram to view Pl\"ucker coordinates in terms of non-intersecting weighted path collections. The LGV diagram we constructed in \Cref{def:LGVdiagram} is tailored to the top cell of the orthogonal Grassmannian, $\mathcal{R}_{{\rm id}, w_0^{[n-1]}}=\OGr^{>0}(n,2n)$. Here, we explain how to extend this construction to arbitrary Richardson cells in the boundary of $\OGr^{\geq 0}(n,2n)$. We begin by defining a helpful manipulation of certain directed graphs. 

\begin{definition}\label{def:signedaction}
    Let $G$ be a graph with $2n$ horizontal labeled strands directed rightwards, and possibly some vertical arrows between strands. Refer to the $k^{\text{th}}$ strand from the bottom as strand $k$. For each $i\in [n]$, we define the \textit{signed $s_i$ action} on $G$ to be the following sequential operations, mapping $G$ to $s_i\cdot G$.

    \begin{enumerate}
        \item\label{signedactionstep1} Let $s_i=(i_1, i_2)(i_3, i_4)$. Permute strands $i_1$ and $i_2$, and also strands $i_3$ and $i_4$. The arrows move with the strands: if an arrow originates on strand $i_1$ in $G$ then it should originate on strand $i_2$ in $s_i\cdot G$.

        \item\label{signedactionstep2} If $i\in [n-1]$, let $\alpha$ be the rightmost vertical arrow incident to strand $i+1$ in $G$. Add a weight of $-1$ to the portion of strand $i+1$ right of $\alpha$. Do the same for strand $n-i+1$. If $i=n$, do the same for strands $n$ and $n+1$.
    \end{enumerate}
\end{definition}

\begin{definition}\label{def:boundaryLGVdiagram}
    We define a weighted directed graph $G_{v,w}$ for $v\leq w$ and $w$ a minimal coset representative. This definition is technical and is most easily understood by an example (see \Cref{ex:generalLGVdiagram}). Let $\underline{w}^{[n-1]}_0 = s_{i_1}\cdots s_{i_N}$ be as in \eqref{eq:w_0ReducedExpr}, with $N=\binom{n}{2}$. Let $\underline{w}^{[n-1]}=\rho_{1}\cdots \rho_{N}$ be the reduced distinguished expression for $w^{[n-1]}$ in $\underline{w}_0^{[n-1]}$, where each $\rho_{j}$ is either $s_{i_j}$ or $1$. Let $\underline{v}=r_{1}\cdots r_{l}$ be the reduced distinguished subexpression for $v$ in $\underline{w}^{[n-1]}$, where each $r_{j}$ is either $s_{i_j}$ or $1$. In particular, for $j\in [N]$, $\rho_j=1$ implies $r_j=1$. We start from the LGV diagram. Recall that as in \Cref{def:LGVdiagram}, each simple reflection in $\underline{w}^{[n-1]}_0$ corresponds to a pair of arrows in the LGV diagram. 
    \begin{enumerate}
        \item \label{LGVstep1}For each $j\in [N]$ such that $\rho_j=1$, delete the arrows corresponding to $s_{i_j}$.
        \smallskip
        \item \label{LGVstep2}For each $j\in [N]$ such that $r_j\neq 1$, in increasing order of $j$, cut the graph in half with a vertical split between the arrows corresponding to $s_{i_j}$ and the closest arrow to the right of them (if there are no more arrows, then anywhere to the right of them). In the left piece of the graph, remove the arrows corresponding to $s_{i_j}$. Viewing the strands as being labeled by the source vertex, preform a signed $s_{i_j}$ action. Then reattach the two sides of the graph.
        
    \end{enumerate}
\end{definition}

\begin{example}\label{ex:generalLGVdiagram}
    Let us explicitly construct $G_{v,w}$ for $w=s_4s_2s_1s_3,\, v=s_1$. Recall the LGV diagram for $n=4$ from \Cref{ex:n=4LGV}. Note that the reduced distinguished subexpression $\underline{w}^+ \prec \underline{w}_0^{[n-1]}$ is $s_4s_2s_1s_311$ and the reduced distinguished subexpression $\underline{v}^+\prec \underline{w}$ is $11s_1111$. Observe that $\rho_5=\rho_6=1$, so step \ref{LGVstep1} of the construction of $G_{v,w}$ is to remove the last four arrows in the LGV diagram.

     
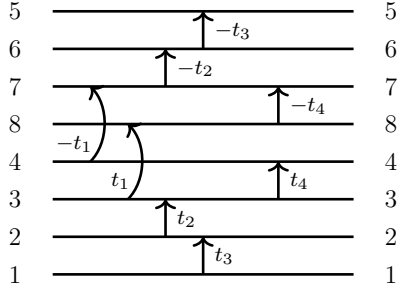
\begin{center}
    \begin{tikzpicture}  
    \coordinate (l7) at (0,4);
    \coordinate (l8) at (0,3.5);
    \coordinate (l9) at (0,3);
    \coordinate (l10) at (0,2.5);        
    \coordinate (l5) at (0,2);        
    \coordinate (l4) at (0,1.5);    
    \coordinate (l3) at (0,1);
    \coordinate (l2) at (0,0.5);
       
    \coordinate (r7) at (4,4);
    \coordinate (r8) at (4,3.5);
    \coordinate (r9) at (4,3);
    \coordinate (r10) at (4,2.5);
    \coordinate (r5) at (4,2);        
    \coordinate (r4) at (4,1.5);    
    \coordinate (r3) at (4,1);
    \coordinate (r2) at (4,0.5);

    \coordinate (v21) at (2,0.5);

    \coordinate (v31) at (1.5,1);
    \coordinate (v32) at (2,1);
    \coordinate (v33) at (3.5,1);
    
    \coordinate (v41) at (1,1.5);
    \coordinate (v42) at (1.5,1.5);
    \coordinate (v43) at (3,1.5);
    \coordinate (v44) at (3.5,1.5);
    \coordinate (v45) at (5,1.5);
    
    \coordinate (v51) at (0.5,2);
    \coordinate (v52) at (3,2);
    \coordinate (v53) at (4.5,2);

    \coordinate (v101) at (1,2.5);    
    \coordinate (v102) at (3,2.5);
    \coordinate (v103) at (5,2.5);

    \coordinate (v91) at (0.5,3);
    \coordinate (v92) at (1.5,3);
    \coordinate (v93) at (3,3);
    \coordinate (v94) at (3.5,3);
    \coordinate (v95) at (4.5,3);

    \coordinate (v81) at (1.5,3.5);
    \coordinate (v82) at (2,3.5);
    \coordinate (v83) at (3.5,3.5);
    \coordinate (v84) at (4,3.5);

    \coordinate (v71) at (2,4);

    \draw[black] (-0.5,0.5) node  [xscale = 0.8, yscale = 0.8] {$1$};
    \draw[black] (-0.5,1) node  [xscale = 0.8, yscale = 0.8] {$2$};
    \draw[black] (-0.5,1.5) node  [xscale = 0.8, yscale = 0.8] {$3$};
    \draw[black] (-0.5,2) node  [xscale = 0.8, yscale = 0.8] {$4$};
    \draw[black] (-0.5,2.5) node  [xscale = 0.8, yscale = 0.8] {$8$};
    \draw[black] (-0.5,3) node  [xscale = 0.8, yscale = 0.8] {$7$};
    \draw[black] (-0.5,3.5) node  [xscale = 0.8, yscale = 0.8] {$6$};
    \draw[black] (-0.5,4) node  [xscale = 0.8, yscale = 0.8] {$5$};
    
    \draw[black] (4.5,0.5) node  [xscale = 0.8, yscale = 0.8] {$1$};
    \draw[black] (4.5,1) node  [xscale = 0.8, yscale = 0.8] {$2$};
    \draw[black] (4.5,1.5) node  [xscale = 0.8, yscale = 0.8] {$3$};
    \draw[black] (4.5,2) node  [xscale = 0.8, yscale = 0.8] {$4$};
    \draw[black] (4.5,2.5) node  [xscale = 0.8, yscale = 0.8] {$8$};
    \draw[black] (4.5,3) node  [xscale = 0.8, yscale = 0.8] {$7$};
    \draw[black] (4.5,3.5) node  [xscale = 0.8, yscale = 0.8] {$6$};
    \draw[black] (4.5,4) node  [xscale = 0.8, yscale = 0.8] {$5$};

    \draw[draw=black, line width=1pt] (l2) -- (r2);
    \draw[draw=black, line width=1pt] (l3) -- (r3);
    \draw[draw=black, line width=1pt] (l4) -- (r4);
    \draw[draw=black, line width=1pt] (l5) -- (r5);
    \draw[draw=black, line width=1pt] (l10) -- (r10);
    \draw[draw=black, line width=1pt] (l9) -- (r9);
    \draw[draw=black, line width=1pt] (l8) -- (r8);
    \draw[draw=black, line width=1pt] (l7) -- (r7);

    \draw[draw=black, line width=1pt, ->]  (v41) .. controls (1.25,1.75) and (1.25,2.25) .. (v101) node[midway, below left] {\tiny $t_{1}$};
    
    \draw[draw=black, line width=1pt, ->]  (v51) .. controls (0.75,2.25) and (0.75,2.75) .. (v91) node[midway, below left] {\tiny $-t_{1}$};
    
    \draw[draw=black, line width=1pt, ->]  (v31) -- (v42) node[below right] {\tiny $t_{2}$};
    \draw[draw=black, line width=1pt, ->]  (v92) -- (v81) node[below right] {\tiny $-t_{2}$};
    
    \draw[draw=black, line width=1pt, ->]  (v21) -- (v32) node[below right] {\tiny $t_{3}$};
    \draw[draw=black, line width=1pt, ->]  (v82) -- (v71) node[below right] {\tiny $-t_{3}$};

    \draw[draw=black, line width=1pt, ->]  (v43) -- (v52) node[below right] {\tiny $t_{4}$};
    \draw[draw=black, line width=1pt, ->]  (v102) -- (v93) node[below right] {\tiny $-t_{4}$};

\end{tikzpicture}
    \end{center}

    We next apply Step \ref{LGVstep2}. The only $j\in[N]$ for which $r_j\neq 1$ is $j=3$. We begin by dividing the graph above into two pieces, between the third and fourth pairs of arrows.

     
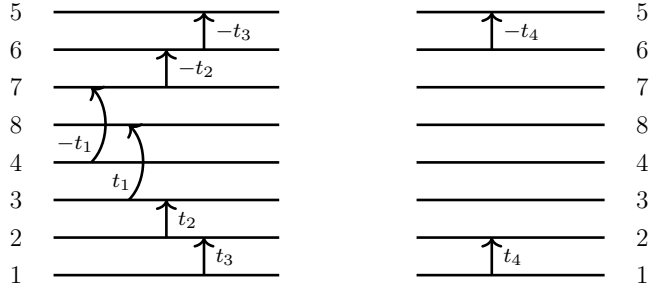
\begin{center}
    \begin{tikzpicture}  
    \coordinate (l7) at (0,4);
    \coordinate (l8) at (0,3.5);
    \coordinate (l9) at (0,3);
    \coordinate (l10) at (0,2.5);        
    \coordinate (l5) at (0,2);        
    \coordinate (l4) at (0,1.5);    
    \coordinate (l3) at (0,1);
    \coordinate (l2) at (0,0.5);
       
    \coordinate (r7) at (3,4);
    \coordinate (r8) at (3,3.5);
    \coordinate (r9) at (3,3);
    \coordinate (r10) at (3,2.5);
    \coordinate (r5) at (3,2);        
    \coordinate (r4) at (3,1.5);    
    \coordinate (r3) at (3,1);
    \coordinate (r2) at (3,0.5);

    \coordinate (v21) at (2,0.5);

    \coordinate (v31) at (1.5,1);
    \coordinate (v32) at (2,1);
    \coordinate (v33) at (3.5,1);
    
    \coordinate (v41) at (1,1.5);
    \coordinate (v42) at (1.5,1.5);
    \coordinate (v43) at (3,1.5);
    \coordinate (v44) at (3.5,1.5);
    \coordinate (v45) at (5,1.5);
    
    \coordinate (v51) at (0.5,2);
    \coordinate (v52) at (3,2);
    \coordinate (v53) at (4.5,2);

    \coordinate (v101) at (1,2.5);    
    \coordinate (v102) at (3,2.5);
    \coordinate (v103) at (5,2.5);

    \coordinate (v91) at (0.5,3);
    \coordinate (v92) at (1.5,3);
    \coordinate (v93) at (3,3);
    \coordinate (v94) at (3.5,3);
    \coordinate (v95) at (4.5,3);

    \coordinate (v81) at (1.5,3.5);
    \coordinate (v82) at (2,3.5);
    \coordinate (v83) at (3.5,3.5);
    \coordinate (v84) at (4,3.5);

    \coordinate (v71) at (2,4);

    \draw[black] (-0.5,0.5) node  [xscale = 0.8, yscale = 0.8] {$1$};
    \draw[black] (-0.5,1) node  [xscale = 0.8, yscale = 0.8] {$2$};
    \draw[black] (-0.5,1.5) node  [xscale = 0.8, yscale = 0.8] {$3$};
    \draw[black] (-0.5,2) node  [xscale = 0.8, yscale = 0.8] {$4$};
    \draw[black] (-0.5,2.5) node  [xscale = 0.8, yscale = 0.8] {$8$};
    \draw[black] (-0.5,3) node  [xscale = 0.8, yscale = 0.8] {$7$};
    \draw[black] (-0.5,3.5) node  [xscale = 0.8, yscale = 0.8] {$6$};
    \draw[black] (-0.5,4) node  [xscale = 0.8, yscale = 0.8] {$5$};

    \draw[draw=black, line width=1pt] (l2) -- (r2);
    \draw[draw=black, line width=1pt] (l3) -- (r3);
    \draw[draw=black, line width=1pt] (l4) -- (r4);
    \draw[draw=black, line width=1pt] (l5) -- (r5);
    \draw[draw=black, line width=1pt] (l10) -- (r10);
    \draw[draw=black, line width=1pt] (l9) -- (r9);
    \draw[draw=black, line width=1pt] (l8) -- (r8);
    \draw[draw=black, line width=1pt] (l7) -- (r7);

    \draw[draw=black, line width=1pt, ->]  (v41) .. controls (1.25,1.75) and (1.25,2.25) .. (v101) node[midway, below left] {\tiny $t_{1}$};
    
    \draw[draw=black, line width=1pt, ->]  (v51) .. controls (0.75,2.25) and (0.75,2.75) .. (v91) node[midway, below left] {\tiny $-t_{1}$};
    
    \draw[draw=black, line width=1pt, ->]  (v31) -- (v42) node[below right] {\tiny $t_{2}$};
    \draw[draw=black, line width=1pt, ->]  (v92) -- (v81) node[below right] {\tiny $-t_{2}$};
    
    \draw[draw=black, line width=1pt, ->]  (v21) -- (v32) node[below right] {\tiny $t_{3}$};
    \draw[draw=black, line width=1pt, ->]  (v82) -- (v71) node[below right] {\tiny $-t_{3}$};
\end{tikzpicture} \qquad \qquad
\begin{tikzpicture}  
    \coordinate (l7) at (0,4);
    \coordinate (l8) at (0,3.5);
    \coordinate (l9) at (0,3);
    \coordinate (l10) at (0,2.5);        
    \coordinate (l5) at (0,2);        
    \coordinate (l4) at (0,1.5);    
    \coordinate (l3) at (0,1);
    \coordinate (l2) at (0,0.5);
       
    \coordinate (r7) at (2.5,4);
    \coordinate (r8) at (2.5,3.5);
    \coordinate (r9) at (2.5,3);
    \coordinate (r10) at (2.5,2.5);
    \coordinate (r5) at (2.5,2);        
    \coordinate (r4) at (2.5,1.5);    
    \coordinate (r3) at (2.5,1);
    \coordinate (r2) at (2.5,0.5);

    \coordinate (v1) at (1,0.5);
    \coordinate (v2) at (1,1);
    \coordinate (v3) at (1,3.5);
    \coordinate (v4) at (1,4);

    \draw[black] (3,0.5) node  [xscale = 0.8, yscale = 0.8] {$1$};
    \draw[black] (3,1) node  [xscale = 0.8, yscale = 0.8] {$2$};
    \draw[black] (3,1.5) node  [xscale = 0.8, yscale = 0.8] {$3$};
    \draw[black] (3,2) node  [xscale = 0.8, yscale = 0.8] {$4$};
    \draw[black] (3,2.5) node  [xscale = 0.8, yscale = 0.8] {$8$};
    \draw[black] (3,3) node  [xscale = 0.8, yscale = 0.8] {$7$};
    \draw[black] (3,3.5) node  [xscale = 0.8, yscale = 0.8] {$6$};
    \draw[black] (3,4) node  [xscale = 0.8, yscale = 0.8] {$5$};

    \draw[draw=black, line width=1pt] (l2) -- (r2);
    \draw[draw=black, line width=1pt] (l3) -- (r3);
    \draw[draw=black, line width=1pt] (l4) -- (r4);
    \draw[draw=black, line width=1pt] (l5) -- (r5);
    \draw[draw=black, line width=1pt] (l10) -- (r10);
    \draw[draw=black, line width=1pt] (l9) -- (r9);
    \draw[draw=black, line width=1pt] (l8) -- (r8);
    \draw[draw=black, line width=1pt] (l7) -- (r7);

    \draw[draw=black, line width=1pt, ->]  (v1) -- (v2) node[below right] {\tiny $t_{4}$};
    \draw[draw=black, line width=1pt, ->]  (v3) -- (v4) node[below right] {\tiny $-t_{4}$};

\end{tikzpicture}
    \end{center}

We remove the arrows with weights $\pm t_3$ and then apply the signed $s_1$ action to the left piece, interchanging strands $1$ and $2$, as well as strands $5$ and $6$. Note that the arrow with weight $t_2$ originates from the second strand from the bottom before applying the signed $s_1$ action but, after permuting strands $1$ and $2$, it originates on the bottom strand.  We use red to indicate the parts of strands which get a weight of $-1$. 
    
     
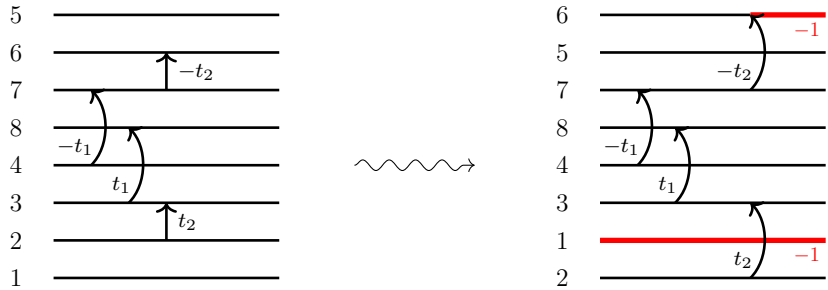
\begin{center}
    \begin{tikzpicture}  
    \coordinate (l7) at (0,4);
    \coordinate (l8) at (0,3.5);
    \coordinate (l9) at (0,3);
    \coordinate (l10) at (0,2.5);        
    \coordinate (l5) at (0,2);        
    \coordinate (l4) at (0,1.5);    
    \coordinate (l3) at (0,1);
    \coordinate (l2) at (0,0.5);
       
    \coordinate (r7) at (3,4);
    \coordinate (r8) at (3,3.5);
    \coordinate (r9) at (3,3);
    \coordinate (r10) at (3,2.5);
    \coordinate (r5) at (3,2);        
    \coordinate (r4) at (3,1.5);    
    \coordinate (r3) at (3,1);
    \coordinate (r2) at (3,0.5);

    \coordinate (v21) at (2,0.5);

    \coordinate (v31) at (1.5,1);
    \coordinate (v32) at (2,1);
    \coordinate (v33) at (3.5,1);
    
    \coordinate (v41) at (1,1.5);
    \coordinate (v42) at (1.5,1.5);
    \coordinate (v43) at (3,1.5);
    \coordinate (v44) at (3.5,1.5);
    \coordinate (v45) at (5,1.5);
    
    \coordinate (v51) at (0.5,2);
    \coordinate (v52) at (3,2);
    \coordinate (v53) at (4.5,2);

    \coordinate (v101) at (1,2.5);    
    \coordinate (v102) at (3,2.5);
    \coordinate (v103) at (5,2.5);

    \coordinate (v91) at (0.5,3);
    \coordinate (v92) at (1.5,3);
    \coordinate (v93) at (3,3);
    \coordinate (v94) at (3.5,3);
    \coordinate (v95) at (4.5,3);

    \coordinate (v81) at (1.5,3.5);
    \coordinate (v82) at (2,3.5);
    \coordinate (v83) at (3.5,3.5);
    \coordinate (v84) at (4,3.5);

    \coordinate (v71) at (2,4);

    \draw[black] (-0.5,0.5) node  [xscale = 0.8, yscale = 0.8] {$1$};
    \draw[black] (-0.5,1) node  [xscale = 0.8, yscale = 0.8] {$2$};
    \draw[black] (-0.5,1.5) node  [xscale = 0.8, yscale = 0.8] {$3$};
    \draw[black] (-0.5,2) node  [xscale = 0.8, yscale = 0.8] {$4$};
    \draw[black] (-0.5,2.5) node  [xscale = 0.8, yscale = 0.8] {$8$};
    \draw[black] (-0.5,3) node  [xscale = 0.8, yscale = 0.8] {$7$};
    \draw[black] (-0.5,3.5) node  [xscale = 0.8, yscale = 0.8] {$6$};
    \draw[black] (-0.5,4) node  [xscale = 0.8, yscale = 0.8] {$5$};

    \path[draw, ->, decorate, decoration ={snake, amplitude = 2}] (4,2) -- (5.6,2);
    
    \draw[draw=black, line width=1pt] (l2) -- (r2);
    \draw[draw=black, line width=1pt] (l3) -- (r3);
    \draw[draw=black, line width=1pt] (l4) -- (r4);
    \draw[draw=black, line width=1pt] (l5) -- (r5);
    \draw[draw=black, line width=1pt] (l10) -- (r10);
    \draw[draw=black, line width=1pt] (l9) -- (r9);
    \draw[draw=black, line width=1pt] (l8) -- (r8);
    \draw[draw=black, line width=1pt] (l7) -- (r7);

    \draw[draw=black, line width=1pt, ->]  (v41) .. controls (1.25,1.75) and (1.25,2.25) .. (v101) node[midway, below left] {\tiny $t_{1}$};
    
    \draw[draw=black, line width=1pt, ->]  (v51) .. controls (0.75,2.25) and (0.75,2.75) .. (v91) node[midway, below left] {\tiny $-t_{1}$};
    
    \draw[draw=black, line width=1pt, ->]  (v31) -- (v42) node[below right] {\tiny $t_{2}$};
    \draw[draw=black, line width=1pt, ->]  (v92) -- (v81) node[below right] {\tiny $-t_{2}$};
    
\end{tikzpicture}  \qquad 
 \begin{tikzpicture}  
    \coordinate (l7) at (0,4);
    \coordinate (l8) at (0,3.5);
    \coordinate (l9) at (0,3);
    \coordinate (l10) at (0,2.5);        
    \coordinate (l5) at (0,2);        
    \coordinate (l4) at (0,1.5);    
    \coordinate (l3) at (0,1);
    \coordinate (l2) at (0,0.5);
       
    \coordinate (r7) at (3,4);
    \coordinate (r8) at (3,3.5);
    \coordinate (r9) at (3,3);
    \coordinate (r10) at (3,2.5);
    \coordinate (r5) at (3,2);        
    \coordinate (r4) at (3,1.5);    
    \coordinate (r3) at (3,1);
    \coordinate (r2) at (3,0.5);

    \coordinate (v21) at (2,0.5);

    \coordinate (v31) at (1.5,1);
    \coordinate (v32) at (2,1);
    \coordinate (v33) at (3.5,1);
    
    \coordinate (v41) at (1,1.5);
    \coordinate (v42) at (1.5,1.5);
    \coordinate (v43) at (3,1.5);
    \coordinate (v44) at (3.5,1.5);
    \coordinate (v45) at (5,1.5);
    
    \coordinate (v51) at (0.5,2);
    \coordinate (v52) at (3,2);
    \coordinate (v53) at (4.5,2);

    \coordinate (v101) at (1,2.5);    
    \coordinate (v102) at (3,2.5);
    \coordinate (v103) at (5,2.5);

    \coordinate (v91) at (0.5,3);
    \coordinate (v92) at (1.5,3);
    \coordinate (v93) at (3,3);
    \coordinate (v94) at (3.5,3);
    \coordinate (v95) at (4.5,3);

    \coordinate (v81) at (1.5,3.5);
    \coordinate (v82) at (2,3.5);
    \coordinate (v83) at (3.5,3.5);
    \coordinate (v84) at (4,3.5);

    \coordinate (v71) at (2,4);

    \coordinate (vmod1) at (2,1.5);
    \coordinate (vmod2) at (2,3);

    \draw[black] (-0.5,0.5) node  [xscale = 0.8, yscale = 0.8] {$2$};
    \draw[black] (-0.5,1) node  [xscale = 0.8, yscale = 0.8] {$1$};
    \draw[black] (-0.5,1.5) node  [xscale = 0.8, yscale = 0.8] {$3$};
    \draw[black] (-0.5,2) node  [xscale = 0.8, yscale = 0.8] {$4$};
    \draw[black] (-0.5,2.5) node  [xscale = 0.8, yscale = 0.8] {$8$};
    \draw[black] (-0.5,3) node  [xscale = 0.8, yscale = 0.8] {$7$};
    \draw[black] (-0.5,3.5) node  [xscale = 0.8, yscale = 0.8] {$5$};
    \draw[black] (-0.5,4) node  [xscale = 0.8, yscale = 0.8] {$6$};

    \draw[red] (2.75,0.8) node [xscale=0.65,yscale=0.65]{$-1$};
    \draw[red] (2.75,3.8) node [xscale=0.65,yscale=0.65]{$-1$};

    \draw[draw=black, line width=1pt] (l2) -- (r2);
    \draw[draw=red, line width=2pt] (l3) -- (r3);
    \draw[draw=black, line width=1pt] (l4) -- (r4);
    \draw[draw=black, line width=1pt] (l5) -- (r5);
    \draw[draw=black, line width=1pt] (l10) -- (r10);
    \draw[draw=black, line width=1pt] (l9) -- (r9);
    \draw[draw=black, line width=1pt] (l8) -- (r8);
    \draw[draw=black, line width=1pt] (l7) -- (v71);
    \draw[draw=red, line width=2pt] (v71) -- (r7);

    \draw[draw=black, line width=1pt, ->]  (v41) .. controls (1.25,1.75) and (1.25,2.25) .. (v101) node[midway, below left] {\tiny $t_{1}$};
    
    \draw[draw=black, line width=1pt, ->]  (v51) .. controls (0.75,2.25) and (0.75,2.75) .. (v91) node[midway, below left] {\tiny $-t_{1}$};
    
    \draw[draw=black, line width=1pt, ->]  (v21)  .. controls (2.25,0.75) and (2.25,1.25) .. (vmod1) node[midway, below left] {\tiny $t_{2}$};
    \draw[draw=black, line width=1pt, ->]  (vmod2) .. controls (2.25,3.25) and (2.25,3.75) .. (v71) node[midway, below left] {\tiny $-t_{2}$};
\end{tikzpicture}
    \end{center}

Finally we reattach the two pieces of the graph to obtain $G_{v,w}$.

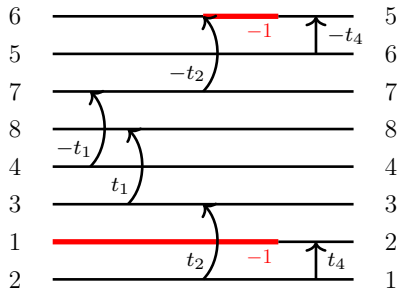
     \begin{center}

 \begin{tikzpicture}  
    \coordinate (l7) at (0,4);
    \coordinate (l8) at (0,3.5);
    \coordinate (l9) at (0,3);
    \coordinate (l10) at (0,2.5);        
    \coordinate (l5) at (0,2);        
    \coordinate (l4) at (0,1.5);    
    \coordinate (l3) at (0,1);
    \coordinate (l2) at (0,0.5);
       
    \coordinate (r7') at (3,4);
    \coordinate (r7) at (4,4);
    \coordinate (r8) at (4,3.5);
    \coordinate (r9) at (4,3);
    \coordinate (r10) at (4,2.5);
    \coordinate (r5) at (4,2);        
    \coordinate (r4) at (4,1.5);    
    \coordinate (r3') at (3,1);
    \coordinate (r3) at (4,1);
    \coordinate (r2) at (4,0.5);

    \coordinate (v21) at (2,0.5);

    \coordinate (v31) at (1.5,1);
    \coordinate (v32) at (2,1);
    \coordinate (v33) at (3.5,1);
    
    \coordinate (v41) at (1,1.5);
    \coordinate (v42) at (1.5,1.5);
    \coordinate (v43) at (3,1.5);
    \coordinate (v44) at (3.5,1.5);
    \coordinate (v45) at (5,1.5);
    
    \coordinate (v51) at (0.5,2);
    \coordinate (v52) at (3,2);
    \coordinate (v53) at (4.5,2);

    \coordinate (v101) at (1,2.5);    
    \coordinate (v102) at (3,2.5);
    \coordinate (v103) at (5,2.5);

    \coordinate (v91) at (0.5,3);
    \coordinate (v92) at (1.5,3);
    \coordinate (v93) at (3,3);
    \coordinate (v94) at (3.5,3);
    \coordinate (v95) at (4.5,3);

    \coordinate (v81) at (1.5,3.5);
    \coordinate (v82) at (2,3.5);
    \coordinate (v83) at (3.5,3.5);
    \coordinate (v84) at (4,3.5);

    \coordinate (v71) at (2,4);

    \coordinate (vmod1) at (2,1.5);
    \coordinate (vmod2) at (2,3);

    \coordinate (vmod3) at (3.5,0.5);
    \coordinate (vmod4) at (3.5,1);
    \coordinate (vmod5) at (3.5,3.5);
    \coordinate (vmod6) at (3.5,4);

    \draw[black] (-0.5,0.5) node  [xscale = 0.8, yscale = 0.8] {$2$};
    \draw[black] (-0.5,1) node  [xscale = 0.8, yscale = 0.8] {$1$};
    \draw[black] (-0.5,1.5) node  [xscale = 0.8, yscale = 0.8] {$3$};
    \draw[black] (-0.5,2) node  [xscale = 0.8, yscale = 0.8] {$4$};
    \draw[black] (-0.5,2.5) node  [xscale = 0.8, yscale = 0.8] {$8$};
    \draw[black] (-0.5,3) node  [xscale = 0.8, yscale = 0.8] {$7$};
    \draw[black] (-0.5,3.5) node  [xscale = 0.8, yscale = 0.8] {$5$};
    \draw[black] (-0.5,4) node  [xscale = 0.8, yscale = 0.8] {$6$};

    \draw[black] (4.5,0.5) node  [xscale = 0.8, yscale = 0.8] {$1$};
    \draw[black] (4.5,1) node  [xscale = 0.8, yscale = 0.8] {$2$};
    \draw[black] (4.5,1.5) node  [xscale = 0.8, yscale = 0.8] {$3$};
    \draw[black] (4.5,2) node  [xscale = 0.8, yscale = 0.8] {$4$};
    \draw[black] (4.5,2.5) node  [xscale = 0.8, yscale = 0.8] {$8$};
    \draw[black] (4.5,3) node  [xscale = 0.8, yscale = 0.8] {$7$};
    \draw[black] (4.5,3.5) node  [xscale = 0.8, yscale = 0.8] {$6$};
    \draw[black] (4.5,4) node  [xscale = 0.8, yscale = 0.8] {$5$};

    \draw[red] (2.75,0.8) node [xscale=0.65,yscale=0.65]{$-1$};
    \draw[red] (2.75,3.8) node [xscale=0.65,yscale=0.65]{$-1$};

    \draw[draw=black, line width=1pt] (l2) -- (r2);
    \draw[draw=red, line width=2pt] (l3) -- (r3');
    \draw[draw=black, line width=1pt] (r3') -- (r3);
    \draw[draw=black, line width=1pt] (l4) -- (r4);
    \draw[draw=black, line width=1pt] (l5) -- (r5);
    \draw[draw=black, line width=1pt] (l10) -- (r10);
    \draw[draw=black, line width=1pt] (l9) -- (r9);
    \draw[draw=black, line width=1pt] (l8) -- (r8);
    \draw[draw=black, line width=1pt] (l7) -- (v71);
    \draw[draw=red, line width=2pt] (v71) -- (r7');
    \draw[draw=black, line width=1pt] (r7') -- (r7);

    \draw[draw=black, line width=1pt, ->]  (v41) .. controls (1.25,1.75) and (1.25,2.25) .. (v101) node[midway, below left] {\tiny $t_{1}$};
    
    \draw[draw=black, line width=1pt, ->]  (v51) .. controls (0.75,2.25) and (0.75,2.75) .. (v91) node[midway, below left] {\tiny $-t_{1}$};
    
    \draw[draw=black, line width=1pt, ->]  (v21)  .. controls (2.25,0.75) and (2.25,1.25) .. (vmod1) node[midway, below left] {\tiny $t_{2}$};
    \draw[draw=black, line width=1pt, ->]  (vmod2) .. controls (2.25,3.25) and (2.25,3.75) .. (v71) node[midway, below left] {\tiny $-t_{2}$};

    \draw[draw=black, line width=1pt, ->]  (vmod3) -- (vmod4) node[below right] {\tiny $t_{4}$};
    \draw[draw=black, line width=1pt, ->]  (vmod5) -- (vmod6) node[below right] {\tiny $-t_{4}$};

\end{tikzpicture}
    \end{center}

\end{example}

Each arrow in $G_{v,w}$ corresponds to some arrow in the original LGV diagram. We will continue to use the notations $a_i^{(j)}$ and $b_i^{(j)}$ to refer to the arrows in $G_{v,w}$.

\begin{proposition}\label{LGVlemmavw}
    Let $I \in \binom{[2n]}{n}$ and $\mathcal{P}_I$ be the set of non-intersecting path collections from $\{1,\ldots,n\}$ to $I$ in $G_{v,w}$. Then, the maximal minor of $X \in \mathcal{R}_{v,w}^{>0}$ corresponding to the columns indexed by $I$ is given by
    \[
        \Delta^I(X) :=\det\left(X_I\right)= \sum_{P=(P_1,\ldots,P_n) \in \mathcal{P}_I} \sgn(\pi_P)\prod_{i=1}^n \prod_{e \in P_i} w(e),
    \]
    where $\pi_P$ is the permutation that sends $i$ to the end point of $P_i$, and the weight of the edge $e$ is $w(e) \in \{\pm t_1, \ldots, \pm t_{\ell(w)-\ell(v)}\}$.
\end{proposition}

\begin{proof}
    The proof of this statement is identical to the proof of \Cref{LGVlemma}, except that we now look at the parametrization of $\mathcal{R}_{v,w}^{>0}$ given by the positive part of the Deodhar component $\mathcal{R}_{\underline{v}^+,\underline{w}}^{>0}$, as in \eqref{eq:Deodharparam}. Consequently, we now need to account for the matrices $\dot{s}_i$ in addition to the matrices $x_{i_l}(t_l)$, which can also be treated as adjacency matrices.
\end{proof}

\subsection{Identification of cells}

Our goal in this section is to identify the Richardson cell $\mathcal{R}_{v,w}^{>0}$ containing a given point $X\in \OGr^{\geq 0}(n, 2n)$. Let $\mathscr{M}_X$ denote the matroid associated to $X$, that is, the rank $n$ matroid on $[2n]$ with bases 
\[
    \left \{ I \in \binom{[2n]}{n} \colon \quad \Delta^{I}(X) \neq 0\right \}.
\] 
 
In the nonnegative Grassmannian, Richardson cells are in bijection with \emph{positroids}. As we shall see, the Richardson cell that contains a given point $X \in \OGr^{\geq 0}(n,2n)$ is also determined by its matroid $\mathscr{M}_X$. The following definition is motivated by \cite[Definition 3.26]{boretsky}.

\begin{definition}\label{def:lowering}
    Let $I$ be a basis of a rank $k$ matroid $\mathscr{M}$ on a ground set $E$. Let $<$ be a total order on $E$. We inductively define a sequence of bases of $\mathscr{M}$ as follows. Let $I^{(0)}=I$. For $1\leq j\leq k$, suppose $I^{(j-1)}=\left\{i^{(j-1)}_1<\cdots <i^{(j-1)}_k\right\}$. Then, $I^{(j)}$ is obtained by replacing $i^{(j-1)}_{j}$ in $I^{(j-1)}$ by the element
    \[
        \min_{<}\left\{ t \in E \colon \quad \left( I^{(j-1)}\setminus \{i^{(j-1)}_{j}\} \right) \cup \{t\} \text{ is a basis of } \mathscr{M} \right\}.
    \] 
    In words, $I^{(j)}$ is the basis obtained from $I^{(j-1)}$ by decreasing the $j^{\text{th}}$ element of $I^{(j)}$ as much as possible with respect to $<$. For $0\leq j\leq k$, we say $I^{(j)}$ is a \textit{lowering} of $I$ with respect to~$<$.
\end{definition}

Recall that permutations $\sigma$ in the Weyl group $W$ of $\SO(2n)$ satisfy $\sigma(n+j)=n+\sigma(j)\text{ modulo } 2n$ for all $i\in [n]$. Thus, to determine $v$ and $w$, it suffices to specify $v(j)$ and $w(j)$ for all $j\in[n]$. The total order $\prec$ we shall use on $[2n]$ is the following:
\[
    1 \prec 2 \prec \dots \prec n \prec 2n \prec 2n-1 \prec \dots \prec n+1.
\]

Recall that for a fixed total order $<$ of $E$, the \textit{Gale order} on $\binom{E}{k}$, also denoted $<$, is defined~by
\begin{equation*}
    A < B \qquad \text{if} \qquad a_i < b_i \ \text{ for all } i \in [k],
\end{equation*}
where $A=\{a_1<a_2<\cdots < a_k\}$ and $B=\{b_1<b_2<\cdots < b_k\}$.

We make use of the following definition of a matroid: a rank $k$ matroid on a ground set $E$ is a collection of subsets of $E$ of size $k$, that has a unique Gale-maximal element with respect to any total order on $E$ \cite[Theorem 4]{Gale}. We denote by $I_X$ the unique maximal element of the matroid $\mathscr{M}_X$ in the Gale order $\prec$.

\begin{definition}\label{def:Iw}
    For $w\in W$, define $I_w=w^{-1} \big([2n]\setminus[n] \big)\cap [n]=\{i\in[n]\colon w(i)>n\}$. A defining property of the Weyl group $W$ of type D is that $|I_w|$ is even.
\end{definition}

\begin{remark}\label{rem:wAsSet}
    Viewing $w\in W$ as a permutation of $[2n]$ written in one-line notation, we note that multiplying $w$ on the left by $s_i$ for $i<n$ permutes separately the values in $[n]$ and those in $[2n]\setminus[n]$ in $w$. Thus, for $w\in W$, the coset $W_{[n-1]} w$ is determined by the (unordered) positions of $[2n]\setminus [n]$ in $w$. Since $w(i)$ determines $w(i+n)$ for $i\in [n]$, the coset $W_{[n-1]}w$ is determined by $I_w$. The minimal coset representative of $w$ will be the shortest permutation $w'$ satisfying $I_{w'} = I_w$. Specifically, this is the permutation which, in one-line notation, has $1,2,\ldots, n-|I_w|$ in order in positions $[n]\setminus I_w$ and $2n, 2n-1, \ldots, 2n-|I_w|+1$ in positions $I_w$. Thus, the right cosets $W^{[n-1]}$ are in bijection with subsets of $[n]$ of even size, with the bijection given by $w^{[n-1]}\leftrightarrow I_{w^{[n-1]}}$. 
\end{remark}

\begin{lemma}\label{lem:bijectionsetstoexpressions}
    Fix an element $w\in W$. A reduced subexpression $\underline{w}^{[n-1]}$ for the minimal coset representative $w^{[n-1]}$ can be extracted from the expression $\underline{w}_0^{[n-1]}$ in \eqref{eq:w_0ReducedExpr} as follows.

    \begin{enumerate}
        \item Locate the first $|I_w| / 2$ factors of $s_n$ appearing in \eqref{eq:w_0ReducedExpr}. These appear in~$\underline{w}^{[n-1]}$.

        \smallskip
        
        \item If $I_w = \{i_1<i_2<\cdots <i_{|I_w|}\}$, then for each $1\leq j\leq |I_w|$ and $j$ odd (resp. even), the reflections $s_{n-2}\cdots s_{i_{(j+1)}}$ (resp. $s_{n-1}\cdots s_{i_j}$) from the $j^\text{th}$ set of parentheses of \eqref{eq:w_0ReducedExpr} appear in $\underline{w}^{[n-1]}$.

        \smallskip

        \item All the other reflections in \eqref{eq:w_0ReducedExpr} are replaced with $1$.
    \end{enumerate}
\end{lemma}

\begin{proof}
    Consider the subexpression $\underline{u}$ defined in the lemma statement. Our goal is to show that $u=w^{[n-1]}$. Ignoring factors of $1$, we can write it as $\underline{u} = \underline{w}_1 \cdots \underline{w}_{{|I_w|}/{2}}$, where $\underline{w}_j=s_n(s_{n-2}\cdots s_{i_{(2j-1)}})(s_{n-1}\cdots s_{i_{(2j)}})$. Recall that acting on the right, simple transpositions act on positions. Note that $I_{s_n}=\{n-1,n\}$. Applying the sequence of permutations $s_{n-2}\cdots s_{i_1}$ on the right moves the entry at position $n-1$ to position $i_1$. Thus, $I_{s_n(s_{n-2}\cdots s_{i_1})}=\{i_1,n\}$. Similarly, applying the sequence of permutations $s_{n-1}\cdots s_{i_2}$ on the right moves the entry at position $n$ to position $i_2$. Thus, $I_{w_1}=\{i_1,i_2\}$. Also, $w_1$ is the shortest element of $W$ with this property. To see this, observe that for any $u\in W$, $I_{us_n}$ equals the symmetric difference $I_u \ \triangle \{n-1,n\}$, while, for $i\in [n-1]$, $I_{us_i}$ equals $I_w$ with $i$ replaced by $i+1$ and vice versa. By our construction, each simple reflection in $\underline{w}_1$ acts non trivially: the factors of $s_n$ act by adding $\{n-1,n\}$ whereas the other $s_i$ act by swapping $i+1$ for $i$. Thus, it uses the fewest possible simple transpositions to obtain the set $I_{w_1}=\{i_1, i_2\}$. 

    \smallskip
    
    Using a similar argument, one can verify inductively that $I_{w_1\cdots w_j}=\{i_1,\ldots, i_{2j}\}$ and that $w_1\cdots w_j$ is the shortest element of $W$ with this property. Thus, we deduce $u=w^{[n-1]}$, and the given expression is reduced.
\end{proof}

\begin{example}
    Let $n=5$. Then, from \eqref{eq:w_0ReducedExpr}, $\underline{w_0}^{[n-1]}=s_5(s_3s_2s_1)(s_4s_3s_2)s_5(s_3)(s_4)$. Let $w=1\; 10\; 2\; 9\; 3$, which is a shortest coset representative. Then, $I_w=\{2,4\}$. The subexpression for $\underline{w}$ is thus $s_5(s_3s_21)(s_411)1(1)(1)$. 
\end{example}

The following proposition proves a number of helpful technical facts about $G_{v,w}$. While many of its claims seem unrelated, they are all proven in the same induction and so are most easily proven together.

\begin{proposition}\label{prop:pathcollectionproperties}
        
        Let $v\leq w$ in $W$ such that $w$ is a minimal coset representative. In $G_{v,w}$, there is a unique non-intersecting path collection from $[n]$ to $I_X$. Moreover, 
        \begin{enumerate}[wide = 30pt, leftmargin = 50pt]
            \item \label{property1} $I_X$ depends only on $w$,
            \item \label{property2} In this path collection, the path originating from source vertex $i$ terminates below the path originating from source vertex $i+1$ for all $i\in[n-1]$,
            \item \label{property3} The path collection is left greedy, and
            \item \label{property4} The path collection uses every vertical edge of $G_{v,w}$.
        \end{enumerate}
\end{proposition}

\begin{proof}
    Fix $w \in W^{[n-1]}$ and $\underline{w}$ the expression given by \Cref{lem:bijectionsetstoexpressions}. We first show the result for $v = {\rm id}$ and proceed by induction on $\ell(v)$.
    
    \smallskip
    
    Let $v={\rm id}$. We give a graphical interpretation of $I_w$. Let $I'_w$ be the set $I'_w := w^{-1}([2n]\setminus [n])$. Note that $I_w$ determines $I'_w$ since for $w\in W$ and $i\in [n]$, $w(i+n)=w(i)+n\text{ modulo }2n$. Define $G'_{v,w}$ to be identical to $G_{v,w}$ but with unoriented vertical arrows. We consider path collections in $G'_{v,w}$ (not necessarily non-intersecting). We define the greedy path originating at $i\in[2n]$ in $G'_{v,w}$ to be the path which travels rightwards along strands and greedily turns onto every vertical edge it passes. The difference from a left greedy path in $G_{v,w}$ is that a greedy path in $G'_{v,w}$ can travel either way along vertical edges. Given a set $S\subset [2n]$, we define the greedy path collection originating at $S$ in $G'_{v,w}$ to be given by the greedy paths originating at each $i\in S$. In general, this is not non-intersecting. 
    
    In the proof of \Cref{lem:bijectionsetstoexpressions}, there is an action of $s_i$ on sets $I_w$. We give a graph theoretic interpretation of this action. For each $i\in[n]$, observe that the greedy path originating from $S\subset[2n]$ in $G'_{{\rm id}, s_i}$ has sink set $s_i(S)$. Thus, the greedy path collection originating from $S$ in $G'_{{\rm id}, w}$ has sink set $w^{-1}(S)$. In particular, the sink set of the greedy path collection $P$ originating from $[2n]\setminus[n]$ in $G'_{{\rm id}, w}$ is $I'_w$. It follows from the proof of \Cref{lem:bijectionsetstoexpressions} that if $\underline{w}=s_{i_1}\cdots s_{i_l}$ we have:
    \begin{align*}
    I'_{s_{i_1}\cdots s_{i_k}} &= \Big(I'_{s_{i_1}\cdots s_{i_{k-1}}}\setminus \{i_k+1,n+i_k+1\} \Big)\cup\{i_k,n+i_k\}, \quad     \text{if }1\leq i_k\leq n-1, \\
    I'_{s_{i_1}\cdots s_{i_{k}}} &= \Big( I'_{s_{i_1}\cdots s_{i_{k-1}}}\setminus \{2n-1,2n\}\Big)\cup\{n-1,n\}, \quad \text{if } i_k=n.
    \end{align*}
    In particular for $i\in [n]$, $I'_{s_{i_1}\cdots s_{i_{k}}}\prec I'_{s_{i_1}\cdots s_{i_{k-1}}}$ and  acts by replacing two elements of $I'_{s_{i_1}\cdots s_{i_{k-1}}}$. This reflects the movement of paths in $P$ downwards along the two arrows corresponding to $s_i$ in $G'_{{\rm id},w}$. Thus paths in $P$ travel downwards along all edges of $G'_{{\rm id}, w}$.

    Observe that the greedy path originating from $[2n]$ in $G'_{{\rm id}, w}$ uses each vertical edge twice, once going upwards and once going downwards. Thus, the greedy path collection $Q$ originating from $[n]$ in $G'_{{\rm id}, w}$ uses all the vertical edges in the upwards direction.
    
    \smallskip
    
    We claim that the paths in $Q$ do not cross. Suppose $p_1$ and $p_2$ cross. Immediately before they first meet, one of them, say $p_1$, must travel along a vertical edge $e$ (necessarily upwards) and the other, $p_2$, along a horizontal strand. But then, by greediness, $p_2$ would be forced to use the vertical edge $e$ going downwards, which is a contradiction.  
    
    Since paths in $Q$ only use vertical edges in the upwards direction and do not cross one another, we may view it as a non-intersecting path collection in $G_{{\rm id}, w}$ instead of in $G'_{{\rm id}, w}$. By the definition of a greedy path collection in $G'_{{\rm id}, w}$, each individual path in $Q$ is left greedy in $G_{{\rm id}, w}$, proving (\ref{property3}). Since $Q$ uses all edges of $G'_{{\rm id},w}$, (\ref{property4}) is satisfied as well. Each path is directed from $i$ to $w^{-1}(i)$ for some $i\in [n]$. Since $w\in W^{[n-1]}$, by \Cref{rem:wAsSet}, $w^{-1}(1)\prec\cdots \prec w^{-1}(n)$, proving (\ref{property2}). In this base case, (\ref{property1}) does not assert anything. Finally, the next lemma implies it is a unique path collection. 
        
    \begin{lemma}\label{lem:usesalledges}
        Suppose $G_{v,w}$ has all vertical edges pointing upwards. Let $P$ be a non-intersecting path collection in $G_{v,w}$ using every vertical edge. Then $P$ is a unique path~collection.
    \end{lemma}
    \begin{proof}
        For a vertical edge $e$ that skips over $k-1$ strands, define $h_e \coloneqq k$. For a path collection $Q$, define $h(Q)\coloneqq \sum_e h_e$ where the sum is over all vertical edges used by $Q$. Observe that 
        \[
            h(P) = \max \Big\{ h(Q) \colon Q\text{ is a non-intersecting path collection in }G_{v,w} \Big\}.
        \]
        
        We also claim that a non-intersecting path collection is uniquely determined by the vertical arrows it uses, which implies that $P$ is the unique non-intersecting path collection maximizing $h$. To prove the claim, let $V$ be the set of vertical edges of a path collection $Q$ in $G_{v,w}$. The path collection $Q$ is obtained from $V$ as follows. Let $\nu_1$ be the rightmost arrow in $V$ (if there are multiple, choose the top one). The path $p_1$ in $Q$ which uses $\nu_1$ must terminate vertically from the head of $\nu_1$ and must approach $\nu_1$ along the strand $\sigma_{\nu_1}$ where $\nu_1$ originates. If there are no other arrows in $V$ whose head lies on $\sigma_{\nu_1}$, $p_1$ must originate on strand $\sigma_{\nu_1}$. Otherwise, let $\nu_2$ be the rightmost arrow in $V$ terminating on strand $\sigma_{\nu_1}$. By the non-intersecting condition, $p_1$ must use $\nu_2$. Continuing in this way, we see that the path $p_1$ which uses $\nu_1$ is uniquely determined. Having determined $p_1$, remove the vertical arrows used by $p_1$ from the set $V$ and run the same argument again to determine another path $p_2$. Repeating this until $V$ is empty, we obtain a path collection. Observe that since we assume that $V$ is the set of arrows of a path collection, the algorithm will in fact terminate. Since every step in the construction was forced, the resulting path collection is unique. 
        
        For $q$ a path in $Q$ whose sink lies $k'$ strands above its source, define $h_q\coloneqq k'$. Observe that $h(Q)=\sum _q h_q$ where the sum is over all paths of $Q$. Thus, $h(Q)$ is fully determined by the source and sink set of $Q$. Since $P$ uniquely maximizes $h$, there is a unique path collection with the same source and sink set as $P$, that is, $P$ is a unique path collection.
    \end{proof}
        
    This completes the proof for $G_{{\rm id}, w}$. We now prove the lemma for $G_{v,w}$ with $\ell(v)>0$ by induction. Suppose the result holds for all $G_{u,w}$ with $\ell(u)<\ell(v)$. Let $u$ be the product of the first $\ell(v)-1$ factors of the distinguished subexpression $\underline{v}^+ \prec \underline{w}$. The corresponding expression $\underline{u}$ is again distinguished in $\underline{w}$. Suppose $\underline{v}^+=\underline{u}s_i$. By \Cref{def:boundaryLGVdiagram} we obtain $G_{v,w}$ from $G_{u,w}$ by cutting the graph to the right of some pair of arrows $(\alpha_1,\alpha_2)$ in $G_{u,w}$, deleting those arrows, and performing a signed $s_i$ action on the left part. As in \Cref{def:signedaction}, we will always use \textit{strand $k$} to refer to the strand that is $k$ from the bottom. We claim $(\alpha_1,\alpha_2)$ is the leftmost pair of vertical arrow between strands $i$ and $i+1$ and strands $n+i+1$ and $n+i$ if $i \in [n-1]$ or strands $n-1$ and $2n$ and strands $n$ and $2n-1$ if $i=n$. We first show the following.

\begin{lemma}\label{lem:removingarrows}
    Let $\underline{u}$ be a reduced distinguished subexpression of $\underline{w}$. Suppose in $G_{u,w}$ there are arrows $\alpha_1,\ldots, \alpha_p$ from strand $i$ to strand $i+1$. Recall that every $\alpha_j$ corresponds to some simple reflection $s_{i_j}$ of $\underline{w}$ which is replaced by $1$ in $\underline{u}$. For $j\in [p]$, let $\underline{u}_j$ be the subexpression of $\underline{w}$ obtained by using $s_{i_j}$ instead of the corresponding $1$ in $\underline{u}$. Then, $u_1=\cdots = u_p$. An analogous result holds for arrows from strand $n-1$ to strand $n+1$. 
\end{lemma}

\begin{proof}
    Since, in \Cref{def:boundaryLGVdiagram}, we swap the source vertices in $G_{u,w}$ for each simple transposition in some expression $\underline{u}$, we can read off $u$ from the source vertices of $G_{u,w}$. From bottom to top, they are $u^{-1}(1), u^{-1}(2), \cdots,u^{-1}(n), u^{-1}(2n), u^{-1}(2n-1),\cdots, u^{-1}(n+1).$ The effect of adding a transposition from $\underline{w}$ to $\underline{u}$ is to permute in $G_{u,w}$ the source labels of the strands connected by the corresponding arrows.  
\end{proof}

By \Cref{lem:removingarrows}, the reduced distinguished expression for $v$ will be the one obtained using the leftmost possible transposition, and so $(\alpha_1,\alpha_2)$ is the leftmost pair of arrows between the appropriate strands. Since, in step \ref{LGVstep2} of \Cref{def:boundaryLGVdiagram} we only modify vertical arrows to the left of the arrows that we remove, this implies the following.

\begin{cor}\label{cor:arrowspointup}
    All arrows in $G_{v,w}$ are pointing upwards. 
\end{cor}

Suppose $(\alpha_1,\alpha_2)=a_i^{(j)}, a_{n+i+1}^{(j)}$ for some $i\in [n-1]$ (the argument is identical for $(\alpha_1,\alpha_2)=b_{n-1}^{(j)},b_{n}^{(j)}$). We compare the left greedy path collections $Q_u$ and $Q_v$ originating from $[n]$ in three different sections of the graphs $G_{u,w}$ and $G_{v,w}$, respectively.

\begin{enumerate}[label=(\roman*)]
    \item \label{collection1} Truncate the graphs $G_{v,w}$ and $G_{u,w}$ immediately to the left of  $a_i^{(j)}$. There is a bijection between path collections originating from $[n]$ in the truncated $G_{u,w}$ terminating at $J$ and those in the truncated $G_{v,w}$ terminating at $J'$, where $J'$ is obtained from $J$ by swapping $i \leftrightarrow i+1$ and $n+i\leftrightarrow n+i+1$. Moreover this bijection preserves left greediness: the truncation of $Q_u$ maps to the truncation of $Q_v$.

    \item \label{collection2} Near $a_{i}^{(j)}$ in $G_{u,w}$, a path in $Q_u$ enters on strand $i$ but not on $i+1$, uses the vertical arrow, and leaves on strand $i+1$ (any other configuration of paths entering on strands $i$ and $i+1$ is not possible since all vertical edges must be used in the left greedy path collection in $G_{u,w}$, by induction). Let $J$ and $J'$ be as in the previous bullet for the truncations of $Q_u$ and $Q_v$, respectively. Then, $i\in J\setminus J'$ and $i+1\in J'\setminus J$. Thus, $Q_v$ just uses strand $i+1$ and we preserve the left greediness of each path around $a_{i}^{(j)}$. Similar considerations apply locally near $a_{n+i+1}^{(j)}$. It follows that if we truncate $G_{u,w}$ and $G_{v,w}$ just to the right of $a_i^{(j)}$, the corresponding truncations of $Q_u$ and $Q_v$ have the same sink set.

    \item \label{collection3} Since $G_{v,w}$ is identical to $G_{u,w}$ to the right of $a_i^{(j)}$, there is a trivial bijection preserving the left greediness of paths between path collections in $G_{u,w}$ and those in $G_{v,w}$ in this part of the graphs.
\end{enumerate}

Taken all together, we can see that the sink sets of $Q_u$ and $Q_v$ are the same and $Q_v$ is non-intersecting and uses all vertical edges of $G_{v,w}$. This proves (\ref{property1}) and (\ref{property4}). The paths in this path collection are all left greedy, proving (\ref{property3}). We next show that for $k\in[n-1]$ the path originating from source vertex $k$ terminates below the path originating from source vertex $k+1$. The only place where this might break is if the left greedy paths paths in the truncation of $G_{u,w}$ in \ref{collection1} originating at $k$ and $k+1$ terminate in strands $i$ and $i+1$; this would imply that both $i,i+1$ belong to $J$, as defined in item \ref{collection1} above. However, this never happens for the left greedy path collection, as explained in item \ref{collection2}. This proves (\ref{property2}). Finally, uniqueness is guaranteed by \Cref{lem:usesalledges}. We have now completed the proof of \Cref{prop:pathcollectionproperties}.
\end{proof}

\begin{lemma}\label{lem:greedyisextremeuv}
    Let $p_i$ be the left greedy path originating from source vertex $i \in [n]$ in $G_{v,w}$. There is no path in $G_{v,w}$ originating from $i$ and achieving a $\prec$-larger sink than $p_i$.

\end{lemma}

\begin{proof}
    We start by proving the result for $G_{{\rm id}, w}$. From \Cref{lem:bijectionsetstoexpressions}, observe that in $G_{{\rm id}, w}$, whenever an arrow $a_{\lambda}^{(\kappa)}$ terminates on strand $s$ and an arrow $b_l^{(k)}$ skips a strand $s$ to the right of $a_{\lambda}^{(\kappa)}$, there is an arrow $\gamma$ originating on strand $s$ between $a_{\lambda}^{(\kappa)}$ and $b_{l}^{(k)}$. Suppose $p'$ were a path originating at $i$ and terminating above the sink of $p$. Then, at some point, $p'$ lies above $p$. Thus, $p'$ at some point uses an arrow which skips a strand $s$ used by $p$. Say this arrow is $b_l^{(k)}$ for some $l,k$. Since $b_l^{(k)}$ originates below strand $s$, source vertex $i$ does not lie on strand $s$. Say $p$ used an arrow $a_{\lambda}^{(\kappa)}$ to reach strand $s$. By our observation at the beginning of the proof, there is an arrow $\gamma$ originating on strand $s$ between $a_{\lambda}^{(\kappa)}$ and $b_l^{(k)}$. Note that the arrows $b_{l}^{(k)}$ in $G_{{\rm id}, w}$ only skip one strand. Since $p'$ gets above $p$ by skipping a single strand, $p$ does not use $\gamma$, which contradicts left greediness.

    For $G_{v,w}$, we argue by induction. Let $\underline{w}^+\prec \underline{w}_)^{[n-1]}$ be reduced distinguished let $\underline{v}^+\prec \underline{w}^+$ be reduced distinguished. Let $u< v$ be such that $\underline{v}^+=\underline{u}s_k$ for some simple transposition $s_k$. Observe that $\underline{u}\prec \underline{w}^+$ is reduced distinguished as well. Let $m_i$ be the sink of the left greedy path originating at $i$ in $G_{v,w}$ and suppose there is a path $p$ originating at $i$ and terminating at $ j\succ m_i$, that is, on a strand above $m_i$. We make two observations. 
    
    \begin{enumerate}
        \item Combining (\ref{property1}) and (\ref{property2}) of \Cref{prop:pathcollectionproperties}, it follows that the left greedy path starting at $i$ in both $G_{u,w}$ and $G_{v,w}$ have the same sink vertex $m_i$. 
        \item The existence of $p$ in $G_{v,w}$ implies that there is a path in $G_{u,w}$ from $i$ to $j'$ for some $j'\succeq j$. To see this, we work case by case, depending on how many of the strands that get permuted when forming $G_{v,w}$ from $G_{u,w}$ are used by $p$. In each case, the observation follows from splitting $G_{u,w}$ and $G_{v,w}$ into three sections each and comparing them section by section, similarly to \ref{collection1} - \ref{collection3} of the proof of \Cref{prop:pathcollectionproperties}. Note that, compared to \ref{collection1}-\ref{collection3}, we are going ``backwards", from $G_{v,w}$ to $G_{u,w}$ where $\ell(v)>\ell(u)$.
    \end{enumerate}  

    Thus, we have a path originating from $i$ in $G_{u,w}$ which terminates at $j'\succ m_i$. Since $m_i$ is the sink of the left greedy path originating at $i$ in $G_{u,w}$, this contradicts our induction hypothesis.
\end{proof}

\begin{remark}
    \Cref{prop:pathcollectionproperties} and \Cref{lem:greedyisextremeuv} with $v={\rm id}$ and $w=w_0^{[n-1]}$ imply \Cref{lem:LGVdiagramLeftGreedy}.
\end{remark}

\begin{theorem} \label{thm:RichardsonRestated}
Let $X\in \OGr^{\geq 0}(n, 2n)$ and $v, w \in W$ such that $X \in \mathcal{R}^{>0}_{v, w}$. Recall 
\[
    t_j :=  \min_{\prec}\left\{t \colon \left(I_X^{(j-1)}\setminus \{i^{(j-1)}_{j}\}\right)\cup\{t\} \text{ is a basis of } \mathscr{M}_X\right\},
\]
where $I_X^{(l)} := \{ i_1^{(l)} \prec \dots \prec i_{n}^{(l)} \}$ is the $l$-th lowering of $I_X=\{i_1\prec i_2\prec\cdots\prec i_n\}$ with respect to $\prec$ as defined in \Cref{def:lowering}. Then
\begin{enumerate}[wide=20pt]
    \item $w \in W$ is given by $w^{-1}(j)=i_j$ for $j\in[n]$, and
    \item $v \in W$ is given by $v(j)=t_j$ for $j\in[n]$.

\end{enumerate}
    In particular, the Richardson containing $X$ is fully determined from $\mathscr{M}_X$
\end{theorem}

\begin{proof}
    By \Cref{prop:pathcollectionproperties}, it suffices to show the first statement for $\mathcal{R}_{{\rm id}, w}^{>0}$. We view the arrows of $G_{{\rm id}, w}$ as representing the transpositions in an expression for $w$. The fact that each vertical edge is used by the left greedy path collection $P$ originating from $[n]$ in $G_{{\rm id}, w}$ means that the sink set of $P$ is obtained by acting on $[n]$ with each transposition in $\underline{w}$ from left to right. This yields precisely $w^{-1}([n])$.
    
    For the second statement, we consider the left greedy path collection $P$ originating from $[n]$ in $G_{v,w}$. We claim that $I^{(l)}_X$ is the sink set of the path collection $P^{(l)}$ originating from $[n]$ obtained from $P$ by replacing the paths originating from $[l]$ with horizontal paths. We show this by induction. Suppose the sink set of $P^{(l-1)}$ is $I^{(l-1)}_X$. Recall from \Cref{prop:pathcollectionproperties} that for $i,j\in [n]$ with $i<j$, the path originating from $i$ in $P$ terminates at a sink smaller than the sink of the path originating at $j$, in $\prec$ order. Thus, $i^{(l-1)}_l$ is the sink of the path originating at $l$ in $P^{(l-1)}$. To obtain $I^{(l)}_X$, we replace this element by the $\prec$ smallest possible thing. Suppose we wish to construct a path collection with sink set $I^{(l)}$. Since each individual path in $P$ is left greedy by \Cref{prop:pathcollectionproperties}, and reaches the highest possible sink by \Cref{lem:greedyisextremeuv}, the subcollection of $P$ from $\{l+1,\ldots, n\}$ to $\{i^{(0)}_{l+1},\ldots, i^{(0)}_n\}=\{i^{(l)}_{l+1},\ldots, i^{(l)}_n\}$ is the unique path collection originating from a subset of $[n]$ with that sink set. Thus, in any path collection from $[n]$ to $I^{(l)}$, these paths appear. To obtain a path collection with sink set $I^{(l)}$, we have to lower the path originating at $l$, that is, modify it to have a lower sink. The lowest possible sink that can be obtained by a collection of paths originating from $[l]$ is a collection of horizontal paths, since there are no vertical edges going downwards in $G_{v,w}$ by \Cref{cor:arrowspointup}. Thus, the sink set of $P^{(l)}$ is indeed $I^{(l)}_X$, as claimed.

    It follows that $t_j$ is the sink of the horizontal path originating from source vertex $j$ in $G_{v,w}$. Since the labels of the source vertices in $G_{v,w}$ are, from bottom to top,
    \[
        v^{-1}(1),\ldots, v^{-1}(n), v^{-1}(2n),\ldots, v^{-1}(n+1),
    \]
    this implies that $t_j=v(j)$.
\end{proof}

\begin{example}
    We use the same $1$-parameter family as in \Cref{ex:NonNegn4} to apply \Cref{thm:Nonnegative} and check that the following matrix is in $\SS_4^{\geq 0}$:
    \[
    A=\begin{bmatrix}
       0&0&0&2\\
       0&0&0&0\\
       0&0&0&2\\
       -2&0&-2&0
       \end{bmatrix}.
    \]
    Using our {\tt Macaulay2} implementation of \Cref{thm:RichardsonRestated} we computed that $\begin{bmatrix} \Id_n | A \end{bmatrix} \in \mathcal{R}^{>0}_{2134, 2385}$ where the permutations are written in one-line notation.
\end{example}

\begin{proposition}\label{prop:BoundSSMat} 
    The following holds
    \[
        \SS_n^{ \geq 0} = \bigsqcup_{v\in W_{[n-1]}} \mathcal{R}_{v,w}^{>0}
    \]
\end{proposition}

\begin{proof}
     The set of positive skew-symmetric matrices $\SS_n^{>0}$ is in bijection with the set of points in $\OGr^{\geq 0}(n,2n)$ represented by matrices $X$ that can be row reduced to the form $\begin{bmatrix} \Id_n | A \end{bmatrix}$. For this, we must have that the minor of the first $n$ columns of $X$ is nonzero. That is, there must be a path collection from $[n]$ to $[n]$ in the graph $G_{v,w}$. By \Cref{cor:arrowspointup}, this is possible if and only if the source vertices on the bottommost $n$ strands are $[n]$. This happens when $v^{-1}([n])=[n]$, equivalently $v([n])=[n]$, which is exactly the condition that $v \in W_{[n-1]}$.  
\end{proof}

\section{Pfaffian signs}\label{sec:5}

    We recall that the determinant $\det(A)$ of a skew-symmetric matrix $A = (a_{ij})_{i,j=1}^{2m}$ of size $2m \times 2m$ is a homogeneous polynomial of degree $2m$ in the entries of $A$. Moreover, it is the square of a degree $m$ homogeneous polynomial $\Pf(A)$ known as the \emph{Pfaffian},
    \begin{equation}
        \Pf(A) := \frac{1}{2^m m!} \sum_{\sigma \in S_{2m}} {\rm sgn}(\sigma) \prod_{i=1}^{m} a_{\sigma(2i-1),\sigma(2i)}.
    \end{equation}
    For any skew-symmetric matrix $A \in \SS_{n}$ and any set $I \subset [n]$ of even size we denote by $\Pf_I(A)$ the pfaffian of the submarix $A_{I}$ of $A$ whose rows and column indices are in $I$. By convention, we set $\Pf_{\emptyset}(A) = 1$. We denote the semi-algebraic set of \emph{Pfaffian-positive} skew-symmetric matrices by 
    \begin{equation}
        \SS_n^{\Pf > 0} \coloneqq \Big\{ A \in \SS_n :  \sgn(I, [n]) \ \Pf_I(A) > 0 \quad \text{for any } I\subset [n] \text{ of even size}\Big\},
    \end{equation}
    and denote by $\SS^{\Pf \geq 0}_n$ its Euclidean closure in $\RR^{n \times n}$.
    
    \smallskip
    
    \begin{remark}
        Let $g = {\rm diag}(1, -1, 1 \dots, (-1)^{n+1})$.  We note that a skew-symmetric matrix $A$ is in $\SS_n^{\Pf > 0}$ if and only if $\Pf_I(A') > 0$ for any subset $I \subset [n]$ of even size where $A'  = g A g^T$. So the set $\SS_n^{\Pf > 0}$ is conjugate to set of matrices with positive pfaffians. This justifies our~notation.
    \end{remark}
    
    A priori, there is no reason why the notion of total positivity in $\SS_n$ we have discussed so far should be compatible with $\SS_n^{\Pf > 0}$. In this section, we briefly discuss how these two notions of positivity are related, namely we prove \Cref{thm:PfaffianSign} and discuss a few consequences. To do so, we shall need some background on the half-spin representation. For more details on Lie theoretic background we refer the reader to \cite{Hall, Procesi}.

    \smallskip

    Recall the vector spaces $E$ and $F$ from \eqref{eq:EFspaces} and let $V = E \oplus F \cong \RR^{2n}$. The space $V$ is equipped with the quadratic form $q$ as in \eqref{eq:quadForm}. The exterior algebra $\bigwedge^\bullet E$ decomposes into a direct sum of two vector spaces \(\bigwedge ^\bullet E = S \oplus S^- \), where $S$ (resp. $S^-$) is the subspace of elements in $\bigwedge^\bullet E$ of degree equal to $n$ (resp. $n+1$) modulo 2. For any $I=\{i_1<\ldots< i_\ell\} \subset [n]$ we denote $e_I := e_{i_1}\wedge \cdots \wedge e_{i_\ell}$. The $2^{n-1}$ vectors $e_{[n] \setminus I}$ where $I$ is a subset of $[n]$ of even size form a basis of $S$. 
    
    The orthogonal Grassmannian $\OGr(n,2n)$ can be minimally embedded as the \emph{Spinor variety} in $\PP(S)$. Here, we briefly explain how this embedding is obtained using Pfaffians. Following \cite[Section 2]{Manivel} and \cite[Section 11.7]{Procesi}, let $\mathscr{C}$ be the \emph{Clifford algebra} of the quadratic space $(V,q)$. This algebra is defined as the quotient of the tensor algebra $T^\bullet V$ modulo the relations 
        \begin{equation}
        u \otimes v + v\otimes u = 2 q(u,v), \quad \text{for } u,v \in  V.
        \end{equation}
    Let $f := f_1 \cdots f_n \in \mathscr{C}$ and note that the left-ideal $\mathscr{C} \cdot f$ is isomorphic (as a vector space) to $\wedge^\bullet E$ via $ e_I \cdot f \mapsto e_I$. Now, let $A \in \SS_n$ and $H=\rowspan \begin{bmatrix} \Id_n | A \end{bmatrix}$. Then, $H$ is spanned by 
    \[
        e_i(A) := e_i + \sum_{j=1}^{n} a_{ij} f_j, \quad \text{for } 1 \leq i \leq n.
    \]  
     Let $h = e_1(A) \cdots e_{n}(A) \in \mathscr{C}$, then the \emph{pure spinor} of $A$, namely the point in $\PP(S)$ that represents $H$, is the one dimensional linear space $h\mathscr{C} \cap \mathscr{C} f$ in $\mathscr{C} f \cong \wedge^\bullet E$, see \cite[III.1.1]{Chevalley}. The line $h\mathscr{C} \cap \mathscr{C} f$ in $S$ is spanned over $\CC$ by $h\cdot f \in \mathscr{C} \cdot f$. In other words, the pure spinor of $A$ is the element $u_H \in \bigwedge ^\bullet E$ such that of $u_H\cdot f = h\cdot f$. 
     
    \begin{lemma}\label{lem:ManivelFormCorrect}
            Given a generic point in $\OGr(n,2n)$ of the form $\begin{bmatrix} \Id_n | A \end{bmatrix}$ with $A \in \SS_n$, the corresponding pure Spinor in $\PP(S)$ is
        \begin{equation}\label{eq:SpinorEmbedding}
            \sum_{I \subset [n] \text{ of even size}} \sgn(I,[n]) \ 2^{|I|/2} \ \Pf_{I}(A) \ e_{[n] \setminus I} \quad \in S.
        \end{equation}
    \end{lemma}
    \begin{proof}
         It is enough to compute $h \cdot f$ in $\mathscr{C}$. We proceed by induction on $n$. Suppose the result holds for all $A \in \SS_n$. Fix $A=(a_{i,j})_{i,j=1}^{n+1} \in \SS_{n+1}$ and let us write 
        \[
            u_i := e_{i} + \sum_{j=1}^{n} a_{ij} f_j \quad
            \text{for } 1 \leq i \leq n+1.
        \]
        Note that for any $1\leq i \leq n+1$ and $1 \leq j \leq n$ we have $f_i^2 = 0$ and $f_{n+1} u_j = - u_j f_{n+1}$. We then have
        \begin{align*}
             e_1(A) \cdots e_{n}(A) \ e_{n+1}(A) \ f 
             &=  e_1(A) \cdots e_{n}(A) \ e_{n+1} \ f \\
             &=  \left(u_1 + a_{1,n+1} f_{n+1}\right)\cdots \left(u_n + a_{n,n+1} f_{n+1}\right) \ e_{n+1} \ f  \\
             & = \left(u_1 \cdots u_n  + \sum_{i=1}^{n} (-1)^{n-i} a_{i,n+1} u_1 \cdots \widehat{u_i} \cdots u_n \ f_{n+1} \right) \ e_{n+1} \ f\\
             & = u_1 \cdots u_n \ e_{n+1} \  f + \sum_{i=1}^{n} (-1)^{n-i} a_{i,n+1} u_1 \cdots \widehat{u_i} \cdots u_n \ f_{n+1} \ e_{n+1} \ f.
        \end{align*}
        By the induction hypothesis, we have 
        \[
        u_1 \cdots u_n e_{n+1} f = \sum \limits_{K \subset [n] \text{ even}} {\rm sgn} (K,[n]) \ 2^{|K|/2} \ \Pf_K(A) e_{[n+1] \setminus K} \ f.
        \]
        Note that $f_{n+1} e_{n+1} \ f = ( - e_{n+1}f_{n+1} \ + \  2 ) \ f = 2 f$. For each $i \in [n]$ we denote $A'_i$ for the submatrix of $A$ with rows and columns $i,n+1$ removed. Then,
        \[
        u_1 \cdots \widehat{u_i} \cdots u_n f = \sum_{K \subset [n]\setminus i \text{ even}} {\rm sgn} (K, [n+1]) \ 2^{|K|/2} \ \Pf_K(A'_i) e_{[n] \setminus (K \cup i)}f.
        \]
        Using the recursive formula for Pfaffians, we obtain
        \begin{align*}
           &\sum_{i=1}^{n} (-1)^{n-i} a_{i,n+1} \ u_1 \cdots \widehat{u_i} \cdots u_n  \ f_{n+1} \ e_{n+1} \ f \\
            &=\sum_{i=1}^{n} (-1)^{n-i} \ 2 \ a_{i,n+1} \ u_1 \cdots \widehat{u_i} \cdots u_n f\\
            &= \sum_{i=1}^{n} \sum_{K \subset [n]\setminus i \text{ even}}(-1)^{n-i} \ 2^{1 + |K|/2} \ a_{i,n+1}{\rm sgn} (K,[n]\setminus i) \ \Pf_K(A'_i) \ e_{[n] \setminus (K \cup i)} \ f\\
            &= \sum_{\substack{K = \{i_1<\ldots< i_l\} \text{ even}\\ i_l=n+1 }} {\rm sgn} (K,[n+1]) 2^{|K|/2} \left( \sum_{j=1}^{l-1} (-1)^{j}a_{i_j,n+1} \ \Pf_{K\setminus \{i_j,n+1\}}(A) \right) \ e_{[n+1] \setminus K} \ f\\
            &=\sum_{\substack{K \subset [n+1] \text{ even} \\ n+1 \in K}} {\rm sgn} (K,[n+1]) \ 2^{|K|/2} \ \Pf_{K}(A) \ e_{[n+1] \setminus K} \ f. \qedhere
        \end{align*} 
    \end{proof}
    We note that \Cref{eq:SpinorEmbedding} appears in \cite[Section 2.3]{Manivel} up to a sign and constant inconsistency, which is not relevant for their purposes.

    \smallskip
    
    The vector space $S$ can be endowed with an action of the spin group $\Spin(2n)$. This is the half-spin representation. It is an irreducible representation corresponding to the highest weight vector $\lambda = \frac{1}{2}(1,\dots, 1) \in \RR^{n}$. The basis $\{e_{[n]\setminus I}: I \subset [n] \text{ of even size}\}$ is a basis of weight vectors of $S$ for this representation, and the weights of $S$ are $ \frac{1}{2} \left(\pm 1,\ldots, \pm 1\right)$ where the number of minus signs is even. Explicitly, the weight of $e_{[n]\setminus I}$ has $-$ in the positions indexed by $I$, and in particular the highest weight vector is $e_{[n]}$ with weight $\lambda$.

    \smallskip
    
    We now describe $S$ as a representation of the Lie algebra $\mathfrak{spin}(2n)\cong \mathfrak{so}(2n)=\bigwedge^2 V$. Following \cite[Chapter II]{Chevalley} and \cite[Section 5]{SpinGeometry} the action of $\mathfrak{so}(2n)$ on $S$ is given as follows. For $v=v' \oplus v'' \in E \oplus F $ and $\omega=\omega_1\wedge \cdots \wedge \omega_k \in \bigwedge^{k} E$, we write
    \[
    v\cdot \omega := v'\wedge \omega_1\wedge \cdots \wedge \omega_k+2\sum_{i=1}^k (-1)^{i-1} q(v'',\omega_i) \ \omega_1\wedge \cdots \omega_{i-1}\wedge\omega_{i+1}\wedge\cdots \wedge \omega_k.
    \]
    
    Note that $v \cdot \omega \in \bigwedge^{k-1} E \oplus \bigwedge^{k+1} E$. Then for $\omega \in S$ we have $v \cdot \omega \in S^{-}$. Acting again with another element $w \in V$ yields $w\cdot(v \cdot \omega) \in S$. Hence we obtain an action of $\mathfrak{so}(2n)$ on $S$:
    \[
        (v_1\wedge v_2)\cdot \omega = \frac{1}{4} [v_1,v_2] \cdot  \omega = \frac{1}{4}(v_1 \cdot (v_2 \cdot \omega) - v_2 \cdot (v_1 \cdot \omega) ) \quad \text{for any } v\wedge w \in \mathfrak{so}(2n) \text{ and } \omega \in S.
    \]
    The factor of $1/4$ comes from the degree $2$ covering map $\Spin(2n) \to \SO(2n)$.

    \smallskip
    
    The $2^{n-1}$ Pfaffians of $A$ are related to certain regular functions on $\OGr(n,2n)$ known as \emph{generalized minors}. These functions are positive on $\OGr^{>0}(n,2n)$, see \cite[Lemma 11.4]{MR}. To define them we need to fix the root space decomposition of $\mathfrak{so}(2n)$. The maximal torus \eqref{eq:maxTorus} of $\SO(2n)$ corresponds to a Cartan subalgebra $\mathfrak{h}$ of $\mathfrak{so}(2n)$. We identify the dual $\mathfrak{h}^\ast$ of $\mathfrak{h}$ with the vector space $\RR^n = \Span_{\RR}(\varepsilon_1,\ldots, \varepsilon_n)$. The action of the Weyl group $W$ on $\RR^n$ is given through its generators \eqref{eq:WeylGroupGens} as follows
    \begin{equation}\label{eq:WeylGroupActionOnWeights}
        s_i\cdot (a_1,\ldots,a_n)=\begin{dcases}
        (a_1,\ldots,a_{i-1},a_{i+1},a_i,a_{i+2},\ldots,a_n) & \text{if } 1\leq i \leq n-1,\\
        (a_1,\ldots, a_{n-2},-a_n,-a_{n-1}) &\text{if } i = n.
    \end{dcases}
    \end{equation}
    We choose a set of simple roots $\Phi = \{\varepsilon_1-\varepsilon_2, \varepsilon_2-\varepsilon_3, \ldots, \varepsilon_{n-1}-\varepsilon_{n}, \varepsilon_{n-1}+\varepsilon_n\}$. The positive roots are then
    \[
    \Phi^{+} = \{\varepsilon_i-\varepsilon_j : 1\leq i < j \leq n\} \cup \{\varepsilon_i+\varepsilon_j: 1\leq i,j \leq n\}.
    \]
    
    Finally, we note that the following roots and their root vectors correspond to the pinning in \Cref{sec:2} by taking the matrix exponential $ \exp: \mathfrak{so}(2n) \cong \wedge^2 \CC^{2n} \to \SO(2n)$: 
    \begin{table}[H]
        \centering
        \begin{tabular}{|c|c|c|c|}
             \hline
            Roots        & $\varepsilon_i-\varepsilon_j$ for  $i \neq j$  & $\varepsilon_i+\varepsilon_j$ for $i \leq j$ & $ -\varepsilon_i-\varepsilon_j$ for $i \leq j$ \\
             \hline
            Root vectors & $e_i\wedge f_j$ & $e_i \wedge e_j$ & $f_j \wedge f_i$\\
             \hline
        \end{tabular}
        \caption{The roots of $\mathfrak{so}(2n)$ and their corresponding root vectors in the half spin representation $S$.}
        \label{tab:rootsAndrootVectors}
    \end{table}
    
    Recall the definition of $\dot{s}_i$ from \eqref{eq:SiDot}.
     \begin{lemma}\label{lem:logSi}
        Let $I$ be a subset of $[n]$ of even size. Then for $1 \leq i \leq n-1$,
        \[
            {\rm Log}(\dot{s}_i) \cdot e_{[n]\setminus I} = \begin{cases}
                 - \frac{\pi}{2} e_{[n]\setminus ( (I \setminus i) \cup i+1 )}   \quad & \text{if } i \in I \text{ and } i+1 \not \in I, \\
                \frac{\pi}{2} e_{[n]\setminus ( (I \setminus i+1) \cup i )}  \quad & \text{if } i \not \in I \text{ and } i+1 \in I, \\
                0  \quad & \text{otherwise,} 
            \end{cases}
        \]
        and
        \[
            {\rm Log}(\dot{s}_n) \cdot e_{[n]\setminus I} =
            \begin{cases}    
                \frac{-\pi}{4}e_{[n]\setminus (I \setminus \{n-1, n\})}  \quad & \text{if } n \in I \text{ and } n+1  \in I, \\
                \pi e_{[n]\setminus (I \cup \{n-1, n\})} \quad & \text{if } n \not \in I \text{ and } n+1 \not \in I, \\
                0  \quad & \text{otherwise.} 
            \end{cases}
        \]
        Here, ${\rm Log} \colon {\rm Spin}(2n) \to \mathfrak{spin}(2n)$ is the logarithm map between a simply connected Lie group and its Lie algebra.
    \end{lemma}
    \begin{proof}
    Using a matrix logarithm for the root subgroup corresponding to $s_i$ in ${\rm Spin}(2n)$ yields the following formula for ${\rm Log}(\dot{s}_i) \in \mathfrak{spin}(2n) \cong \wedge^2 (E \oplus F)$:
    \[
    {\rm Log}(\dot{s}_i) = \begin{dcases}
        \frac{\pi}{2}(e_{i+1}\wedge f_i-e_i\wedge f_{i+1}) & \text{if } 1\leq i \leq n-1,\\
        \frac{\pi}{2}(f_{n}\wedge f_{n-1}-e_{n-1}\wedge e_{n}) &\text{if } i = n. 
    \end{dcases}
    \]
    Note that, for $I \subset [n]$ of even size, we have for any $i,j \in [n]$
    \begin{align*}
    & e_i\cdot e_{[n]\setminus I} =  \begin{dcases}
            (-1)^{\sgn(i,I\setminus i)} e_{[n]\setminus (I\setminus i)}, &\text{ if } i \in I,\\
            0, &\text{otherwise},
        \end{dcases} \\
     & f_j\cdot e_{[n]\setminus I} = 2 \begin{dcases}
            (-1)^{\sgn(j,I)} e_{[n]\setminus (I\cup j)}, &\text{ if }  j\notin I,\\
            0, &\text{otherwise.}
        \end{dcases}    
    \end{align*}
    The following computation gives the desired result for ${\rm Log}(\dot{s}_i) \cdot e_{[n]\setminus I}$. A similar computation yields the result for $\dot{s}_n$.
    \begin{align*}
        &(e_{i+1}\wedge f_i-e_i\wedge f_{i+1})\cdot e_{[n]\setminus I} = \frac{1}{4}\left(e_{i+1}f_i-f_ie_{i+1}-e_if_{i+1}+f_{i+1}e_{i}\right)\cdot e_{[n]\setminus I}\\
        &= \frac{1}{2}\left(e_{i+1}f_i-e_if_{i+1}\right)\cdot e_{[n]\setminus I}=\begin{dcases}
            -e_{[n]\setminus((I\setminus i)\cup i+1)} & \text{if } i \in I \text{ and } i+1 \not \in I, \\
            e_{[n]\setminus ( (I \setminus i+1) \cup i )}  \quad & \text{if } i \not \in I \text{ and } i+1 \in I, \\
            0  \quad & \text{otherwise.} 
        \end{dcases}
    \end{align*}
    \end{proof}

    We are now ready to prove \Cref{thm:PfaffianSign}. 
    
    \begin{proof}[\bf Proof of \Cref{thm:PfaffianSign}]
    Given a skew-symmetric matrix $A \in \SS_n$, the Pfaffians of $A$ are, up to fixed scalars, an instance of \emph{generalized minors} of the point $[\Id_n | A]$ in $\OGr(n,2n)$ for the half-spin representation $S$. Following \cite[Section III, 1.3-1.7]{Chevalley}, we know that if $g \in \Spin(2n)$ is a lift of an element $\begin{bmatrix}
        \Id_n & \ast\\ A & \ast
    \end{bmatrix}$ of $\SO(2n)$ where $A \in \SS_n$, then $g\cdot e_{[n]}$ is the pure spinor corresponding to the rowspan of $\begin{bmatrix}
        \Id_n | A
    \end{bmatrix}$. Thus, by \Cref{lem:ManivelFormCorrect},
    \begin{equation}\label{eq:PfaffianExpansions}
        g \cdot e_{[n]} = \sum_{\substack{I \subset [n]\\ |I| \text{ even}}} \sgn(I, [n]) \  2^{|I|/2} \ \Pf_I(A) \ e_{[n]\setminus I}.
    \end{equation}
    Note that the change from row to columns convention is due to \Cref{rem:MRParamConvention}. The generalized minors corresponding to the half-spin representation $S$ are given by
    \[
        m_w(g)=\langle g \cdot e_{[n]}, \dot{w}\cdot e_{[n]} \rangle,
    \]
    where $\dot{w} \in \Spin(2n)$ is a lift of $w \in W$ to the Spin group, and $\langle  \cdot, \cdot \rangle$ is the standard inner product with respect to the basis $(e_{[n] \setminus I})$ of $S$. We stress that this lifting depends on the choice of a Cartan subalgebra $\mathfrak{h}$, the simple roots, and root vectors made previously for $\mathfrak{so}(2n)$, that is, a pinning of $\Spin(2n)$ compatible with that of $\SO(2n)$ in \Cref{sec:2}.
    
    \smallskip
    
    Recall from \Cref{def:Iw} that $ I_w := \{i \in [n]: w(i)>n \}$.
    
    \begin{lemma}
        For any $w \in W$ we have 
        \[
            \dot{w} \cdot e_{[n]} = c_w \ e_{[n] \setminus I_{w^{-1}}}, \quad \text{for some nonzero scalar } c_w \in \RR.
        \]
    \end{lemma}
    \begin{proof}
        A classical result from the representation theory of Lie algebras \cite[Section 7.2 and Lemma 21.3]{HumphreysBook} guarantees that for any $w \in W$, the vector $\dot{w} \cdot e_{[n]}$ is a weight vector of weight $w \cdot \lambda$, where $\lambda$ is the highest weight. Hence, it is enough to show that $w \cdot \lambda$ has $-\frac{1}{2}$ in the positions indexed by the set $I_{w^{-1}}$.  By equation \eqref{eq:WeylGroupActionOnWeights}, $w\in W$ acts on $\lambda$ by applying a signed permutation. Observe that when we apply $w$ to $\lambda=(\frac{1}{2}, \cdots , \frac{1}{2},-\frac{1}{2}, \cdots , -\frac{1}{2})$, we have $w^{-1}(i)>n$ if and only if $(w\cdot \lambda)_{i}<0$. Since all the weight spaces of $S$ are one-dimensional, we obtain $\dot{w} \cdot e_{[n]} = c_w \ e_{[n]\setminus I_{w^{-1}}}$, where $c_w \in \RR$ as desired.
    \end{proof}

     Hence, from \eqref{eq:PfaffianExpansions}, we deduce that
    \[
        m_{w}(g) = \langle g \cdot e_{[n]} , \dot{w}\cdot e_{[n]} \rangle = c_w \sgn(I_{w^{-1}},[n])\ 2^{|I_{w^{-1}}|/2} \ \Pf_{I_{w^{-1}}}(A).
    \]
    Then, by \cite[Lemma 11.4]{MR} (and, indirectly, by \cite[Theorem 3.4]{Lusztig2}), it is enough to prove that the scalars $c_w$ are positive. To do so we first need the following.
   
    \smallskip
    
    From \Cref{rem:wAsSet} there is a bijection $w W_{[n-1]} \mapsto W_{[n-1]}w^{-1} \mapsto I_{w^{-1}}$ between left cosets of $W_{[n-1]}$ and even subsets of $[n]$. Hence, it is enough to consider the minimal representatives of left-cosets in $W / W_{[n-1]}$ to obtain all Pfaffians as generalized minors.

    \smallskip
    
    We now prove inductively on the length of $w \in W^{[n-1]}$ that $c_w>0$. For $\ell(w)=0$, $w={\rm id}$ and $I_{w^{-1}} = \emptyset$. Then, we have the desired sign because $c_{\rm id}=1$. Assume that $c_{w'}>0$ and let $w= s_i \ w'$, for some $i \in [n]$. Since $w'$ is a minimal left coset representative, $w$ is as well. We have $\dot{w}\cdot e_{[n]} = \dot{s}_i \ \dot{w}'\cdot e_{[n]} = c_{w'} \ \dot{s}_i \cdot e_{[n] \setminus I_{(w')^{-1}}}$. As we are considering a representation of the simply connected group of type $D_n$, $\Spin(2n)$, we have $\dot{s}_i\cdot e_{[n]\setminus I_{(w')^{-1}}}=\exp({\rm Log}(\dot{s}_i))\cdot e_{[n]\setminus I_{(w')^{-1}}}$. 
    
    \smallskip

    Since $\ell(w')=\ell(w)-1$, $(w')^{-1}(\alpha_i)$ is a positive root. If $1\leq i \leq n-1$, this means $(w')^{-1}\varepsilon_i-(w')^{-1}\varepsilon_{i+1}> 0$. However, since $w$ and $w'$ are minimal coset representatives, $I_{w^{-1}} \neq I_{(w')^{-1}}$ and so either $i \in I_{(w')^{-1}}$ or $i+1 \in I_{(w')^{-1}}$, but not both. Thus, $i+1 \in I_{(w')^{-1}}$ and $i \notin I_{(w')^{-1}}$, that is, $I_{w^{-1}}=(I_{(w')^{-1}}\setminus \{i+1\})\cup \{i\}$. If $i=n$, since $|I_{(w')^{-1}}|$ is even, $(w')^{-1}\alpha_n > 0$ implies $n-1,n \notin I_{(w')^{-1}}$. In this case, $I_{w^{-1}} = I_{(w')^{-1}} \cup \{n-1,n\}$.

    \smallskip
    
    By \Cref{lem:logSi}, $\Span_{\RR}(e_{[n]\setminus I_{(w')^{-1}}},e_{[n]\setminus I_{w^{-1}})}$ is invariant under the action of ${\rm Log}(\dot{s}_i)$. We can calculate the corresponding action of $\dot{s}_i$ via the matrix exponential, obtaining $c_w=c_{w'}$ if $i \in [n-1]$ and $c_w=2c_{w'}$ if $i=n$. 
    \end{proof}    

    \begin{example}\label{ex:NegPfaff}
        The matrix $A$ in \Cref{ex:NonNegn4} is manifestly not nonnegative as its Pfaffians do not satisfy the sign pattern in \Cref{thm:PfaffianSign}. Indeed the entry $(3,4)$ of $A$ is $-2<0$. 
    \end{example}

    \begin{example}
    We have shown that $\SS_n^{\geq 0} \subset \SS_{n}^{\Pf \geq 0}$. This inclusion is strict. For example
        \[
            A = \begin{bmatrix}
                0&0&0&2\\
                  0&0&1&0\\
                  0&-1&0&2\\
                  -2&0&-2&0
            \end{bmatrix}
        \]
    belongs in $\SS_{n}^{\Pf \geq 0}$. However, using \Cref{thm:Nonnegative} as in \Cref{ex:NonNegn4}, one can check that the minors $M_{j,k}(B(\epsilon))$ are
    \begin{align*}
        &8\,\epsilon+ o(\epsilon^{2}), \quad  &&16 \, \epsilon + o( \epsilon^{2}),& \quad    & 8\,\epsilon +  o(\epsilon^{2}),\\
        &2 + o(\epsilon^{2}), \quad    &&\mathbf{-2} +  o(\epsilon^{2}),&  \quad  & 2 + o(\epsilon).
    \end{align*}
    Thus, the matrix $A$ is \emph{not} nonnegative.
    \end{example}

\section{Future directions}\label{sec:6}
\subsection{Cluster structure}
The connection between positivity and cluster algebras has arisen frequently since the seminal works of Bernstein, Fomin and Zelevinsky, see for example \cite{BFZIII}. Consequently, whenever there is a positivity test for a semi-algebraic region of a projective variety in terms of regular functions, it is natural to ask whether they form a cluster in some cluster algebra. The cluster algebra structure of partial flag varieties was described in general in \cite{GeissLeclercSchroer} using the language of cluster categories, but it is difficult to infer whether it relates to the minors $M_{j,k}$ in the case of the orthogonal Grassmannian. We briefly comment on how the relationship between positivity and cluster structures on the orthogonal Grassmannian $\OGr(n,2n)$ compares to related varieties. 

\smallskip

Recent work of Bossinger and Li \cite{BossingerLi} explored cluster structures on partial flag varieties of type A, giving an explicit interpretation of the cluster structure in \cite{GeissLeclercSchroer}. Bossinger and Li describe how to obtain the cluster structure on two specific families of partial flag varieties from the cluster structure on a Grassmannian and conjecture this can be done more generally. In particular, starting from an initial quiver for the Grassmannian, they give a sequence of mutations, freezings of mutable vertices, and deletions of vertices that transform it into a quiver of the partial flag variety cluster algebra. 
Since positivity of the variables in a single cluster guarantees positivity of all cluster variables, this work provides a cluster theoretic grounding for the positivity test described in \cite{boretsky} and can be interpreted as conjecturally giving an explicit positivity test for any (type A) partial flag variety. In this setting, all Pl\"ucker coordinates are cluster variables and thus are positive on the positive flag variety. 

\smallskip

We now turn to $\OGr(n,2n)$. Unlike the type A setting, our positivity asks that a set of Pl\"ucker coordinates have a fixed sign pattern, so that some coordinates are forced to be negative. Moreover, there exist Pl\"ucker coordinates whose sign is not fixed on the positive orthogonal Grassmannian. In \Cref{ex:n=4MRparam} we described the open Deodehar component in $\OGr(4,8)$; the Pl\"ucker coordinate $\Delta^{1458}(X)=t_1^2t_2^2t_3t_4-t_1t_2t_3t_4t_5t_6$ has indeterminate sign for $t_i \in \RR_{>0}$. We remark there are similar behaviors in the BCFW tiles of the $(n,k,m=4)$ amplituhedron \cite{BCFWClusterStructure}. Using \textit{signed seeds}, the authors of \cite{BCFWClusterStructure} are still able to obtain the cluster structure on a standard BCFW tile from a seed of the cluster structure on $\Gr(4,n)$ by freezing certain variables, where the Pl\"ucker coordinates whose signs are not fixed on the tiles do not appear as cluster variables. This could offer a hint of how to proceed in the case of the orthogonal Grassmannian.

\subsection{Positroids}
One could try to study \textit{orthogonal positroids} of rank $n$. Positroids are rich combinatorial objects, which can be concisely described in terms of many familiar structures, including (decorated) permutations, (matroid) polytopes, planar (bi-colored) graphs, and (Le-diagrams on) tableaux. Let $\mathbb{k}$ be a field and let $X\in \mathbb{k}^{\binom{n}{k}}$ have coordinates $\Delta^I(X)$ indexed by $\binom{[n]}{k}$. Define the \textit{support} of $X$ as ${\rm supp}(X)=\{I\colon \Delta^I(X)\neq 0\}$. Then, the rank $k$ positroids on $[n]$ are $\{{\rm supp}(X)\colon X\in \Gr^{\geq 0}(k,n)\}$. These are in bijection with cells in the Richardson decomposition of $\Gr(k,n)$. Points of the Grassmannian over the sign hyperfield $\mathscr{S}$, in the sense of 
\cite{BakerBowler}, are called oriented matroids. By \cite{ArdilaRinconWilliams,SpeyerWilliams21}, positively oriented matroids, that is, those with all coordinates equal to $1$ or $0$, are in bijection with positroids; If we extend the definition of support to allow $\mathbb{k}$ to be a hyperfield, then positroids are precisely the supports of positively oriented matroids. To generalize this story, we can analogously define orthogonal positroids of rank $n$ to be $\{{\rm supp}(X)\colon X\in \OGr^{\geq 0}(n,2n)\}$. We note that $\OGr(n,2n)$ is cut out by Pl\"ucker relations and by some additional relations $\mathcal{S}$ equating certain Pl\"ucker coordinates (possibly with a sign). We propose the following definitions.

\begin{definition}
    The \textit{orthogonal Grassmannian $\OGr_{\mathscr{H}}(n,2n)$ over a hyperfield} $\mathscr{H}$ is the subset of the Grassmannian $\Gr(n,2n)$ over $\mathscr{H}$ consisting of points which satisfy $\mathcal{S}$.
\end{definition}

\begin{definition}
    We say a point of $\OGr_{\mathscr{H}}(n,2n)$ is a (rank $n$) \textit{orthogonal oriented matroid} on $2n$ elements.
\end{definition}

As in our consideration of cluster algebra structures, we need to be careful since not all Pl\"ucker coordinates have a fixed sign on the nonnegative orthogonal Grassmannian. For each $\epsilon = (\epsilon_I)\in \{-1,1\}^{\binom{2n}{n}}$ indexed by $\binom{[2n]}{n}$, let

\[
    \OGr_{\mathscr{S}}^\epsilon(n,2n)=\left\{\Delta^I\in \{-1,0,1\}^{\binom{2n}{n}}\in \OGr_{\mathscr{S}}(n,2n) \colon \epsilon_I \Delta^I\in \{0,1\} \text{ for each } I\in \binom{[2n]}{n}\right\}.
\] 

\begin{definition}
    We say that an oriented orthogonal matroid in $\OGr_{\mathscr{S}}^\epsilon(n,2n)$ is an \textit{orthogonal $\epsilon$-oriented matroid.}
\end{definition}

It is straightforward to show that there exists $\epsilon$ such that each orthogonal positroid is the support of an orthogonal $\epsilon$-oriented matroid. However, it is not clear if there are choices of $\epsilon$ for which the converse holds.

\smallskip

Finally, while there are type D Le diagrams as constructed in \cite{LamWilliamsCominiscule}, there are many questions to explore surrounding the combinatorics of orthogonal positroids. Concretely, one may ask if there exist objects generalizing decorated permutations, matroid polytopes, plabic graphs, or other familiar positroid cryptomorphisms which admit nice combinatorial descriptions and are in bijection with the Richardson cells $\mathcal{R}_{v,w^{[n-1]}}$.

\subsection{Related semi-algebraic sets}

 In the beginning of this article, we made a choice of the quadratic form $q$. However, besides connecting the orthogonal Grassmannian to skew-symmetric matrices, this choice was somewhat arbitrary and our result translates nicely to other quadratic forms of signature $(n,n)$. Explicitly, changing quadratic forms is the same as doing a change of basis, and if we also shift the pinning by the same change of bases, we obtain a shifted positivity test. We denote by $\OGr_{q}(n,2n)$ (resp. $\OGr^{\geq 0}_{q}(n,2n)$) the orthogonal Grassmannian (resp. nonnegative orthogonal Grassmannian) for the quadratic form $q$. Another natural choice of quadratic form we may consider is
 \[
     q':\RR^{2n} \xrightarrow[]{} \RR, \quad x \mapsto  \sum_{i=1}^{2n} (-1)^{i+1} x_i^2.
 \]
 The standard component of $\OGr_{q'}(n,2n)$ is cut out by the Pl\"ucker relations and the additional linear relations
 \[
     \Delta^I(X) = \Delta^{[2n]\setminus I} (X), \quad \text{for all } I \in \binom{[2n]}{n}.
 \]
 The Pl\"ucker non-negative semi-algebraic set in $\OGr_{q'}(n,2n)$ emerged as the geometry behind the ABJM amplitudes in \cite{ABJM2,ABJM1} and, in \cite{IsingModel}, it was connected to the Ising model in statistical mechanics. This semi-algebraic set is rather different from $\OGr^{\geq 0}_{q'}(n,2n)$. It particular, it turns out that some  Pl\"ucker coordinates do not even have a fixed sign in $\OGr^{\geq 0}_{q'}(n,2n)$.

\medskip

More generally, the notion of total positivity introduced by Lusztig for flag varieties only coincides with the positivity of Pl\"ucker coordinates in special cases, see \cite{BBEG24}. In many cases, the geometry that is relevant to applications, for example in physics, is the semi-algebraic set where certain regular functions are positive. In the cases where the notion of positivity we are interested in does not coincide with the notion of total positivity it generally becomes hard to find combinatorial structures which explain the boundary stratification of the semi-algebraic set of interest. An example of such a situation is $\SS_n^{\Pf \geq 0}$ which is interesting to study in its own right. As the structure stops being governed by the combinatorics of the Weyl group, other tools are necessary to study its cell decomposition. Also, as briefly explained in \Cref{sec:5}, this is closely related to considering a \emph{positive Spinor variety}. We leave these questions for future work.

\titleformat{\section}{\centering \fontsize{12}{17} \large \bf \scshape }{\thesection}{0mm}{ \hspace{0.00mm}}

\bibliographystyle{acm}
\bibliography{references}

\end{document}